\documentclass[10pt,aps,prb, preprint,reqno]{amsart}
\usepackage{amsmath}
\usepackage{amssymb}
\usepackage{hyperref}

\newtheorem{theorem}{Theorem}[section]
\newtheorem{remark}{Remark}[section]
\newtheorem{lemma}[theorem]{Lemma}
\newtheorem{proposition}[theorem]{Proposition}

\newtheorem{define}{Definition}[section]

\begin{document}
\title[3D MHD system in scaling-invariant spaces]{On the three-dimensional magnetohydrodynamics system in scaling-invariant spaces}
\author{Kazuo Yamazaki}  
\date{}
\maketitle

\begin{abstract}
We study the criterion for the velocity and magnetic vector fields that solve the three-dimensional magnetohydrodynamics system, given any initial data sufficiently smooth, to experience a finite-time blowup. Following the work of [12] and making use of the structure of the system, we obtain a criterion that is imposed on the magnetic vector field and only one of the three components of the velocity vector field, both in scaling-invariant spaces.

\vspace{5mm}

\end{abstract}
\footnote{2000MSC : 35B65, 35Q35, 35Q86}
\footnote{Department of Mathematics, Washington State University, Pullman, WA 99164-3113, USA}

\section{Introduction and statement of results}

We study the following magnetohydrodynamics system in $\mathbb{R}^{3}$: 
\begin{subequations}
\begin{align}
&\frac{du}{dt} + (u\cdot\nabla) u - (b\cdot\nabla) b + \nabla \pi = \nu \Delta u,\\
&\frac{db}{dt} + (u\cdot\nabla) b - (b\cdot\nabla) u = \eta \Delta b,\\
&\nabla\cdot u = \nabla\cdot b = 0, \hspace{5mm} (u,b)(x,0) = (u_{0}, b_{0})(x),
\end{align}
\end{subequations}
where $u, b: \mathbb{R}^{3} \times \mathbb{R}^{+} \mapsto \mathbb{R}^{3}$ represent the velocity and magnetic vector fields respectively while $\pi: \mathbb{R}^{3} \times \mathbb{R}^{+} \mapsto \mathbb{R}$ the scalar pressure field, $\nu \geq 0$ the viscosity and $\eta \geq 0$ the magnetic diffusivity. Without loss of generality, hereafter we assume $\nu = \eta = 1$ and also write $\partial_{t}$ for $\frac{d}{dt}$ and $\partial_{i}$ for $\frac{d}{dx_{i}}, i = 1, 2, 3, x = (x_{1}, x_{2}, x_{3})$. Let us also set for any three-dimensional vector field $f = (f^{1}, f^{2}, f^{3}), f^{h} \triangleq (f^{1}, f^{2}, 0)$, 
\begin{equation*}
\Omega \triangleq \nabla \times u, \hspace{5mm} j \triangleq \nabla \times b, \hspace{5mm}  \omega \triangleq \Omega \cdot e^{3}, \hspace{5mm}  d \triangleq j \cdot e^{3} \hspace{1mm} \text{ where } e^{3} \triangleq (0, 0, 1).
\end{equation*}
When $b\equiv 0$, (1a)-(1c) recovers the Navier-Stokes equations (NSE), for which the question of whether a smooth local solution can preserve its regularity for all time remains unknown. The analogous problem for the MHD system (1a)-(1c) remains just as difficult, if not more. One of the sources of the difficulty of the global regularity issue of the MHD system (1a)-(1c) may be traced back to the rescaling and its known bounded quantities. It can be shown that if 
$(u,b)(x,t)$ solves the system (1a)-(1c), then so does 
$(u_{\lambda}, b_{\lambda})(x,t) \triangleq \lambda(u,b)(\lambda x, \lambda^{2} t)$ while $\lVert u_{\lambda}(x,t) \rVert_{L^{2}}^{2} + \lVert b_{\lambda}(x,t) \rVert_{L^{2}}^{2} = \lambda^{-1} (\lVert u(x,\lambda^{2}t) \rVert_{L^{2}}^{2} + \lVert b(x,\lambda^{2}t) \rVert_{L^{2}}^{2})$. 

In two-dimensional case, both the NSE and the MHD system, if $\nu, \eta > 0$, admit a unique global smooth solution starting from any data sufficiently smooth (cf. [22, 24]). Due to the difficulty in the three-dimensional case, much effort has been devoted to provide regularity and blow-up criterion some of which we review now. 
 
In [25], the author initiated important research direction of regularity criterion which led to, along with others such as [13], that if a weak solution $u$ of the three-dimensional NSE with $\nu > 0$ satisfies 
\begin{equation*}
u \in L^{r}(0, T; L^{p}(\mathbb{R}^{3})), \hspace{3mm} \frac{3}{p} + \frac{2}{r} \leq 1, p \in [3, \infty],
\end{equation*}
then $u$ is smooth. Among many other results, in [3] it was shown that if $u$ solves the NSE with $\nu > 0$  and 
\begin{equation*}
\nabla u \in L^{r}(0, T; L^{p}(\mathbb{R}^{3})), \hspace{3mm} \frac{3}{p} + \frac{2}{r} \leq 2, \hspace{3mm} 1 < r \leq 3,  
\end{equation*}
then $u$ is a regular solution (cf. also [2, 15]). We emphasize that the norm $\lVert \cdot \rVert_{L_{T}^{r}L_{x}^{p}}$ and $\lVert \cdot \rVert_{L_{T}^{r}\dot{W}_{x}^{1,p}}$ are both invariant under the scalings of the solutions to the NSE and the MHD system precisely when $\frac{3}{p} + \frac{2}{r} = 1, \frac{3}{p} + \frac{2}{r} = 2$ respectively. For the MHD system, e.g. the author in [26] showed that if $
\nabla u, \nabla b \in L^{4}(0, T; L^{2}(\mathbb{R}^{3})),$ then no singularity occurs in $[0,T]$. Moreover, the work in [6] in particular showed that if $[0,T^{\ast})$ is the maximal interval of existence of smooth solution and $T^{\ast} < \infty$, then 
\begin{equation*}
\int_{0}^{T^{\ast}} \lVert \Omega \rVert_{L^{\infty}} + \lVert j\rVert_{L^{\infty}} d\tau = \infty.
\end{equation*}
We note that the authors in [16, 36] independently realized that in particular the criterion for the solution to the MHD system may be reduced to just $u$, dropping condition on $b$ completely (see also [17, 27]). Let us also state results that are directly related to our work. In [14] the authors showed that given $u_{0} \in \dot{H}^{\frac{1}{2}}(\mathbb{R}^{3})$, there exists a maximal interval $[0, T^{\ast})$ on which a unique solution $u$ to the NSE exists. Analogous results with $(u_{0}, b_{0}) \in \dot{H}^{\frac{1}{2}}(\mathbb{R}^{3})$ may be found in [23, 35]. 

We now survey some component reduction results of such conditions. The authors in [21] showed that if $u$ solves the NSE with $\nu > 0$ and 
\begin{align*}
u^{3} \in& L^{r}(0, T; L^{p}(\mathbb{R}^{3})), \hspace{1mm} \frac{3}{p} + \frac{2}{r} \leq \frac{5}{8}, \hspace{5mm} r \in [\frac{54}{23}, \frac{18}{5}],\\
\text{ or } \nabla u^{3} \in& L^{r}(0, T; L^{p}(\mathbb{R}^{3})), \hspace{1mm} \frac{3}{p} + \frac{2}{r} \leq \frac{11}{6}, \hspace{3mm} r \in [\frac{24}{5}, \infty],\nonumber
\end{align*}
then the solution is regular (see also [7 , 37] for similar results on $u^{3}, \nabla u^{3}$). This result was successfully extended to the MHD system as the authors in [20] showed that if $u$ solves (1a)-(1c) with $\nu, \eta > 0$ and 
\begin{equation}
u^{3}, b \in L^{r}(0, T; L^{p}(\mathbb{R}^{3})), \hspace{3mm} \frac{3}{p} + \frac{2}{r} \leq \frac{3}{4} + \frac{1}{2p}, \hspace{3mm} p > \frac{10}{3},
\end{equation}
then the solution pair $(u,b)$ remains smooth for all time. In [32], the author reduced this constraint on $u^{3}, b$ to $u^{3}, b^{1}, b^{2}$ in special cases without worsening the upper bound of $\frac{3}{4} + \frac{1}{2p}$ making use of the structure of (1b). We note however that the upper bound of $\frac{3}{4} + \frac{1}{2p}$ does not allow the norm to be scaling-invariant. For more interesting component reduction results of such criterion, we refer to e.g. [8, 9, 30, 31]. In particular, we point out that the author in [33] obtained the following regularity criterion for the solution to the three-dimensional MHD system:
\begin{equation}
 \begin{cases}u^{3} \in L^{r_{1}}(0, T; L^{p_{1}}(\mathbb{R}^{3})), \hspace{5mm} 
\frac{3}{p_{1}} + \frac{2}{r_{1}} \leq \frac{1}{3} + \frac{1}{2p_{1}}, \hspace{1mm} \frac{15}{2} < p_{1},\\
 d \in L^{r_{2}}(0, T; L^{p_{2}}(\mathbb{R}^{3})), \hspace{7mm} \frac{3}{p_{2}} + \frac{2}{r_{2}} \leq 2, \hspace{12mm} \frac{3}{2} < p_{2}. 
\end{cases}  
\end{equation}
We remark that the upper bound of $p_{2}, r_{2}$ for $d$ allows the scaling-invariant case. 

We now motivate our study. In [12], the authors succeeded in showing that given an initial data for the NSE in $\dot{W}^{1, \frac{3}{2}}(\mathbb{R}^{3})$, if blow-up occurs at $T^{\ast}<\infty$, then $u_{3} \notin L^{p}(0, T^{\ast}; \dot{H}^{\frac{1}{2} + \frac{2}{p}}), p \in (4, 6)$. The purpose of this manuscript is to extend this result to the MHD system (1a)-(1c). We emphasize that because the proof in [12] required taking a curl of the NSE and studying its vorticity formulation carefully, such a generalization to the MHD system is non-trivial. A well-known example for this type of difficulty is that although in two-dimensional case, the author in [34] showed that the solution to the Euler equations admits a unique global smooth solution, it remains unknown if such a result may be extended to the two-dimensional MHD system even with full magnetic diffusion (see [10] and references found therein). Similarly, although in [11], the authors obtained a two-vorticity component regularity criterion for the three-dimensional NSE making use of the Biot-Savart law in a special way, to the best of the author's knowledge, an analogous criterion, e.g. in terms of two-vorticity components, does not exist for the MHD system. The difficulty in extending these results is due to the current density formulation upon taking a curl on (1b) (in particular $M(u,b)$ in (19)). This difficulty appears even in two-dimensional case, e.g. [18] in which the author elaborates why the current density formulation is not as simple as that of vorticity. 

We let $\Omega_{0}(x) = \Omega(x, 0), j_{0}(x) = j(x,0)$ and denote the following norm which is invariant under the scaling of the solution to the MHD system
\begin{equation*}
\lVert f\rVert_{SC_{p, p_{1}, p_{2}}} \triangleq \lVert f\rVert_{\dot{H}^{\frac{1}{2} + \frac{2}{p}}}^{p} + \lVert f\rVert_{L^{p_{1}}}^{r_{1}} + \lVert \nabla f\rVert_{L^{p_{2}}}^{r_{2}}, \hspace{3mm} \frac{3}{p_{1}} + \frac{2}{r_{1}} = 1, \frac{3}{p_{2}} + \frac{2}{r_{2}} = 2.
\end{equation*}
We remark that $\int_{0}^{T} \lVert b_{\lambda}(x,\tau) \rVert_{\dot{H}^{\frac{1}{2} + \frac{2}{p}}}^{p} d\tau = \int_{0}^{\lambda^{2} T}\lVert b(x, \tau)\rVert_{\dot{H}^{\frac{1}{2} + \frac{2}{p}}}^{p} d\tau$. 
Our results read
\begin{theorem}
Suppose $\Omega_{0}, j_{0} \in L^{\frac{3}{2}}(\mathbb{R}^{3})$. Then $\exists !$ solution pair $(u,b)$ to the MHD system (1a)-(1c)  such that $u, b \in C(0, T^{\ast}; \dot{H}^{\frac{1}{2}}(\mathbb{R}^{3})) \cap L_{loc}^{2}(0, T^{\ast}; \dot{H}^{\frac{3}{2}}(\mathbb{R}^{3}))$ and  
\begin{align*}
\sup_{t\in [0, T]}(\lVert \Omega\rVert_{L^{\frac{3}{2}}}^{\frac{3}{2}} + \lVert j&\rVert_{L^{\frac{3}{2}}}^{\frac{3}{2}})(t) + \int_{0}^{T} \lvert \nabla (\Omega + j)\rvert^{2} \lvert \Omega + j\rvert^{-\frac{1}{2}} + \lvert \nabla (\Omega - j)\rvert^{2} \lvert \Omega - j\rvert^{-\frac{1}{2}} d\tau\\
\leq& c(1+ \lVert \Omega_{0} \rVert_{L^{\frac{3}{2}}}^{\frac{3}{2}} + \lVert j_{0} \rVert_{L^{\frac{3}{2}}}^{\frac{3}{2}})\exp\left(\int_{0}^{T} \lVert u\rVert_{\dot{H}^{\frac{3}{2}}}^{2} + \lVert b\rVert_{\dot{H}^{\frac{3}{2}}}^{2} d\tau\right) < \infty
\end{align*}
$\forall \hspace{1mm} T < T^{\ast}$. Moreover, let $p \in (4,6), p_{1} > 9, p_{2} > \frac{9}{2}$. If $T^{\ast} < \infty$, then 
\begin{equation}
\int_{0}^{T^{\ast}}\lVert u^{3} \rVert_{\dot{H}^{\frac{1}{2} + \frac{2}{p}}}^{p} + \lVert b \rVert_{SC_{p, p_{1},p_{2}}} d\tau = \infty.
\end{equation}
\end{theorem}

\begin{remark}
\begin{enumerate}
\item In comparison with (2) and (3), the conditions in (4) is at the scaling-invariant level. 
\item Taking $b \equiv 0$ recovers the result in [12].  
\item The difficulty in the estimate of the third component of the curl formulation in contrast to the case of the NSE is in particular $M(u,b)$ in (19). The difficulty in the estimate of the third component of the velocity equation is the pressure term which now involves the quadratic of the magnetic field (see (41a)). We had to take advantage of the structure of the MHD system, in particular make cancellations such as in (21)-(23) (see also (43a)-(43e)) and $(V_{1} + VI_{1})_{1}$ and $(V_{1} + VI_{1})_{4}$ in (82).  

\end{enumerate}
\end{remark}

In the Preliminaries we set up notations and state key lemmas. Thereafter, we prove the second statement of Theorem 1.1, namely (4). Because local existence theory is classical, we sketch it in the Appendix for completeness. 

\section{Preliminaries}
We write $A \lesssim_{a,b} B$ when there exists a constant $c \geq 0$ of significant dependence only on $a, b$ such that $A \leq c B$, similarly $A \approx_{a,b} B$ if $A = cB$. For simplicity, we denote $\int = \int_{\mathbb{R}^{3}}$ and omit $dx$ when no confusion arises. We also denote  
\begin{equation*}
\nabla_{h}^{\bot} \triangleq (-\partial_{2}, \partial_{1}, 0), \hspace{2mm} \Delta_{h} \triangleq \sum_{i=1}^{2}\partial_{i}^{2},
\end{equation*}
with which we may write down the key identity to be used frequently in our proof, namely $\forall f = (f^{h}, f^{3})$ such that $\nabla \cdot f = 0$, 
\begin{align}
f^{h} = f_{\text{curl}}^{h} + f_{\text{div}}^{h} \hspace{2mm} \text{where} \hspace{2mm}
f_{\text{curl}}^{h} \triangleq \nabla_{h}^{\bot}\Delta_{h}^{-1} (\nabla \times f)\cdot e^{3}, \hspace{2mm} f_{\text{div}}^{h} \triangleq -\nabla_{h}\Delta_{h}^{-1} \partial_{3}f^{3}.
\end{align}
We also denote for any scalar function $f^{i}$, $f^{i}_{\alpha} \triangleq \frac{f^{i}}{\lvert f^{i}\rvert} \lvert f^{i}\rvert^{\alpha}, \alpha \in \mathbb{R}^{+}$. 

We recall the anisotropic Lebesgue spaces reminding ourselves that its order matters, i.e. 
$\lVert \lVert f(\cdot, x_{2}) \rVert_{L^{p}(X_{1}, \mu_{1})} \rVert_{L^{q}(X_{2}, \mu_{2})} \leq \lVert \lVert f(x_{1}, \cdot) \rVert_{L^{q}(X_{2}, \mu_{2})}\rVert_{L^{p}(X_{1}, \mu_{1})}$ for any two measure spaces $(X_{1}, \mu_{1}), (X_{2}, \mu_{2})$ with $1 \leq p \leq q \leq \infty$ (cf. [4]). Let us denote by $\mathcal{S}$ the Schwartz space and $\mathcal{S}'$ its dual. We continue to use the following definitions of anisotropic Sobolev spaces from [12] (see also [19, 29]).   
\begin{define}
For $s, s' \in \mathbb{R}, \dot{H}^{s, s'}$ denotes the space of $f \in \mathcal{S}'$ such that 
\begin{equation*}
\lVert f\rVert_{\dot{H}^{s,s'}}^{2} \triangleq \int_{\mathbb{R}^{3}} \lvert \xi_{h} \rvert^{2s} \lvert \xi_{3} \lvert^{2s'} \rvert\hat{f}(\xi) \rvert^{2} d\xi < \infty, \hspace{5mm} \xi_{h} = (\xi_{1}, \xi_{2}, 0).
\end{equation*}
Moreover, for $\theta \in (0, \frac{1}{2})$, we denote $\mathcal{H}_{\theta} \triangleq \dot{H}^{-\frac{1}{2} + \theta, -\theta}$. 
\end{define}

We recall the Littlewood-Paley decomposition; with $\chi, \phi$ smooth functions such that 
\begin{align*}
&\text{supp } \phi \subset \{\xi \in \mathbb{R} : \frac{3}{4} \leq \lvert \xi\rvert \leq \frac{8}{3}\}, \hspace{10mm} \sum_{j\in\mathbb{Z}} \phi(2^{-j} \xi) = 1,\\
&\text{supp } \chi \subset \{\xi \in \mathbb{R}: \lvert \xi\rvert \leq \frac{4}{3}\}, \hspace{6mm} \chi(\xi) + \sum_{j\geq 0} \phi(2^{-j} \xi) = 1,
\end{align*}
we denote the Littlewood-Paley operators, classical and anisotropic, 
\begin{align*}
&\dot{\Delta}_{j} f \triangleq \mathcal{F}^{-1} (\phi(2^{-j} \lvert \xi\rvert)\hat{f}), \hspace{12mm} \dot{S}_{j}f \triangleq \mathcal{F}^{-1} (\chi(2^{-j} \lvert \xi\rvert)\hat{f}),\\
&\dot{\Delta}_{k}^{h} f \triangleq \mathcal{F}^{-1} (\phi(2^{-k} \lvert \xi_{h} \rvert)\hat{f}), \hspace{10mm} \dot{S}_{k}^{h} f \triangleq \mathcal{F}^{-1} (\chi(2^{-k} \lvert \xi_{h} \rvert)\hat{f}),\\
&\dot{\Delta}_{l}^{v} f \triangleq \mathcal{F}^{-1} (\phi(2^{-l} \lvert \xi_{3} \rvert)\hat{f}), \hspace{11mm} \dot{S}_{l}^{v} f \triangleq \mathcal{F}^{-1} (\chi(2^{-l} \lvert \xi_{3} \rvert)\hat{f}).
\end{align*}
We define $\mathcal{S}_{h}'$ to be the subspace of $\mathcal{S}'$ such that $\lim_{j\to -\infty} \lVert \dot{S}_{j} f\rVert_{L^{\infty}} = 0 \hspace{1mm} \forall f \in \mathcal{S}_{h}'$.  
\begin{define}
For $p, q \in [1,\infty], s \in \mathbb{R}, s < \frac{3}{p} (s = \frac{3}{p} \text{ if } q = 1)$, we define the Besov spaces $\dot{B}_{p,q}^{s}(\mathbb{R}^{3}) \triangleq \{f \in \mathcal{S}_{h}' : \lVert f\rVert_{\dot{B}_{p,q}^{s}} < \infty\}$ where 
\begin{equation*}
\lVert f\rVert_{\dot{B}_{p,q}^{s}} \triangleq \left\lVert (2^{js}\lVert \dot{\Delta}_{j} f\rVert_{L^{p}})_{j} \right\rVert_{l^{q}(\mathbb{Z})}.
\end{equation*}
Moreover, for $p \in (1, \infty)$, we shall use the notations $\mathcal{B}_{p} \triangleq \dot{B}_{\infty, \infty}^{-2 + \frac{2}{p}}$.  

We define the anisotropic Besov spaces $(\dot{B}_{p, q_{1}}^{s_{1}})_{h}(\dot{B}_{p, q_{2}}^{s_{2}})_{v}$ as the space of distributions in $\mathcal{S}_{h}'$ endowed with its norm of 
\begin{equation*}
\lVert f\rVert_{(\dot{B}_{p, q_{1}}^{s_{1}})_{h}(\dot{B}_{p, q_{2}}^{s_{2}})_{v}} \triangleq \left( \sum_{k\in\mathbb{Z}}2^{q_{1} ks_{1}}\left(\sum_{l\in\mathbb{Z}} 2^{q_{2} l s_{2}}\lVert \dot{\Delta}_{k}^{h} \dot{\Delta}_{l}^{v} f\rVert_{L^{p}}^{q_{2}}\right)^{\frac{q_{1}}{q_{2}}} \right)^{\frac{1}{q_{1}}}.
\end{equation*}
\end{define}
It is well-known that $\dot{B}_{2,2}^{s} = \dot{H}^{s}$ (cf. [5]). Moreover, $(\dot{B}_{p, q_{1}}^{s_{1}})_{h}(\dot{B}_{p, q_{2}}^{s_{2}})_{v}\rvert_{p = q_{1} = q_{2} = 2} = \dot{H}^{s_{1}, s_{2}}$. We recall the important Bony's para-product decomposition: 
\begin{equation}
fg = T(f,g) + T(g,f) + R(f,g)
\end{equation}
\begin{equation*}
\text{where } T(f,g) \triangleq \sum_{j\in\mathbb{Z}} \dot{S}_{j-1} f \dot{\Delta}_{j} g, \hspace{3mm} R(f,g) \triangleq \sum_{j\in\mathbb{Z}} \dot{\Delta}_{j} f \tilde{\dot{\Delta}}_{j} g, \hspace{3mm} \text{ with } \tilde{\dot{\Delta}}_{j} \triangleq \sum_{l=j-1}^{j+1}\dot{\Delta}_{l}.
\end{equation*}
We also recall the useful anisotropic Bernstein's inequalities:

\begin{lemma}
Let $\mathcal{B}_{h}$ (respectively $\mathcal{B}_{v}$) a ball in $\mathbb{R}_{h}^{2}$ (resp. $\mathbb{R}_{v}$) and $\mathcal{C}_{h}$ (resp. $\mathcal{C}_{v}$) a ring in $\mathbb{R}_{h}^{2}$ (resp. $\mathbb{R}_{v}$). Moreover, let $1 \leq p_{2} \leq p_{1} \leq \infty, 1 \leq q_{2} \leq q_{1} \leq \infty$. Then 
\begin{align*}
\lVert \nabla_{h}^{ \alpha}f\rVert_{L_{h}^{p_{1}}(L_{v}^{q_{1}})} \lesssim& 2^{k(\lvert \alpha \rvert + 2(\frac{1}{p_{2}} - \frac{1}{p_{1}}))} \lVert f\rVert_{L_{h}^{p_{2}}(L_{v}^{q_{1}})}, \hspace{3mm} &\text{ if supp} \hat{f} \subset 2^{k}\mathcal{B}_{h},\\
\lVert \partial_{3}^{ \beta}f\rVert_{L_{h}^{p_{1}}(L_{v}^{q_{1}})} \lesssim& 2^{l(\lvert \beta \rvert + (\frac{1}{q_{2}} - \frac{1}{q_{1}}))} \lVert f\rVert_{L_{h}^{p_{1}}(L_{v}^{q_{2}})}, \hspace{3mm} &\text{ if supp} \hat{f} \subset 2^{l}\mathcal{B}_{v},\\
\lVert f\rVert_{L_{h}^{p_{1}}(L_{v}^{q_{1}})} \lesssim& 2^{-kN} \sup_{\lvert \alpha \rvert = N} \lVert \nabla_{h}^{\alpha} f\rVert_{L_{h}^{p_{1}}(L_{v}^{q_{1}})}, \hspace{3mm} &\text{ if supp} \hat{f} \subset 2^{k}\mathcal{C}_{h},\\
\lVert f\rVert_{L_{h}^{p_{1}}(L_{v}^{q_{1}})} \lesssim& 2^{-lN}  \lVert \partial_{3}^{N} f\rVert_{L_{h}^{p_{1}}(L_{v}^{q_{1}})}, \hspace{3mm} &\text{ if supp} \hat{f} \subset 2^{l}\mathcal{C}_{v}.
\end{align*}

\end{lemma} 
We also recall in relevance the product law in anisotropic spaces (cf. Lemma 4.5 [12]): for $q \geq 1, p_{1} \geq p_{2} \geq 1, \frac{1}{p_{1}} + \frac{1}{p_{2}} \leq 1, s_{1} < \frac{2}{p_{1}}, s_{2} < \frac{2}{p_{2}}$ (resp. $s_{1} \leq \frac{2}{p_{1}}, s_{2} \leq \frac{2}{p_{2}}$ if $q = 1$), $s_{1} + s_{2} > 0$, $\sigma_{1}<  \frac{1}{p_{1}}, \sigma_{2} < \frac{1}{p_{2}}$ (resp. $\sigma_{1} \leq \frac{1}{p_{1}}, \sigma_{2} \leq \frac{1}{p_{2}}$ if $q = 1$), $\sigma_{1} + \sigma_{2} > 0$,
\begin{equation}
\lVert fg\rVert_{(\dot{B}_{p_{1}, q}^{s_{1} + s_{2} - \frac{2}{p_{2}}})_{h}(\dot{B}_{p_{1}, q}^{\sigma_{1} + \sigma_{2} - \frac{1}{p_{2}}})_{v}} \lesssim \lVert f\rVert_{(\dot{B}_{p_{1}, q}^{s_{1}})_{h}(\dot{B}_{p_{1}, q}^{\sigma_{1}})_{v}}\lVert g\rVert_{(\dot{B}_{p_{2}, q}^{s_{2}})_{h}(\dot{B}_{p_{2}, q}^{\sigma_{2}})_{v}}
\end{equation}
(cf. also [19, 29]).

We now recall several results from [12] on which we will rely. Firstly, the proof of Proposition 4.1 in [12] verifies the following inequality: 

\begin{lemma}
Let $f = (f^{1}, f^{2}, f^{3})$ satisfy $\nabla\cdot f = 0$ and $g = (\nabla \times f) \cdot e_{3}$. Then for $\alpha, \theta \in (0, \frac{1}{2})$, 
\begin{equation*}
\lVert f^{h} \rVert_{(\dot{B}_{2,1}^{1})_{h}(\dot{B}_{2,1}^{\frac{1}{2} - \alpha})_{v}} \lesssim \lVert g_{\frac{3}{4}} \rVert_{L^{2}}^{\frac{1}{3} + \alpha}\lVert \nabla g_{\frac{3}{4}} \rVert_{L^{2}}^{1-\alpha} + \lVert \partial_{3} f^{3} \rVert_{\mathcal{H}_{\theta}}^{\alpha} \lVert \nabla \partial_{3} f^{3} \rVert_{\mathcal{H}_{\theta}}^{1-\alpha}.
\end{equation*}
\end{lemma} 

\begin{lemma}
The following inequalities hold for $f = (f^{1}, f^{2}, f^{3})$: 
\begin{equation}
\text{((2.4) [12])} \hspace{2mm} \lVert \partial_{3} f^{3} \rVert_{\mathcal{H}_{\theta}} \lesssim \lVert f \rVert_{\dot{H}^{\frac{1}{2}}} \hspace{5mm} \text{if }  \nabla\cdot f=  0, 
\end{equation}
\begin{equation}
\text{(Lemma 3.2, (3.8) [12])} \hspace{2mm} \lVert \nabla f^{i}\rVert_{L^{\frac{3}{2}}} \lesssim \lVert \nabla f_{\frac{3}{4}}^{i} \rVert_{L^{2}} \lVert f_{\frac{3}{4}}^{i} \rVert_{L^{2}}^{\frac{1}{3}}, \hspace{1mm} \lVert \nabla f^{i} \rVert_{L^{\frac{9}{5}}} \lesssim \lVert \nabla f_{\frac{3}{4}}^{i} \rVert_{L^{2}}^{\frac{4}{3}},
\end{equation}
\begin{equation}
\text{(Lemma 3.2 [12])} \hspace{2mm}  \lVert f^{i}\rVert_{\dot{H}^{s}} \lesssim \lVert f_{\frac{3}{4}}^{i} \rVert_{L^{2}}^{\frac{5}{6} - s} \lVert \nabla f_{\frac{3}{4}}^{i} \rVert_{L^{2}}^{\frac{1}{2}+ s}, \hspace{3mm} s \in [-\frac{1}{2}, \frac{5}{6}],
\end{equation}
\begin{equation}
\text{(pg. 31 [12])}\hspace{2mm}  \lVert f\rVert_{\mathcal{B}_{p}} \lesssim \lVert f\rVert_{L^{\frac{3}{2}}}^{1-\frac{3}{2p}}\lVert \nabla f\rVert_{L^{\frac{9}{5}}}^{\frac{3}{2p}}, \hspace{3mm} p > \frac{3}{2},
\end{equation} 
\begin{equation}
\text{(pg. 22 [12])} \hspace{2mm} \lVert f\rVert_{\dot{H}^{\theta, \frac{1}{2} - \theta - \frac{1}{p}}} \lesssim \lVert f\rVert_{\mathcal{H}_{\theta}}^{\frac{1}{p}}\lVert \nabla f\rVert_{\mathcal{H}_{\theta}}^{1-\frac{1}{p}}, \hspace{3mm} p > 2.
\end{equation}

\end{lemma}

\begin{lemma}
(Lemma 5.2 [12]) Let $\theta \in (0, \frac{1}{6}), \sigma \in (\frac{3}{4}, 1), s = \frac{3}{2} - \frac{2}{3} \sigma$. Then 
\begin{subequations}
\begin{align}
\lvert \int \partial_{h} \Delta_{h}^{-1} f \partial_{h} g h_{\frac{1}{2}} \rvert \lesssim \lVert f\rVert_{L^{\frac{3}{2}}} \lVert g\rVert_{\dot{H}^{s}} \lVert h_{\frac{3}{4}} \rVert_{\dot{H}^{\sigma}}^{\frac{2}{3}},\\
\lvert \int \partial_{h} \Delta_{h}^{-1} f \partial_{h} g h_{\frac{1}{2}} \rvert \lesssim \lVert f\rVert_{\mathcal{H}_{\theta}} \lVert g\rVert_{\dot{H}^{s}} \lVert h_{\frac{3}{4}} \rVert_{\dot{H}^{\sigma}}^{\frac{2}{3}}.
\end{align}
\end{subequations}
\end{lemma}

\begin{lemma}
(Lemma 6.1 [12]) Let $A$ be a bounded Fourier multiplier, $f, g, h$ any scalar-valued functions. Then for $p, \theta$ such that $0 < \theta < \frac{1}{2} - \frac{1}{p}$, 
\begin{equation}
\lvert (A(D) (fg) \lvert \partial_{3}h)_{\mathcal{H}_{\theta}} \rvert \lesssim \lVert f\rVert_{\dot{H}^{\theta, \frac{1}{2} - \theta - \frac{1}{p}}}\lVert g\rVert_{\dot{H}^{\theta, \frac{1}{2} - \theta - \frac{1}{p}}}\lVert h\rVert_{\dot{H}^{\frac{1}{2} + \frac{2}{p}}}.
\end{equation}
Moreover, for such $\theta$ and $p$, (pg. 22 [12])
\begin{equation}
\lVert f\rVert_{\dot{H}^{\theta, \frac{1}{2} - \theta - \frac{1}{p}}} \lesssim \lVert f_{\frac{3}{4}} \rVert_{L^{2}}^{\frac{p+3}{3p}}\lVert \nabla f_{\frac{3}{4}} \rVert_{L^{2}}^{1-\frac{1}{p}}.
\end{equation}
\end{lemma}

We need the following generalized version of Lemma 6.2 from [CZ14] for our purpose; for completeness we sketch its proof: 

\begin{lemma}
Let $A$ be a bounded Fourier multiplier and $f = (f^{1}, f^{2}, f^{3}), g = (g^{1}, g^{2}, g^{3}), h = (h^{1}, h^{2}, h^{3})$. If $p, \theta$ satisfy $0 < \theta < \frac{1}{2} - \frac{1}{p}, 0 < \frac{2}{p}$, then for $l \in \{1, 2\}$, 
\begin{align*}
&\lvert (A(D) (f^{l} \partial_{l} \partial_{3} g^{3}) \lvert \partial_{3} h^{3})_{\mathcal{H}_{\theta}}\rvert\lesssim \lVert f^{l} \rVert_{(\dot{B}_{2,1}^{1})_{h}(\dot{B}_{2,1}^{\frac{1}{2} - \frac{2}{p}})_{v}} \lVert \nabla \partial_{3} g^{3} \rVert_{\mathcal{H}_{\theta}}\lVert h^{3} \rVert_{\dot{H}^{\frac{1}{2} + \frac{2}{p}}}.
\end{align*}
In particular, if $p > 4$ and $(u,b)$ solves the MHD system (1a)-(1c), then  
\begin{align}
&\lvert (A(D)(u^{l} \partial_{l} \partial_{3} g^{3}) \lvert \partial_{3} h^{3})_{\mathcal{H}_{\theta}}\rvert\\
\lesssim& \left(\lVert \omega_{\frac{3}{4}} \rVert_{L^{2}}^{\frac{1}{3} + \frac{2}{p}} \lVert \nabla \omega_{\frac{3}{4}}\rVert_{L^{2}}^{1-\frac{2}{p}} + \lVert \partial_{3} u^{3} \rVert_{\mathcal{H}_{\theta}}^{\frac{2}{p}}\lVert \nabla \partial_{3} u^{3} \rVert_{\mathcal{H}_{\theta}}^{1-\frac{2}{p}}\right)\lVert \nabla \partial_{3} g^{3} \rVert_{\mathcal{H}_{\theta}}\lVert h^{3} \rVert_{\dot{H}^{\frac{1}{2} + \frac{2}{p}}},\nonumber\\
&\lvert (A(D)(b^{l} \partial_{l} \partial_{3} g^{3}) \lvert \partial_{3} h^{3})_{\mathcal{H}_{\theta}}\rvert\\
\lesssim& \left(\lVert d_{\frac{3}{4}} \rVert_{L^{2}}^{\frac{1}{3} + \frac{2}{p}} \lVert \nabla d_{\frac{3}{4}}\rVert_{L^{2}}^{1-\frac{2}{p}} + \lVert \partial_{3} b^{3} \rVert_{\mathcal{H}_{\theta}}^{\frac{2}{p}}\lVert \nabla \partial_{3} b^{3} \rVert_{\mathcal{H}_{\theta}}^{1-\frac{2}{p}}\right)\lVert \nabla \partial_{3} g^{3} \rVert_{\mathcal{H}_{\theta}}\lVert h^{3} \rVert_{\dot{H}^{\frac{1}{2} + \frac{2}{p}}}.\nonumber
\end{align}

\end{lemma}

\begin{proof}

We estimate  
\begin{align*}
&\lvert (A(D)(f^{l} \partial_{l}\partial_{3} g^{3}) \lvert \partial_{3} h^{3})_{\mathcal{H}_{\theta}}\rvert\\
\lesssim& \lVert f^{l} \partial_{l}\partial_{3} g^{3} \rVert_{\dot{H}^{-1 + 2\theta, \frac{1}{2} - \frac{2}{p} - 2\theta}}\lVert \partial_{3} h^{3} \rVert_{\dot{H}^{0, -\frac{1}{2} + \frac{2}{p}}}\\
\lesssim& \lVert f^{l} \rVert_{(\dot{B}_{2,1}^{1})_{h}(\dot{B}_{2,1}^{\frac{1}{2} - \frac{2}{p}})_{v}}\lVert \partial_{l} \partial_{3} g^{3} \rVert_{(\dot{B}_{2,2}^{-1 + 2\theta})_{h}(\dot{B}_{2,2}^{\frac{1}{2} - 2\theta})_{v}}\lVert h^{3} \rVert_{\dot{H}^{\frac{1}{2} + \frac{2}{p}}}\\
\lesssim& \lVert f^{l} \rVert_{(\dot{B}_{2,1}^{1})_{h}(\dot{B}_{2,1}^{\frac{1}{2} - \frac{2}{p}})_{v}}\lVert \partial_{3} g^{3} \rVert_{(\dot{B}_{2,2}^{ 2\theta})_{h}(\dot{B}_{2,2}^{\frac{1}{2} - 2\theta})_{v}}\lVert h^{3} \rVert_{\dot{H}^{\frac{1}{2} + \frac{2}{p}}}
\end{align*}
where we used the fact that $
\lVert h^{3} \rVert_{\dot{H}^{0, \frac{1}{2} + \frac{2}{p}}}^{2} \lesssim \int \lvert \xi\rvert^{1+\frac{4}{p}}\lvert \hat{h^{3}}(\xi) \rvert^{2}d\xi \approx \lVert h^{3} \rVert_{\dot{H}^{\frac{1}{2} + \frac{2}{p}}}^{2}$. Moreover, we remark that the second inequality actually cannot be an application of the product law in anisotropic spaces (7). Nevertheless, it can be justified by a standard technique of anisotropic space estimate (94). Now since 
\begin{align*}
\lVert \partial_{3} g^{3} \rVert_{(\dot{B}_{2,2}^{ 2\theta})_{h}(\dot{B}_{2,2}^{\frac{1}{2} - 2\theta})_{v}}^{2} 
=& \int \lvert \xi_{h} \rvert^{4\theta} \lvert \xi_{3} \rvert^{1-4\theta}\lvert \widehat{\partial_{3} g^{3}}(\xi) \rvert^{2} d\xi \lesssim \lVert \nabla \partial_{3} g^{3} \rVert_{\mathcal{H}_{\theta}}^{2} 
\end{align*}
which can be verified using that  $\theta < \frac{1}{2}$, we obtain 
\begin{align*}
\lvert (A(D)(f^{l} \partial_{l}\partial_{3} g^{3}) \lvert \partial_{3} h^{3})_{\mathcal{H}_{\theta}}\rvert
\lesssim \lVert f^{l} \rVert_{(\dot{B}_{2,1}^{1})_{h}(\dot{B}_{2,1}^{\frac{1}{2} - \frac{2}{p}})_{v}}\lVert \nabla \partial_{3} g^{3} \rVert_{\mathcal{H}_{\theta}}^{2} \lVert h^{3} \rVert_{\dot{H}^{\frac{1}{2} + \frac{2}{p}}}.
\end{align*}
The particular cases are just consequences Lemma 2.2 with $\alpha = \frac{2}{p}.$ This completes the proof of Lemma 2.6. 
\end{proof}

We end this Preliminaries with the following lemma:

\begin{lemma}
(Lemma 4.3 [12]) Let $s > 0, \alpha \in (0, s)$. Then for $f \in \dot{B}_{p, q}^{s}$,
\begin{equation*}
\lVert f\rVert_{(\dot{B}_{p,q}^{s-\alpha})_{h} (\dot{B}_{p,1}^{\alpha})_{v}}\lesssim \lVert f\rVert_{\dot{B}_{p, q}^{s}}.
\end{equation*}
\end{lemma}

\section{Three Propositions}

\begin{proposition}
Under the hypothesis of Theorem 1.1, for $\theta \in (0, \frac{1}{6})$ the solution to the MHD system (1a)-(1c) satisfies for any $t< T^{\ast}$ 
\begin{align*}
&\frac{2}{3} (\lVert \omega_{\frac{3}{4}} \rVert_{L^{2}}^{2} + \lVert d_{\frac{3}{4}} \rVert_{L^{2}}^{2})(t) + \frac{5}{9} \int_{0}^{t} \lVert \nabla \omega_{\frac{3}{4}} \rVert_{L^{2}}^{2} + \lVert \nabla d_{\frac{3}{4}} \rVert_{L^{2}}^{2} d\tau \\
\lesssim& e^{c\int_{0}^{t}\lVert u^{3}\rVert_{\dot{H}^{\frac{1}{2} + \frac{2}{p}}}^{p} + \lVert b\rVert_{SC_{p, p_{1}, p_{2}}} d\tau}\\
&\times \left(\frac{2}{3} (\lVert \omega_{\frac{3}{4}} \rVert_{L^{2}}^{2} + \lVert d_{\frac{3}{4}} \rVert_{L^{2}}^{2})(0) + \left(\int_{0}^{t} \lVert \nabla \partial_{3} u^{3} \rVert_{\mathcal{H}_{\theta}}^{2} + \lVert \nabla \partial_{3} b^{3} \rVert_{\mathcal{H}_{\theta}}^{2}d\tau \right)^{\frac{3}{4}}\right).
\end{align*}
\end{proposition}

\begin{remark}
The fact that the MHD system (1a)-(1c) forces a worse bound in terms of $\int_{0}^{t} \lVert \nabla \partial_{3} u^{3} \rVert_{\mathcal{H}_{\theta}}^{2} + \lVert \nabla \partial_{3} b^{3} \rVert_{\mathcal{H}_{\theta}}^{2}d\tau$ rather than $\int_{0}^{t} \lVert \partial_{3}^{2} u^{3} \rVert_{\mathcal{H}_{\theta}}^{2} d\tau$ in [12] is in particular due to the matrix $M(u,b)$ in (19). 
\end{remark}

\begin{proof}

We take a curl on (1a), (1b) to obtain 
\begin{subequations}
\begin{align}
&\partial_{t} \Omega - \Delta \Omega + (u\cdot\nabla) \Omega - (\Omega\cdot\nabla) u - (b\cdot\nabla) j + (j\cdot\nabla) b = 0\\
&\partial_{t} j - \Delta j + (u\cdot\nabla) j - (j\cdot\nabla) u - (b\cdot\nabla) \Omega + (\Omega\cdot\nabla) b = 2 M(u,b)
\end{align}
\end{subequations}
where
\begin{equation}
M(u,b) \triangleq 
\begin{pmatrix}
\partial_{2} b \cdot\partial_{3} u - \partial_{3} b \cdot\partial_{2} u\\
\partial_{3} b \cdot\partial_{1} u - \partial_{1} b\cdot\partial_{3} u\\
\partial_{1} b\cdot\partial_{2} u - \partial_{2} b \cdot\partial_{1} u 
\end{pmatrix}.
\end{equation}
In particular, the third components of this system reads
\begin{subequations}
\begin{align}
&\partial_{t} \omega + (u\cdot\nabla) \omega - (b\cdot\nabla) d - \Delta \omega = (\Omega \cdot\nabla) u^{3} - (j\cdot\nabla) b^{3},\\
&\partial_{t} d + (u\cdot\nabla) d - (b\cdot\nabla) \omega - \Delta d = (j\cdot\nabla) u^{3} - (\Omega\cdot\nabla) b^{3} + 2[\partial_{1} b \cdot\partial_{2} u - \partial_{2} b \cdot\partial_{1} u].
\end{align}
\end{subequations}
We make a few important cancellations: 
\begin{align}
(\Omega \cdot\nabla) u^{3} = \partial_{3} u^{3} \omega + \partial_{2} u^{3} \partial_{3} u^{1} - \partial_{1} u^{3} \partial_{3} u^{2},
\end{align}
\begin{align}
(j\cdot\nabla) b^{3}
= \partial_{3}b^{3} d + \partial_{2}b^{3} \partial_{3}b^{1} - \partial_{1}b^{3}\partial_{3} b^{2}.
\end{align}
We make cancellations within $2[\partial_{1} b \cdot\partial_{2} u - \partial_{2} b \cdot\partial_{1} u]$ as well: 
\begin{align}
&(j\cdot\nabla) u^{3} - (\Omega \cdot\nabla)b^{3} + 2[\partial_{1} b \cdot\partial_{2} u - \partial_{2} b \cdot\partial_{1} u]\\
=& \partial_{3}u^{3} d - \partial_{3} b^{3} \omega - \partial_{3} b^{2} \partial_{1} u^{3} + \partial_{3} b^{1} \partial_{2} u^{3} + \partial_{3} u^{2} \partial_{1} b^{3} - \partial_{3} u^{1} \partial_{2} b^{3}\nonumber \\
&+ 2[\partial_{1} b^{1} \partial_{2} u^{1} - \partial_{2} b^{1} \partial_{1} u^{1} + \partial_{1} b^{2} \partial_{2} u^{2} - \partial_{2} b^{2} \partial_{1} u^{2} ].\nonumber
\end{align}
Therefore, we have from (20a)-(20b), (21)-(23), 
\begin{subequations}
\begin{align}
&\partial_{t} \omega + (u\cdot\nabla) \omega - (b\cdot\nabla) d - \Delta \omega = F_{1},\\
&\partial_{t} d + (u\cdot\nabla) d - (b\cdot\nabla) \omega - \Delta d = F_{2},\\
& F_{1} \triangleq (\partial_{3} u^{3} \omega + \partial_{2} u^{3} \partial_{3} u^{1} - \partial_{1} u^{3} \partial_{3} u^{2}) - (\partial_{3}b^{3} d + \partial_{2} b^{3} \partial_{3} b^{1} - \partial_{1} b^{3} \partial_{3} b^{2}),\\
& F_{2} \triangleq \partial_{3}u^{3} d - \partial_{3} b^{3} \omega - \partial_{3} b^{2} \partial_{1} u^{3} + \partial_{3} b^{1} \partial_{2} u^{3} + \partial_{3} u^{2} \partial_{1} b^{3} - \partial_{3} u^{1} \partial_{2} b^{3} \\
& \hspace{5mm} + 2[\partial_{1}b^{h} \cdot \partial_{2}u^{h} - \partial_{2}b^{h} \cdot \partial_{1}u^{h} ].\nonumber
\end{align}
\end{subequations} 
Taking $L^{\frac{3}{2}}$-norm estimate, using divergence-free conditions and that $\lvert \nabla f_{\frac{3}{4}}^{i} \rvert = \frac{3}{4} \lvert \nabla f^{i} \rvert \lvert f^{i} \rvert^{-\frac{1}{4}}$, integrating in time we obtain with $\omega(0) = \omega_{0}, d(0) = d_{0}$, 
\begin{align}
&\frac{2}{3}(\lVert \omega \rVert_{L^{\frac{3}{2}}}^{\frac{3}{2}} + \lVert d \rVert_{L^{\frac{3}{2}}}^{\frac{3}{2}})(t) + \frac{8}{9}\int_{0}^{t} \lVert \nabla \omega_{\frac{3}{4}} \rVert_{L^{2}}^{2} + \lVert \nabla d_{\frac{3}{4}} \rVert_{L^{2}}^{2} d\tau\\
=& \frac{2}{3} \left(\lVert \omega_{0} \rVert_{L^{\frac{3}{2}}}^{\frac{3}{2}} + \lVert d_{0} \rVert_{L^{\frac{3}{2}}}^{\frac{3}{2}}\right)+ \int_{0}^{t} \int (b\cdot\nabla) d \omega_{\frac{1}{2}} + (b\cdot\nabla) \omega d_{\frac{1}{2}} + F_{1} \omega_{\frac{1}{2}} + F_{2} d_{\frac{1}{2}}  d\tau.\nonumber
\end{align}
We first estimate 
\begin{align}
&\int (b\cdot\nabla) d \omega_{\frac{1}{2}} + (b\cdot\nabla) \omega d_{\frac{1}{2}}\\
\lesssim& \int \lvert b\rvert \lvert \nabla d_{\frac{3}{4}}\rvert \lvert d\rvert^{\frac{1}{4}} \lvert \omega\rvert^{\frac{1}{2}} + \lvert b\rvert \lvert \nabla \omega_{\frac{3}{4}} \rvert \lvert \omega \rvert^{\frac{1}{4}} \lvert d\rvert^{\frac{1}{2}}\nonumber\\
\lesssim& \lVert b\rVert_{L^{p_{1}}} \left(\lVert \nabla d_{\frac{3}{4}} \rVert_{L^{2}} \lVert d_{\frac{3}{4}} \rVert_{L^{\frac{2p_{1}}{p_{1} - 6}}}^{\frac{1}{3}}\lVert \omega_{\frac{3}{4}} \rVert_{L^{2}}^{\frac{2}{3}} +  \lVert \nabla \omega_{\frac{3}{4}} \rVert_{L^{2}} \lVert \omega_{\frac{3}{4}} \rVert_{L^{\frac{2p_{1}}{p_{1} - 6}}}^{\frac{1}{3}}\lVert d_{\frac{3}{4}} \rVert_{L^{2}}^{\frac{2}{3}}\right)\nonumber\\
\lesssim& \lVert b\rVert_{L^{p_{1}}} \left(\lVert \nabla d_{\frac{3}{4}} \rVert_{L^{2}}^{\frac{p_{1} + 3}{p_{1}}} \lVert d_{\frac{3}{4}} \rVert_{L^{2}}^{\frac{1}{3}(\frac{p_{1} - 9}{p_{1}})}\lVert \omega_{\frac{3}{4}} \rVert_{L^{2}}^{\frac{2}{3}} + 
 \lVert \nabla \omega_{\frac{3}{4}} \rVert_{L^{2}}^{\frac{p_{1} + 3}{p_{1}}} \lVert \omega_{\frac{3}{4}} \rVert_{L^{2}}^{\frac{1}{3}(\frac{p_{1} - 9}{p_{1}})}\lVert d_{\frac{3}{4}} \rVert_{L^{2}}^{\frac{2}{3}}\right)\nonumber\\
\leq& \frac{1}{18} \left(\lVert \nabla \omega_{\frac{3}{4}} \rVert_{L^{2}}^{2} + \lVert \nabla d_{\frac{3}{4}} \rVert_{L^{2}}^{2}\right) + c \lVert b\rVert_{L^{p_{1}}}^{\frac{2p_{1}}{p_{1} - 3}}(\lVert \omega_{\frac{3}{4}} \rVert_{L^{2}}^{2} + \lVert d_{\frac{3}{4}} \rVert_{L^{2}}^{2} )\nonumber
\end{align}
where we used H$\ddot{o}$lder's, Gagliardo-Nirenberg and Young's inequalities. 

Next, we rearrange terms carefully and estimate differently as follows: 
\begin{align}
\int_{0}^{t} \int F_{1}\omega_{\frac{1}{2}} + F_{2} d_{\frac{1}{2}} d\tau \triangleq \sum_{i=1}^{5}I_{i}
\end{align}
where
\begin{align*}
&I_{1} = \int_{0}^{t} \int \partial_{3} u^{3} \omega\omega_{\frac{1}{2}} + \partial_{3}u^{3} dd_{\frac{1}{2}}d\tau, \hspace{5mm} I_{2} = - \int_{0}^{t} \int \partial_{3}b^{3} d \omega_{\frac{1}{2}} + \partial_{3} b^{3} \omega d_{\frac{1}{2}}d\tau,\nonumber\\
&I_{3} = \int_{0}^{t} \int [\partial_{2} u^{3} \partial_{3} u^{1} - \partial_{1} u^{3} \partial_{3} u^{2} - \partial_{2} b^{3} \partial_{3} b^{1} + \partial_{1} b^{3} \partial_{3}b^{2}]\omega_{\frac{1}{2}}d\tau,\nonumber\\
&I_{4} = \int_{0}^{t} \int[- \partial_{3} b^{2} \partial_{1} u^{3} + \partial_{3} b^{1} \partial_{2} u^{3} + \partial_{3} u^{2} \partial_{1} b^{3} - \partial_{3} u^{1} \partial_{2} b^{3}]d_{\frac{1}{2}}d\tau,\\ 
&I_{5} = 2\int_{0}^{t} \int(\partial_{1}b^{h} \cdot \partial_{2}u^{h} - \partial_{2}b^{h} \cdot \partial_{1}u^{h} ) d_{\frac{1}{2}} 
d\tau.
\end{align*}
Firstly, after integrating by parts we estimate 
\begin{align}
I_{1} \lesssim& \int_{0}^{t} \int \lvert u^{3} \rvert (\lvert \omega \rvert^{\frac{1}{2}} \lvert \partial_{3} \omega \rvert + \lvert d\rvert^{\frac{1}{2}} \lvert \partial_{3} d\rvert ) d\tau\\
\lesssim& \int_{0}^{t} \lVert u^{3} \rVert_{L^{\frac{3p}{p-2}}}\left(\left(\int (\lvert \omega \rvert^{\frac{1}{2}} )^{\frac{3p}{2}} \right)^{\frac{2}{3p}} \lVert \partial_{3} \omega \rVert_{L^{\frac{3}{2}}} +
\left(\int (\lvert d\rvert^{\frac{1}{2}} )^{\frac{3p}{2}} \right)^{\frac{2}{3p}}\lVert \partial_{3} d\rVert_{L^{\frac{3}{2}}}\right)d\tau\nonumber\\
\lesssim& \int_{0}^{t} \lVert u^{3} \rVert_{\dot{H}^{\frac{1}{2} + \frac{2}{p}}}(\lVert \omega_{\frac{3}{4}} \rVert_{L^{p}}^{\frac{2}{3}} \lVert \nabla \omega \rVert_{L^{\frac{3}{2}}} + \lVert d_{\frac{3}{4}} \rVert_{L^{p}}^{\frac{2}{3}}\lVert \nabla d\rVert_{L^{\frac{3}{2}}})d\tau\nonumber\\
\lesssim& \int_{0}^{t} \lVert u^{3} \rVert_{\dot{H}^{\frac{1}{2} + \frac{2}{p}}} (\lVert \nabla \omega_{\frac{3}{4}} \rVert_{L^{2}}^{2(\frac{p-1}{p})} + \lVert \nabla d_{\frac{3}{4}} \rVert_{L^{2}}^{2(\frac{p-1}{p})})(\lVert \omega_{\frac{3}{4}} \rVert_{L^{2}}^{\frac{2}{p}} + \lVert d_{\frac{3}{4}} \rVert_{L^{2}}^{\frac{2}{p}})d\tau\nonumber\\
\leq& \frac{1}{18} \int_{0}^{t} \lVert \nabla \omega_{\frac{3}{4}} \rVert_{L^{2}}^{2} + \lVert \nabla d_{\frac{3}{4}} \rVert_{L^{2}}^{2} d\tau + c \int_{0}^{t} \lVert u^{3} \rVert_{\dot{H}^{\frac{1}{2} + \frac{2}{p}}}^{p}(\lVert \omega_{\frac{3}{4}} \rVert_{L^{2}}^{2} + \lVert d_{\frac{3}{4}} \rVert_{L^{2}}^{2})d\tau\nonumber 
\end{align}
by H$\ddot{o}$lder's inequalities, Sobolev embedding of $\dot{H}^{\frac{1}{2} + \frac{2}{p}}(\mathbb{R}^{3}) \hookrightarrow L^{\frac{3p}{p-2}}(\mathbb{R}^{3})$, Gagliardo-Nirenberg inequalities, (9) and Young's inequalities.  

Next, we estimate 
\begin{align}
I_{2}\lesssim& \int_{0}^{t} \lVert \nabla b\rVert_{L^{p_{2}}}(\lVert d\rVert_{L^{\frac{3}{2}}} \left(\int \lvert \omega \rvert^{\frac{1}{2}(\frac{3p_{2}}{p_{2} - 3})}\right)^{\frac{p_{2} - 3}{3p_{2}}} + 
\lVert \omega\rVert_{L^{\frac{3}{2}}} \left(\int \lvert d \rvert^{\frac{1}{2}(\frac{3p_{2}}{p_{2} - 3})}\right)^{\frac{p_{2} - 3}{3p_{2}}})d\tau\\
\lesssim&\int_{0}^{t} \lVert \nabla b\rVert_{L^{p_{2}}} (
\lVert d_{\frac{3}{4}}\rVert_{L^{2}}^{\frac{4}{3}} \lVert \omega_{\frac{3}{4}} \rVert_{L^{2}}^{\frac{2p_{2}-9}{3p_{2}}}\lVert \nabla \omega_{\frac{3}{4}} \rVert_{L^{2}}^{\frac{3}{p_{2}}} + 
\lVert \omega_{\frac{3}{4}}\rVert_{L^{\frac{3}{2}}}^{\frac{4}{3}} \lVert d_{\frac{3}{4}} \rVert_{L^{2}}^{\frac{2p_{2}-9}{3p_{2}}}\lVert \nabla d_{\frac{3}{4}} \rVert_{L^{2}}^{\frac{3}{p_{2}}} )d\tau\nonumber\\
\leq& \frac{1}{18} \int_{0}^{t} 
\lVert \nabla \omega_{\frac{3}{4}} \rVert_{L^{2}}^{2} + \lVert \nabla d_{\frac{3}{4}} \rVert_{L^{2}}^{2} d\tau + c \int_{0}^{t} \lVert \nabla b\rVert_{L^{p_{2}}}^{\frac{2p_{2}}{2p_{2} - 3}}\left(\lVert d_{\frac{3}{4}} \rVert_{L^{2}}^{2} + \lVert \omega_{\frac{3}{4}} \rVert_{L^{2}}^{2}\right)d\tau\nonumber
\end{align}
by H$\ddot{o}$lder's, Gagliardo-Nirenberg and Young's inequalities. 

Next, we first write by (5) 
\begin{align}
I_{3} =& \int_{0}^{t} \int [\partial_{2} u^{3} \partial_{2} \Delta_{h}^{-1} (-\partial_{3} \omega) + \partial_{1} u^{3} \partial_{1} \Delta_{h}^{-1} (-\partial_{3} \omega)\nonumber\\
& \hspace{35mm} + \partial_{2} b^{3} \partial_{2} \Delta_{h}^{-1} \partial_{3} d + \partial_{1} b^{3} \partial_{1} \Delta_{h}^{-1} \partial_{3} d] \omega_{\frac{1}{2}}\nonumber\\
&+ [\partial_{2} u^{3}\partial_{1} \Delta_{h}^{-1} (-\partial_{3}^{2} u^{3}) + \partial_{1} u^{3} \partial_{2} \Delta_{h}^{-1} (\partial_{3}^{2} u^{3})\nonumber\\
& \hspace{35mm} + \partial_{2} b^{3} \partial_{1} \Delta_{h}^{-1} (\partial_{3}^{2} b^{3}) - \partial_{1} b^{3} \partial_{2} \Delta_{h}^{-1} (\partial_{3}^{2} b^{3})] \omega_{\frac{1}{2}} d\tau.
\end{align}
We bound the first and second terms of (30) as follows:
\begin{align}
&\int_{0}^{t} \int [\partial_{2} u^{3} \partial_{2} \Delta_{h}^{-1} (-\partial_{3} \omega) + \partial_{1} u^{3} \partial_{1} \Delta_{h}^{-1} (-\partial_{3} \omega)] \omega_{\frac{1}{2}}d\tau\\
\lesssim& \int_{0}^{t} \lVert \partial_{3} \omega \rVert_{L^{\frac{3}{2}}} \lVert u^{3} \rVert_{\dot{H}^{\frac{3}{2} - \frac{2}{3}\sigma}} \lVert \omega_{\frac{3}{4}} \rVert_{\dot{H}^{\sigma}}^{\frac{2}{3}}d\tau\nonumber\\
\lesssim& \int_{0}^{t} \lVert u^{3} \rVert_{\dot{H}^{\frac{1}{2} + \frac{2}{p}}} \lVert \omega_{\frac{3}{4}} \rVert_{L^{2}}^{\frac{2}{p}} \lVert \nabla \omega_{\frac{3}{4}} \rVert_{L^{2}}^{2(1-\frac{1}{p})}d\tau\nonumber\\
\leq& \frac{1}{36} \int_{0}^{t} \lVert \nabla \omega_{\frac{3}{4}} \rVert_{L^{2}}^{2} d\tau + c\int_{0}^{t} \lVert u^{3} \rVert_{\dot{H}^{\frac{1}{2} + \frac{2}{p}}}^{p}\lVert \omega_{\frac{3}{4}} \rVert_{L^{2}}^{2} d\tau \nonumber
\end{align}
by (13a) with $f = \partial_{3} \omega, g = u^{3}, h = \omega$, $\sigma = 3(\frac{1}{2} - \frac{1}{p})$, (9), Gagliardo-Nirenberg and Young's inequalities. Similarly we bound third and fourth terms of (30) by 
\begin{align}
&\int_{0}^{t} \int [\partial_{2} b^{3} \partial_{2} \Delta_{h}^{-1} (\partial_{3} d) + \partial_{1} b^{3} \partial_{1} \Delta_{h}^{-1} \partial_{3} d] \omega^{\frac{1}{2}}d\tau\\
\lesssim& \int_{0}^{t} \lVert \partial_{3} d\rVert_{L^{\frac{3}{2}}} \lVert b^{3} \rVert_{\dot{H}^{\frac{3}{2} - \frac{2}{3} \sigma}}\lVert \omega_{\frac{3}{4}} \rVert_{\dot{H}^{\sigma}}^{\frac{2}{3}}d\tau\nonumber\\
\lesssim& \int_{0}^{t} \lVert b^{3} \rVert_{\dot{H}^{\frac{1}{2} + \frac{2}{p}}} \lVert \nabla d_{\frac{3}{4}} \rVert_{L^{2}}\lVert  d_{\frac{3}{4}} \rVert_{L^{2}}^{\frac{1}{3}} \lVert \omega_{\frac{3}{4}} \rVert_{L^{2}}^{\frac{2}{p} - \frac{1}{3}}\lVert \nabla \omega_{\frac{3}{4}} \rVert_{L^{2}}^{1-\frac{2}{p}}  d\tau\nonumber\\
\leq& \frac{1}{72} \int_{0}^{t} \lVert \nabla \omega_{\frac{3}{4}} \rVert_{L^{2}}^{2} + \lVert \nabla d_{\frac{3}{4}} \rVert_{L^{2}}^{2} d\tau + c\int_{0}^{t} \lVert b^{3} \rVert_{\dot{H}^{\frac{1}{2} + \frac{2}{p}}}^{p} (\lVert \omega_{\frac{3}{4}} \rVert_{L^{2}}^{2} + \lVert d_{\frac{3}{4}} \rVert_{L^{2}}^{2})d\tau\nonumber
\end{align}
by (13a) with $f = \partial_{3} d, g = b^{3}, h = \omega$, (9), $\sigma = 3 (\frac{1}{2} - \frac{1}{p})$, Gagliardo-Nirenberg and Young's inequalities. Next, we bound the fifth, sixth, seventh and eighth terms of (30) by 
\begin{align}
&\int_{0}^{t} \int 
[\partial_{2} u^{3}\partial_{1} \Delta_{h}^{-1} (-\partial_{3}^{2} u^{3}) + \partial_{1} u^{3} \partial_{2} \Delta_{h}^{-1} (\partial_{3}^{2} u^{3})\\
& \hspace{35mm} + \partial_{2} b^{3} \partial_{1} \Delta_{h}^{-1} (\partial_{3}^{2} b^{3}) - \partial_{1} b^{3} \partial_{2} \Delta_{h}^{-1} (\partial_{3}^{2} b^{3})] \omega_{\frac{1}{2}} d\tau\nonumber\\
\lesssim& \int_{0}^{t} (\lVert \partial_{3}^{2} u^{3} \rVert_{\mathcal{H}_{\theta}} + \lVert \partial_{3}^{2} b^{3} \rVert_{\mathcal{H}_{\theta}})(\lVert u^{3} \rVert_{\dot{H}^{\frac{3}{2} - \frac{2}{3} \sigma}} + \lVert b^{3} \rVert_{\dot{H}^{\frac{3}{2} - \frac{2}{3}\sigma}})\lVert \omega_{\frac{3}{4}} \rVert_{\dot{H}^{\sigma}}^{\frac{2}{3}}d\tau\nonumber\\
\lesssim& \int_{0}^{t} (\lVert \partial_{3}^{2} u^{3} \rVert_{\mathcal{H}_{\theta}} + \lVert \partial_{3}^{2} b^{3} \rVert_{\mathcal{H}_{\theta}})(\lVert u^{3} \rVert_{\dot{H}^{\frac{1}{2} + \frac{2}{p}}} + \lVert b^{3} \rVert_{\dot{H}^{\frac{1}{2} + \frac{2}{p}}})\lVert \omega_{\frac{3}{4}} \rVert_{L^{2}}^{\frac{2}{3}(1-\sigma)}\lVert \nabla \omega_{\frac{3}{4}} \rVert_{L^{2}}^{\frac{2\sigma}{3}}d\tau\nonumber\\
\lesssim& \left(\int_{0}^{t} \lVert \partial_{3}^{2} u^{3} \rVert_{\mathcal{H}_{\theta}}^{2} + \lVert \partial_{3}^{2} b^{3} \rVert_{\mathcal{H}_{\theta}}^{2} d\tau \right)^{\frac{1}{2}} \left(\int_{0}^{t} \lVert u^{3} \rVert_{\dot{H}^{\frac{1}{2} + \frac{2}{p}}}^{p} + \lVert b^{3} \rVert_{\dot{H}^{\frac{1}{2} + \frac{2}{p}}}^{p} d\tau \right)^{\frac{1}{6}}\nonumber\\
&\times \left(\int_{0}^{t} ( \lVert u^{3} \rVert_{\dot{H}^{\frac{1}{2} + \frac{2}{p}}}^{p} + \lVert b^{3} \rVert_{\dot{H}^{\frac{1}{2} + \frac{2}{p}}}^{p}) \lVert \omega_{\frac{3}{4}} \rVert_{L^{2}}^{2} d\tau \right)^{\frac{1}{p} - \frac{1}{6}}\left(\int_{0}^{t} \lVert \nabla \omega_{\frac{3}{4}} \rVert_{L^{2}}^{2} d\tau\right)^{\frac{1}{2} - \frac{1}{p}}\nonumber\\
\leq& \frac{1}{72} \int_{0}^{t} \lVert \nabla \omega_{\frac{3}{4}} \rVert_{L^{2}}^{2} d\tau + c \int_{0}^{t} (\lVert u^{3} \rVert_{\dot{H}^{\frac{1}{2} + \frac{2}{p}}}^{p} + \lVert b^{3} \rVert_{\dot{H}^{\frac{1}{2} + \frac{2}{p}}}^{p}) \lVert \omega_{\frac{3}{4}} \rVert_{L^{2}}^{2} d\tau\nonumber \\
& \hspace{5mm} + c \left(\int_{0}^{t} \lVert \partial_{3}^{2} u^{3} \rVert_{\mathcal{H}_{\theta}}^{2} + \lVert \partial_{3}^{2} b^{3} \rVert_{\mathcal{H}_{\theta}}^{2} d\tau \right)^{\frac{3}{4}} \left(\int_{0}^{t} \lVert u^{3} \rVert_{\dot{H}^{\frac{1}{2} + \frac{2}{p}}}^{p} +  \lVert b^{3} \rVert_{\dot{H}^{\frac{1}{2} + \frac{2}{p}}}^{p}d\tau\right)^{\frac{1}{4}}\nonumber
\end{align}
by (13b) with $(f,g, h) = (\partial_{3}^{2}u^{3}, u^{3}, \omega), (\partial_{3}^{2}b^{3}, b^{3}, \omega), \sigma = 3(\frac{1}{2} - \frac{1}{p})$, Gagliardo-Nirenberg, H$\ddot{o}$lder's and Young's inequalities. Therefore, in sum of (31)-(33) in (30), we have 
\begin{align}
I_{3} \leq& \frac{1}{18} \int_{0}^{t} \lVert \nabla \omega_{\frac{3}{4}} \rVert_{L^{2}}^{2} + \lVert \nabla d_{\frac{3}{4}} \rVert_{L^{2}}^{2} d\tau\\
&+ c\int_{0}^{t} (\lVert u^{3} \rVert_{\dot{H}^{\frac{1}{2} + \frac{2}{p}}}^{p} + \lVert b^{3} \rVert_{\dot{H}^{\frac{1}{2} + \frac{2}{p}}}^{p}) (\lVert \omega_{\frac{3}{4}} \rVert_{L^{2}}^{2} + \lVert d_{\frac{3}{4}} \rVert_{L^{2}}^{2}) d\tau\nonumber \\
&  + c \left(\int_{0}^{t} \lVert \partial_{3}^{2} u^{3} \rVert_{\mathcal{H}_{\theta}}^{2} + \lVert \partial_{3}^{2} b^{3} \rVert_{\mathcal{H}_{\theta}}^{2} d\tau \right)^{\frac{3}{4}} \left(\int_{0}^{t} \lVert u^{3} \rVert_{\dot{H}^{\frac{1}{2} + \frac{2}{p}}}^{p} +  \lVert b^{3} \rVert_{\dot{H}^{\frac{1}{2} + \frac{2}{p}}}^{p}d\tau\right)^{\frac{1}{4}}.\nonumber
\end{align}
Similarly, we can rewrite by (5) and then estimate 
\begin{align}
I_{4} 
=& \int_{0}^{t} \int [-\partial_{1} u^{3} \partial_{1} \Delta_{h}^{-1} \partial_{3} d - \partial_{2} u^{3} \partial_{2} \Delta_{h}^{-1} \partial_{3} d\\
& \hspace{35mm} + \partial_{1} b^{3} \partial_{1} \Delta_{h}^{-1} \partial_{3} \omega + \partial_{2} b^{3} \partial_{2} \Delta_{h}^{-1} \partial_{3} \omega] d_{\frac{1}{2}} d\tau\nonumber\\
&+ [\partial_{1} u^{3} \partial_{2} \Delta_{h}^{-1} \partial_{3}^{2} b^{3} - \partial_{2} u^{3} \partial_{1} \Delta_{h}^{-1} \partial_{3}^{2} b^{3}\nonumber\\
& \hspace{35mm} - \partial_{1} b^{3} \partial_{2} \Delta_{h}^{-1} \partial_{3}^{2} u^{3} + \partial_{2}b^{3} \partial_{1} \Delta_{h}^{-1} \partial_{3}^{2} u^{3} ] d_{\frac{1}{2}} d\tau\nonumber\\
\lesssim& \int_{0}^{t} (\lVert \partial_{3} d\rVert_{L^{\frac{3}{2}}} + \lVert \partial_{3} \omega \rVert_{L^{\frac{3}{2}}} )(\lVert u^{3} \rVert_{\dot{H}^{\frac{3}{2} - \frac{2}{3} \sigma}} + \lVert b^{3} \rVert_{\dot{H}^{\frac{3}{2} - \frac{2}{3} \sigma}})\lVert d_{\frac{3}{4}} \rVert_{\dot{H}^{\sigma}}^{\frac{2}{3}}\nonumber\\
&+ (\lVert \partial_{3}^{2} b^{3} \rVert_{\mathcal{H}_{\theta}} + \lVert \partial_{3}^{2} u^{3} \rVert_{\mathcal{H}_{\theta}})(\lVert u^{3} \rVert_{\dot{H}^{\frac{3}{2} - \frac{2}{3} \sigma}} + \lVert b^{3} \rVert_{\dot{H}^{\frac{3}{2} - \frac{2}{3} \sigma}})\lVert d_{\frac{3}{4}} \rVert_{\dot{H}^{\sigma}}^{\frac{2}{3}} d\tau\nonumber\\
\lesssim& \int_{0}^{t} \left( \lVert \nabla d_{\frac{3}{4}} \rVert_{L^{2}}^{2} + \lVert \nabla \omega_{\frac{3}{4}} \rVert_{L^{2}}^{2} \right)^{1-\frac{1}{p}} \left( \lVert \omega_{\frac{3}{4}} \rVert_{L^{2}} + \lVert d_{\frac{3}{4}} \rVert_{L^{2}}\right)^{\frac{1}{3}}\nonumber\\
&\times \left( \lVert u^{3} \rVert_{\dot{H}^{\frac{1}{2} + \frac{2}{p}}} + \lVert b^{3} \rVert_{\dot{H}^{\frac{1}{2} + \frac{2}{p}}}\right) \lVert d_{\frac{3}{4}} \rVert_{L^{2}}^{-\frac{1}{3} + \frac{2}{p}}\nonumber\\
&+ (\lVert \partial_{3}^{2} u^{3} \rVert_{\mathcal{H}_{\theta}} + \lVert \partial_{3}^{2} b^{3} \rVert_{\mathcal{H}_{\theta}})(\lVert u^{3} \rVert_{\dot{H}^{\frac{1}{2} _+ \frac{2}{p}}} + \lVert b^{3} \rVert_{\dot{H}^{\frac{1}{2} _+ \frac{2}{p}}}) \lVert d_{\frac{3}{4}} \rVert_{L^{2}}^{-\frac{1}{3} + \frac{2}{p}}\lVert \nabla d_{\frac{3}{4}} \rVert_{L^{2}}^{1-\frac{2}{p}}d\tau\nonumber\\
\lesssim& \int_{0}^{t} \left( \lVert \nabla d_{\frac{3}{4}} \rVert_{L^{2}}^{2} + \lVert \nabla \omega_{\frac{3}{4}} \rVert_{L^{2}}^{2} \right)^{1-\frac{1}{p}} \left( \lVert \omega_{\frac{3}{4}} \rVert_{L^{2}}^{2} + \lVert d_{\frac{3}{4}} \rVert_{L^{2}}^{2}\right)^{\frac{1}{p}}\nonumber\\
&\times \left( \lVert u^{3} \rVert_{\dot{H}^{\frac{1}{2} + \frac{2}{p}}} + \lVert b^{3} \rVert_{\dot{H}^{\frac{1}{2} + \frac{2}{p}}}\right)\nonumber\\
&+ (\lVert \partial_{3}^{2} u^{3} \rVert_{\mathcal{H}_{\theta}} + \lVert \partial_{3}^{2} b^{3} \rVert_{\mathcal{H}_{\theta}})\left( \lVert u^{3} \rVert_{\dot{H}^{\frac{1}{2} + \frac{2}{p}}}^{p} + \lVert b^{3} \rVert_{\dot{H}^{\frac{1}{2} + \frac{2}{p}}}^{p}\right)^{\frac{1}{6}}\nonumber\\
&\times[(\lVert u^{3} \rVert_{\dot{H}^{\frac{1}{2} + \frac{2}{p}}}^{p} + \lVert b^{3} \rVert_{\dot{H}^{\frac{1}{2} + \frac{2}{p}}}^{p})(\lVert \omega_{\frac{3}{4}} \rVert_{L^{2}}^{2} + \lVert d_{\frac{3}{4}} \rVert_{L^{2}}^{2})  ]^{\frac{1}{p} - \frac{1}{6}}\lVert \nabla d_{\frac{3}{4}} \rVert_{L^{2}}^{2(\frac{1}{2} - \frac{1}{p})}d\tau\nonumber\\
\leq& \frac{1}{36} \int_{0}^{t} \lVert \nabla \omega_{\frac{3}{4}}\rVert_{L^{2}}^{2} + \lVert \nabla d_{\frac{3}{4}} \rVert_{L^{2}}^{2} d\tau\nonumber\\
&+ c \int_{0}^{t} 
(\lVert u^{3} \rVert_{\dot{H}^{\frac{1}{2} + \frac{2}{p}}}^{p} + \lVert b^{3} \rVert_{\dot{H}^{\frac{1}{2} + \frac{2}{p}}}^{p})(\lVert \omega_{\frac{3}{4}} \rVert_{L^{2}}^{2} + \lVert d_{\frac{3}{4}} \rVert_{L^{2}}^{2})  d\tau\nonumber\\ 
&+ c\left(\int_{0}^{t} \lVert \partial_{3}^{2} u^{3} \rVert_{\mathcal{H}_{\theta}}^{2} + \lVert \partial_{3}^{2} b^{3} \rVert_{\mathcal{H}_{\theta}}^{2} d\tau\right)^{\frac{1}{2}} \left(\int_{0}^{t} \lVert u^{3} \rVert_{\dot{H}^{\frac{1}{2} + \frac{2}{p}}}^{p} + \lVert b^{3} \rVert_{\dot{H}^{\frac{1}{2} + \frac{2}{p}}}^{p} d\tau\right)^{\frac{1}{6}}\nonumber\\
&\times \left(\int_{0}^{t} (\lVert u^{3} \rVert_{\dot{H}^{\frac{1}{2} + \frac{2}{p}}} + \lVert b^{3} \rVert_{\dot{H}^{\frac{1}{2} + \frac{2}{p}}}^{p}) (\lVert \omega_{\frac{3}{4}} \rVert_{L^{2}}^{2} + \lVert d_{\frac{3}{4}} \rVert_{L^{2}}^{2}) d\tau \right)^{\frac{1}{p} - \frac{1}{6}}\left(\int_{0}^{t} \lVert \nabla d_{\frac{3}{4}} \rVert_{L^{2}}^{2} d\tau \right)^{\frac{1}{2} - \frac{1}{p}}\nonumber\\
\leq& \frac{1}{18} \int_{0}^{t} \lVert \nabla \omega_{\frac{3}{4}}\rVert_{L^{2}}^{2} + \lVert \nabla d_{\frac{3}{4}} \rVert_{L^{2}}^{2} d\tau\nonumber\\
&+ c \int_{0}^{t} 
(\lVert u^{3} \rVert_{\dot{H}^{\frac{1}{2} + \frac{2}{p}}}^{p} + \lVert b^{3} \rVert_{\dot{H}^{\frac{1}{2} + \frac{2}{p}}}^{p})(\lVert \omega_{\frac{3}{4}} \rVert_{L^{2}}^{2} + \lVert d_{\frac{3}{4}} \rVert_{L^{2}}^{2})  d\tau\nonumber\\ 
&+ 
c\left(\int_{0}^{t} \lVert \partial_{3}^{2} u^{3} \rVert_{\mathcal{H}_{\theta}}^{2} + \lVert \partial_{3}^{2} b^{3} \rVert_{\mathcal{H}_{\theta}}^{2} d\tau \right)^{\frac{3}{4}}
\left(\int_{0}^{t} \lVert u^{3} \rVert_{\dot{H}^{\frac{1}{2} + \frac{2}{p}}}^{p} + \lVert b^{3} \rVert_{\dot{H}^{\frac{1}{2} + \frac{2}{p}}}^{p} d\tau \right)^{\frac{1}{4}} \nonumber
\end{align}
by (13a) with $(f,g,h) = (\partial_{3}d, u^{3}, d), (\partial_{3} \omega, b^{3}, d)$, (13b) with $(f, g, h) = (\partial_{3}^{2} b^{3}, u^{3}, d), (\partial_{3}^{2} u^{3}, b^{3}, d)$, $\sigma = 3(\frac{1}{2} - \frac{1}{p})$, (9),  Gagliardo-Nirenberg, Young's and H$\ddot{o}$lder's inequalities. 

We now work on the term with all index being one's and two's. We write by (5) 
\begin{align}
I_{5} 
=& 2\int_{0}^{t} \int [\partial_{1}b^{h} \cdot \partial_{2} (\nabla_{h}^{\bot}\Delta_{h}^{-1} \omega - \nabla_{h} \Delta_{h}^{-1} \partial_{3}u^{3})\\
&- \partial_{2}b^{h} \cdot \partial_{1}(\nabla_{h}^{\bot}\Delta_{h}^{-1} \omega - \nabla_{h} \Delta_{h}^{-1} \partial_{3}u^{3})] d_{\frac{1}{2}}  d\tau.\nonumber
\end{align}
We bound by identical estimates applied in (29) to obtain 
\begin{align}
&2\int_{0}^{t} \int [\partial_{1}b^{h} \cdot \partial_{2} \nabla_{h}^{\bot} \Delta_{h}^{-1} \omega - \partial_{2} b^{h} \cdot \partial_{1} \nabla_{h}^{\bot} \Delta_{h}^{-1} \omega ] d_{\frac{1}{2}} d\tau\\
\lesssim& \int_{0}^{t} \lVert \nabla b\rVert_{L^{p_{2}}} \lVert \omega \rVert_{L^{\frac{3}{2}}} \lVert d_{\frac{3}{4}} \rVert_{L^{\frac{2p_{2}}{p_{2} - 3}}}^{\frac{2}{3}}d\tau\nonumber\\
\lesssim& \int_{0}^{t} \lVert \nabla b\rVert_{L^{p_{2}}} \lVert \omega_{\frac{3}{4}} \rVert_{L^{2}}^{\frac{4}{3}} \lVert d_{\frac{3}{4}} \rVert_{L^{2}}^{\frac{2p_{2} - 9}{3p_{2}}}\lVert \nabla d_{\frac{3}{4}} \rVert_{L^{2}}^{\frac{3}{p_{2}}} d\tau\nonumber\\
\leq& \frac{1}{36} \int_{0}^{t} \lVert \nabla d_{\frac{3}{4}} \rVert_{L^{2}}^{2}+ c \int_{0}^{t} \lVert \nabla b\rVert_{L^{p_{2}}}^{\frac{2p_{2}}{2p_{2} - 3}}(\lVert \omega_{\frac{3}{4}} \rVert_{L^{2}}^{2} + \lVert d_{\frac{3}{4}} \rVert_{L^{2}}^{2})d\tau\nonumber
\end{align}
by H$\ddot{o}$lder's, Gagliardo-Nirenberg and Young's inequalities while we bound 
\begin{align}
&2\int_{0}^{t} \int [-\partial_{1} b^{h} \cdot \partial_{2} \nabla_{h} \Delta_{h}^{-1}\partial_{3}u^{3} + \partial_{2}b^{h} \cdot \partial_{1} \nabla_{h} \Delta_{h}^{-1}\partial_{3}u^{3}] d_{\frac{1}{2}} d\tau\\
\lesssim& \int_{0}^{t} \lVert \nabla_{h}\partial_{3} u^{3} \rVert_{\mathcal{H}_{\theta}}\lVert b^{h} \rVert_{\dot{H}^{\frac{3}{2} - \frac{2}{3}\sigma}} \lVert d_{\frac{3}{4}} \rVert_{\dot{H}^{\sigma}}^{\frac{2}{3}}d\tau\nonumber\\
\lesssim& \int_{0}^{t} \lVert \nabla \partial_{3} u^{3} \rVert_{\mathcal{H}_{\theta}} \lVert b\rVert_{\dot{H}^{\frac{1}{2} + \frac{2}{p}}}^{\frac{p}{6}} 
(\lVert b\rVert_{\dot{H}^{\frac{1}{2} + \frac{2}{p}}}^{p} 
\lVert d_{\frac{3}{4}} \rVert_{L^{2}}^{2})^{\frac{1}{p} - \frac{1}{6}} \lVert \nabla d_{\frac{3}{4}} \rVert_{L^{2}}^{2(\frac{1}{2} - \frac{1}{p})}d\tau\nonumber\\
\lesssim& \left(\int_{0}^{t} \lVert \nabla \partial_{3} u^{3} \rVert_{\mathcal{H}_{\theta}}^{2}d\tau \right)^{\frac{1}{2}} \left(\int_{0}^{t} \lVert b\rVert_{\dot{H}^{\frac{1}{2} + \frac{2}{p}}}^{p} d\tau \right)^{\frac{1}{6}}\nonumber\\
&\times \left(\int_{0}^{t} \lVert b\rVert_{\dot{H}^{\frac{1}{2} + \frac{2}{p}}}^{p} \lVert d_{\frac{3}{4}} \rVert_{L^{2}}^{2} d\tau\right)^{\frac{1}{p} - \frac{1}{6}} \left(\int_{0}^{t} \lVert \nabla d_{\frac{3}{4}} \rVert_{L^{2}}^{2} d\tau \right)^{\frac{1}{2} - \frac{1}{p}}\nonumber \\
\leq& \frac{1}{36} \int_{0}^{t} \lVert \nabla d_{\frac{3}{4}} \rVert_{L^{2}}^{2} d\tau\nonumber\\
&+ c \left(\int_{0}^{t} \lVert \nabla \partial_{3} u^{3} \rVert_{\mathcal{H}_{\theta}}^{2}d\tau \right)^{\frac{3}{4}} \left(\int_{0}^{t} \lVert b\rVert_{\dot{H}^{\frac{1}{2} + \frac{2}{p}}}^{p} d\tau \right)^{\frac{1}{4}} +  c\int_{0}^{t} \lVert b\rVert_{\dot{H}^{\frac{1}{2} + \frac{2}{p}}}^{p} \lVert d_{\frac{3}{4}} \rVert_{L^{2}}^{2} d\tau\nonumber
\end{align}
by (13b) with $(f,g,h) = (\nabla_{h} \partial_{3}u^{3}, b^{h}, d), \sigma = 3(\frac{1}{2} - \frac{1}{p})$, Gagliardo-Nirenberg, H$\ddot{o}$lder's and Young's inequalities. Thus, due to (37), (38) in (36) 
\begin{align}
I_{5} \leq& \frac{1}{18} \int_{0}^{t} \lVert \nabla d_{\frac{3}{4}} \rVert_{L^{2}}^{2} d\tau + 
c\int_{0}^{t} (\lVert b\rVert_{\dot{H}^{\frac{1}{2} + \frac{2}{p}}}^{p} + \lVert \nabla b\rVert_{L^{p_{2}}}^{\frac{2p_{2}}{2p_{2} - 3}}) (\lVert \omega_{\frac{3}{4}} \rVert_{L^{2}}^{2} + \lVert d_{\frac{3}{4}} \rVert_{L^{2}}^{2})d\tau\\
&+ c\left(\int_{0}^{t} \lVert \nabla \partial_{3} u^{3} \rVert_{\mathcal{H}_{\theta}}^{2} d\tau \right)^{\frac{3}{4}} \left(\int_{0}^{t} \lVert b\rVert_{\dot{H}^{\frac{1}{2} + \frac{2}{p}}}^{p} d\tau \right)^{\frac{1}{4}}.\nonumber
\end{align}
In sum of (26)-(29), (34), (35), (39), we have 
\begin{align}
&\frac{2}{3} (\lVert \omega_{\frac{3}{4}} \rVert_{L^{2}}^{2} + \lVert d_{\frac{3}{4}} \rVert_{L^{2}}^{2})(t) + \frac{5}{9} \int_{0}^{t} \lVert \nabla \omega_{\frac{3}{4}} \rVert_{L^{2}}^{2} + \lVert \nabla d_{\frac{3}{4}} \rVert_{L^{2}}^{2} d\tau \\
\leq& \frac{2}{3} (\lVert \omega_{\frac{3}{4}} \rVert_{L^{2}}^{2} + \lVert d_{\frac{3}{4}} \rVert_{L^{2}}^{2})(0)\nonumber \\
&+ c \int_{0}^{t} (\lVert \omega_{\frac{3}{4}} \rVert_{L^{2}}^{2} + \lVert d_{\frac{3}{4}} \rVert_{L^{2}}^{2}) (\lVert u^{3}\rVert_{\dot{H}^{\frac{1}{2} + \frac{2}{p}}}^{p} + \lVert b\rVert_{\dot{H}^{\frac{1}{2} + \frac{2}{p}}}^{p} + \lVert b\rVert_{L^{p_{1}}}^{\frac{2p_{1}}{p_{1} - 3}} + \lVert \nabla b\rVert_{L^{p_{2}}}^{\frac{2p_{2}}{2p_{2} - 3}})d\tau\nonumber\\
&+ c \left(\int_{0}^{t} \lVert \nabla \partial_{3} u^{3} \rVert_{\mathcal{H}_{\theta}}^{2} + \lVert \nabla \partial_{3} b^{3} \rVert_{\mathcal{H}_{\theta}}^{2}d\tau \right)^{\frac{3}{4}} \left(\int_{0}^{t} \lVert u^{3} \rVert_{\dot{H}^{\frac{1}{2} + \frac{2}{p}}}^{p} + \lVert b \rVert_{\dot{H}^{\frac{1}{2} + \frac{2}{p}}}^{p}  d\tau \right)^{\frac{1}{4}}\nonumber
\end{align}
so that Gronwall's type inequality argument using 
\begin{align*}
&e^{c\int_{0}^{t}\lVert u^{3}\rVert_{\dot{H}^{\frac{1}{2} + \frac{2}{p}}}^{p} + \lVert b\rVert_{\dot{H}^{\frac{1}{2} + \frac{2}{p}}}^{p} + \lVert b\rVert_{L^{p_{1}}}^{\frac{2p_{1}}{p_{1} - 3}} + \lVert \nabla b\rVert_{L^{p_{2}}}^{\frac{2p_{2}}{2p_{2} - 3}} d\tau}\left(\int_{0}^{t} \lVert u^{3} \rVert_{\dot{H}^{\frac{1}{2} + \frac{2}{p}}}^{p} + \lVert b \rVert_{\dot{H}^{\frac{1}{2} + \frac{2}{p}}}^{p}  d\tau \right)^{\frac{1}{4}}\\
\lesssim& e^{c\int_{0}^{t}\lVert u^{3}\rVert_{\dot{H}^{\frac{1}{2} + \frac{2}{p}}}^{p} + \lVert b\rVert_{\dot{H}^{\frac{1}{2} + \frac{2}{p}}}^{p} + \lVert b\rVert_{L^{p_{1}}}^{\frac{2p_{1}}{p_{1} - 3}} + \lVert \nabla b\rVert_{L^{p_{2}}}^{\frac{2p_{2}}{2p_{2} - 3}} d\tau}
\end{align*}
completes the proof of Proposition 3.1.
\end{proof}

\begin{proposition}
Under the hypothesis of Theorem 1.1, for $\theta \in (\frac{1}{2} - \frac{2}{p}, \frac{1}{6})$, the solution to the MHD system (1a)-(1c) satisfies for any $t < T^{\ast}$
\begin{align*}
&(\lVert \partial_{3} u^{3} \rVert_{\mathcal{H}_{\theta}}^{2} + \lVert \partial_{3} b^{3} \rVert_{\mathcal{H}_{\theta}}^{2})(t) + \int_{0}^{t} \lVert \nabla \partial_{3} u^{3} \rVert_{\mathcal{H}_{\theta}}^{2} + \lVert \nabla \partial_{3} b^{3} \rVert_{\mathcal{H}_{\theta}}^{2}d\tau\\
\lesssim& e^{c\int_{0}^{t}\lVert u^{3} \rVert_{\dot{H}^{\frac{1}{2} + \frac{2}{p}}}^{p} + \lVert b^{3} \rVert_{\dot{H}^{\frac{1}{2} + \frac{2}{p}}}^{p} d\tau }\\
\times&[\lVert \Omega_{0} \rVert_{L^{\frac{3}{2}}}^{2} + \lVert j_{0} \rVert_{L^{\frac{3}{2}}}^{2}\\
&+ \int_{0}^{t} \lVert u^{3} \rVert_{\dot{H}^{\frac{1}{2} + \frac{2}{p}}}(\lVert \omega_{\frac{3}{4}} \rVert_{L^{2}}^{2(\frac{p+3}{3p})}\lVert \nabla \omega_{\frac{3}{4}} \rVert_{L^{2}}^{2(1-\frac{1}{p})} + \lVert d_{\frac{3}{4}} \rVert_{L^{2}}^{2(\frac{p+3}{3p})}\lVert \nabla d_{\frac{3}{4}} \rVert_{L^{2}}^{2(1-\frac{1}{p})})\\
&+ (\lVert u^{3} \rVert_{\dot{H}^{\frac{1}{2} + \frac{2}{p}}}^{2} + \lVert b^{3} \rVert_{\dot{H}^{\frac{1}{2} + \frac{2}{p}}}^{2})(\lVert \omega_{\frac{3}{4}} \rVert_{L^{2}}^{2} + \lVert d_{\frac{3}{4}} \rVert_{L^{2}}^{2})^{\frac{p+6}{3p}}(\lVert \nabla \omega_{\frac{3}{4}} \rVert_{L^{2}}^{2} + \lVert \nabla d_{\frac{3}{4}} \rVert_{L^{2}}^{2})^{1-\frac{2}{p}}d\tau].
\end{align*}

\end{proposition}

\begin{proof}
Applying $\partial_{3}$ on the third components of (1a), (1b), we obtain 
\begin{subequations}
\begin{align}
&\partial_{t} \partial_{3} u^{3} - \Delta \partial_{3} u^{3} = -\partial_{3} u \cdot\nabla u^{3} - (u\cdot\nabla) \partial_{3} u^{3} + \partial_{3} b\cdot\nabla b^{3} + (b\cdot\nabla) \partial_{3} b^{3}\\
& \hspace{25mm} - \partial_{3}^{2} (-\Delta)^{-1} \sum_{l,m=1}^{3} (\partial_{l} u^{m} \partial_{m} u^{l} - \partial_{l} b^{m} \partial_{m} b^{l}), \nonumber\\
&\partial_{t} \partial_{3}b^{3} - \Delta \partial_{3} b^{3} = -\partial_{3} u\cdot\nabla b^{3} - (u\cdot\nabla) \partial_{3} b^{3} + \partial_{3} b\cdot\nabla u^{3} + (b\cdot\nabla) \partial_{3} u^{3}.
\end{align}
\end{subequations}
We write 
\begin{align*}
&\sum_{l,m=1}^{3} \partial_{l} u^{m} \partial_{m} u^{l} - \partial_{l} b^{m} \partial_{m} b^{l} \\
=& \sum_{l,m=1}^{2} \partial_{l} u^{m} \partial_{m} u^{l} - \partial_{l} b^{m} \partial_{m} b^{l} + 2\sum_{l=1}^{2} \partial_{l} u^{3} \partial_{3} u^{l} - \partial_{l} b^{3} \partial_{3} b^{l} + (\partial_{3} u^{3})^{2} - (\partial_{3} b^{3})^{2},
\end{align*}
and 
\begin{align*}
&-\partial_{3} u\cdot\nabla u^{3} + \partial_{3} b\cdot\nabla b^{3}  = \sum_{l=1}^{2} -\partial_{3} u^{l} \partial_{l} u^{3} + \partial_{3} b^{l} \partial_{l} b^{3} - (\partial_{3} u^{3})^{2} + (\partial_{3} b^{3})^{2},\\
&-\partial_{3} u \cdot\nabla b^{3} + \partial_{3} b\cdot\nabla u^{3} =\sum_{l=1}^{2} -\partial_{3} u^{l} \partial_{l} b^{3} + \partial_{3} b^{l} \partial_{l} u^{3} 
\end{align*}
so that we can take $\mathcal{H}_{\theta}$-inner products on (41a), (41b) and sum to obtain 
\begin{align}
&\frac{1}{2} \partial_{t} (\lVert \partial_{3} u^{3} \rVert_{\mathcal{H}_{\theta}}^{2} + \lVert \partial_{3} b^{3} \rVert_{\mathcal{H}_{\theta}}^{2}) + \lVert \nabla \partial_{3} u^{3} \rVert_{\mathcal{H}_{\theta}}^{2} + \lVert \nabla \partial_{3} b^{3} \rVert_{\mathcal{H}_{\theta}}^{2}\\
=& \sum_{n=1}^{3}(Q_{n} (u, b) \lvert \partial_{3} u^{3} )_{\mathcal{H}_{\theta}} + \sum_{n=1}^{2} (R_{n} (u, b) \lvert \partial_{3} b^{3})_{\mathcal{H}_{\theta}}\nonumber
\end{align}
where
\begin{subequations}
\begin{align}
&Q_{1}(u,b) \triangleq (-Id - \partial_{3}^{2} (-\Delta)^{-1} ) ((\partial_{3} u^{3})^{2} - (\partial_{3} b^{3})^{2})\\
& \hspace{15mm} - \partial_{3}^{2} (-\Delta)^{-1} \sum_{l,m=1}^{2} (\partial_{l} u^{m} \partial_{m} u^{l} - \partial_{l} b^{m} \partial_{m} b^{l} ), \nonumber\\
&Q_{2}(u,b) \triangleq (Id + 2\partial_{3}^{2} (-\Delta)^{-1} ) (\sum_{l=1}^{2} -\partial_{3} u^{l} \partial_{l} u^{3} + \partial_{3} b^{l} \partial_{l} b^{3}),\\
&Q_{3}(u,b) \triangleq -(u\cdot\nabla) \partial_{3} u^{3} + (b\cdot\nabla) \partial_{3} b^{3},\\
&R_{1}(u,b) \triangleq \sum_{l=1}^{2} -\partial_{3} u^{l} \partial_{l} b^{3} + \partial_{3} b^{l} \partial_{l} u^{3},\\
&R_{2}(u,b) \triangleq -(u\cdot\nabla) \partial_{3} b^{3} + (b\cdot\nabla) \partial_{3}u^{3}.
\end{align}
\end{subequations}
We estimate first 
\begin{align}
&\lvert (Q_{1}(u, b) \lvert \partial_{3} u^{3})_{\mathcal{H}_{\theta}}\rvert\\
\leq& \lvert ( (-Id - \partial_{3}^{2} (-\Delta)^{-1} ) ((\partial_{3} u^{3})^{2} - (\partial_{3} b^{3})^{2}) \lvert \partial_{3} u^{3})_{\mathcal{H}_{\theta}}\rvert\nonumber\\
&+ \lvert (\partial_{3}^{2} (-\Delta)^{-1} \sum_{l,m=1}^{2} (\partial_{l} u^{m} \partial_{m} u^{l} - \partial_{l} b^{m} \partial_{m} b^{l}) \lvert \partial_{3} u^{3})_{\mathcal{H}_{\theta}}\rvert \triangleq II_{1} + II_{2}\nonumber
\end{align}
where 
\begin{align}
II_{1} \lesssim& \lVert u^{3} \rVert_{\dot{H}^{\frac{1}{2} + \frac{2}{p}}}(\lVert \partial_{3} u^{3} \rVert_{\dot{H}^{\theta, \frac{1}{2} - \theta - \frac{1}{p}}}^{2} + \lVert \partial_{3} b^{3} \rVert_{\dot{H}^{\theta, \frac{1}{2} - \theta - \frac{1}{p}}}^{2})\\
\lesssim& \lVert u^{3} \rVert_{\dot{H}^{\frac{1}{2} + \frac{2}{p}}}(\lVert \partial_{3} u^{3} \rVert_{\mathcal{H}_{\theta}}^{\frac{2}{p}}\lVert \nabla \partial_{3} u^{3} \rVert_{\mathcal{H}_{\theta}}^{2(1-\frac{1}{p})} + \lVert \partial_{3}b^{3} \rVert_{\mathcal{H}_{\theta}}^{\frac{2}{p} }\lVert \nabla \partial_{3} b^{3} \rVert_{\mathcal{H}_{\theta}}^{2(1-\frac{1}{p})})\nonumber\\
\leq& \frac{1}{36} (\lVert \nabla \partial_{3} u^{3} \rVert_{\mathcal{H}_{\theta}}^{2} + \lVert \nabla \partial_{3} b^{3} \rVert_{\mathcal{H}_{\theta}}^{2}) + c\lVert u^{3} \rVert_{\dot{H}^{\frac{1}{2} + \frac{2}{p}}}^{p} (\lVert \partial_{3} u^{3} \rVert_{\mathcal{H}_{\theta}}^{2} + \lVert \partial_{3} b^{3} \rVert_{\mathcal{H}_{\theta}}^{2} )\nonumber
\end{align}
by (14) with $A = -Id - \partial_{3}^{2} (-\Delta)^{-1}, f = g = \partial_{3}u^{3}, h = u^{3}$ and again with $f = g = \partial_{3}b^{3}, h = u^{3}$ for $\theta$ such that $0 < \theta < \frac{1}{2} - \frac{1}{p}$, (12) and Young's inequality. Similarly, 
\begin{align}
II_{2} 
\lesssim& \lVert u^{3} \rVert_{\dot{H}^{\frac{1}{2} + \frac{2}{p}}}\\
&\times (\lVert \omega \rVert_{\dot{H}^{\theta, \frac{1}{2} - \theta - \frac{1}{p}}}^{2} + \lVert \partial_{3} u^{3} \rVert_{\dot{H}^{\theta, \frac{1}{2} - \theta - \frac{1}{p}}}^{2} + \lVert  d\rVert_{\dot{H}^{\theta, \frac{1}{2} - \theta - \frac{1}{p}}}^{2} + \lVert\partial_{3} b^{3} \rVert_{\dot{H}^{\theta, \frac{1}{2} - \theta - \frac{1}{p} }}^{2})\nonumber\\
\lesssim& \lVert u^{3} \rVert_{\dot{H}^{\frac{1}{2} + \frac{2}{p}}}(\lVert \omega_{\frac{3}{4}} \rVert_{L^{2}}^{2(\frac{p+3}{3p})}\lVert \nabla \omega_{\frac{3}{4}} \rVert_{L^{2}}^{2(1-\frac{1}{p})} + \lVert d_{\frac{3}{4}} \rVert_{L^{2}}^{2(\frac{p+3}{3p})}\lVert \nabla d_{\frac{3}{4}} \rVert_{L^{2}}^{2(1-\frac{1}{p})}\nonumber\\
&\hspace{15mm} + \lVert \partial_{3} u^{3} \rVert_{\mathcal{H}_{\theta}}^{\frac{2}{p}} \lVert \nabla \partial_{3} u^{3} \rVert_{\mathcal{H}_{\theta}}^{2-\frac{2}{p}} + \lVert \partial_{3} b^{3} \rVert_{\mathcal{H}_{\theta}}^{\frac{2}{p}} \lVert \nabla \partial_{3} b^{3} \rVert_{\mathcal{H}_{\theta}}^{2-\frac{2}{p}})\nonumber\\
\leq& \frac{1}{36} (\lVert \nabla \partial_{3} u^{3} \rVert_{\mathcal{H}_{\theta}}^{2} + \lVert \nabla \partial_{3} b^{3} \rVert_{\mathcal{H}_{\theta}}^{2})+ c \lVert u^{3} \rVert_{\dot{H}^{\frac{1}{2} + \frac{2}{p}}}^{p} (\lVert \partial_{3} u^{3} \rVert_{\mathcal{H}_{\theta}}^{2} + \lVert \partial_{3} b^{3} \rVert_{\mathcal{H}_{\theta}}^{2})\nonumber\\
&+ c \lVert u^{3} \rVert_{\dot{H}^{\frac{1}{2} + \frac{2}{p}}}(\lVert \omega_{\frac{3}{4}}\rVert_{L^{2}}^{2(\frac{p+3}{3p})} \lVert \nabla \omega_{\frac{3}{4}} \rVert_{L^{2}}^{2(1-\frac{1}{p})} + \lVert d_{\frac{3}{4}}\rVert_{L^{2}}^{2(\frac{p+3}{3p})} \lVert \nabla d_{\frac{3}{4}} \rVert_{L^{2}}^{2(1-\frac{1}{p})} ) \nonumber
\end{align}
by (5), (14), (15), (12) and Young's inequality. 

We now consider 
\begin{align}
(Q_{2}(u,b)\lvert \partial_{3}u^{3})_{\mathcal{H}_{\theta}}
=& ((Id + 2\partial_{3}^{2} (-\Delta)^{-1})(\sum_{l=1}^{2} u^{l} \partial_{l} u^{3}- b^{l} \partial_{l} b^{3} ) \lvert \partial_{3}^{2} u^{3})_{\mathcal{H}_{\theta}}\nonumber\\
&+ ((Id + 2\partial_{3}^{2} (-\Delta)^{-1}) (\sum_{l=1}^{2} u^{l} \partial_{l} \partial_{3} u^{3} - b^{l} \partial_{l} \partial_{3} b^{3}) \lvert \partial_{3} u^{3})_{\mathcal{H}_{\theta}}\nonumber\\
\triangleq& (Q_{2,1}(u,b) \lvert \partial_{3} u^{3})_{\mathcal{H}_{\theta}} + (Q_{2,2}(u,b) \lvert \partial_{3} u^{3})_{\mathcal{H}_{\theta}}
\end{align}
due to integration by parts. We estimate 
\begin{align}
&\lvert (Q_{2,1} (u,b) \lvert \partial_{3}u^{3})_{\mathcal{H}_{\theta}}\rvert\\
\lesssim& \sum_{l=1}^{2} (\lVert u^{l} \partial_{l} u^{3} \rVert_{\dot{H}^{-\frac{1}{2} + \theta, -\theta}} + \lVert b^{l} \partial_{l} b^{3} \rVert_{\dot{H}^{-\frac{1}{2} + \theta, -\theta}})\lVert \partial_{3}^{2} u^{3} \rVert_{\mathcal{H}_{\theta}}\nonumber\\
\lesssim& \sum_{l=1}^{2} (\lVert u^{l} \rVert_{(\dot{B}_{2,1}^{1})_{h}(\dot{B}_{2,1}^{\frac{1}{2} - \frac{2}{p}})_{v}}\lVert \partial_{l} u^{3} \rVert_{\dot{H}^{-\frac{1}{2} + \theta, \frac{2}{p} - \theta}}\nonumber\\
& \hspace{35mm} + \lVert b^{l} \rVert_{(\dot{B}_{2,1}^{1})_{h}(\dot{B}_{2,1}^{\frac{1}{2} - \frac{2}{p}})_{v}}\lVert \partial_{l} b^{3} \rVert_{\dot{H}^{-\frac{1}{2} + \theta, \frac{2}{p} - \theta}})\lVert \partial_{3}^{2} u^{3} \rVert_{\mathcal{H}_{\theta}}\nonumber\\
\lesssim& \sum_{l=1}^{2} (\lVert u^{l} \rVert_{(\dot{B}_{2,1}^{1})_{h}(\dot{B}_{2,1}^{\frac{1}{2} - \frac{2}{p}})_{v}}\lVert u^{3} \rVert_{\dot{H}^{\frac{1}{2} + \frac{2}{p}}} + \lVert b^{l} \rVert_{(\dot{B}_{2,1}^{1})_{h}(\dot{B}_{2,1}^{\frac{1}{2} - \frac{2}{p}})_{v}}\lVert b^{3} \rVert_{\dot{H}^{\frac{1}{2} + \frac{2}{p}}})\lVert \partial_{3}^{2} u^{3} \rVert_{\mathcal{H}_{\theta}}\nonumber\\
\lesssim& \left(\lVert \omega_{\frac{3}{4}} \rVert_{L^{2}}^{\frac{1}{3} + \frac{2}{p}} \lVert \nabla \omega_{\frac{3}{4}} \rVert_{L^{2}}^{1-\frac{2}{p}} + \lVert \partial_{3} u^{3} \rVert_{\mathcal{H}_{\theta}}^{\frac{2}{p}}\lVert \nabla \partial_{3} u^{3} \rVert_{\mathcal{H}_{\theta}}^{1-\frac{2}{p}}\right) \lVert u^{3} \rVert_{\dot{H}^{\frac{1}{2} + \frac{2}{p}}}\lVert \partial_{3}^{2} u^{3} \rVert_{\mathcal{H}_{\theta}}\nonumber\\
&+ \left(\lVert d_{\frac{3}{4}} \rVert_{L^{2}}^{\frac{1}{3} + \frac{2}{p}} \lVert \nabla d_{\frac{3}{4}} \rVert_{L^{2}}^{1-\frac{2}{p}} + \lVert \partial_{3} b^{3} \rVert_{\mathcal{H}_{\theta}}^{\frac{2}{p}}\lVert \nabla \partial_{3} b^{3} \rVert_{\mathcal{H}_{\theta}}^{1-\frac{2}{p}}\right) \lVert b^{3} \rVert_{\dot{H}^{\frac{1}{2} + \frac{2}{p}}}\lVert \partial_{3}^{2} u^{3} \rVert_{\mathcal{H}_{\theta}}\nonumber\\
\leq& \frac{1}{36} (\lVert \nabla \partial_{3} u^{3} \rVert_{\mathcal{H}_{\theta}}^{2} + \lVert \nabla \partial_{3} b^{3} \rVert_{\mathcal{H}_{\theta}}^{2})\nonumber\\
&+ c( \lVert u^{3} \rVert_{\dot{H}^{\frac{1}{2} + \frac{2}{p}}}^{2} \lVert \omega_{\frac{3}{4}} \rVert_{L^{2}}^{2(\frac{p+6}{3p})}\lVert \nabla \omega_{\frac{3}{4}} \rVert_{L^{2}}^{2(1-\frac{2}{p})} + \lVert u^{3} \rVert_{\dot{H}^{\frac{1}{2} + \frac{2}{p}}}^{p}\lVert \partial_{3} u^{3} \rVert_{\mathcal{H}_{\theta}}^{2}\nonumber\\
&+ \lVert b^{3} \rVert_{\dot{H}^{\frac{1}{2} + \frac{2}{p}}}^{2}\lVert d_{\frac{3}{4}} \rVert_{L^{2}}^{2(\frac{p+6}{3p})}\lVert \nabla d_{\frac{3}{4}} \rVert_{L^{2}}^{2(1-\frac{2}{p})} + \lVert b^{3} \rVert_{\dot{H}^{\frac{1}{2} + \frac{2}{p}}}^{p} \lVert \partial_{3} b^{3} \rVert_{\mathcal{H}_{\theta}}^{2})\nonumber
\end{align}
where we used continuity of Riesz transform, (94), the fact that 
$\sum_{l=1}^{2} \lVert \partial_{l} f\rVert_{\dot{H}^{-\frac{1}{2} + \theta, \frac{2}{p} - \theta}} \lesssim \lVert f\rVert_{\dot{H}^{\frac{1}{2} + \frac{2}{p}}}$, 
Lemma 2.2 with $\alpha = \frac{2}{p}$ and Young's inequalities. Next, 
\begin{align}
&\lvert (Q_{2,2}(u,b) \lvert \partial_{3}u^{3})_{\mathcal{H}_{\theta}}\rvert\\
\lesssim& (\lVert \omega_{\frac{3}{4}} \rVert_{L^{2}}^{\frac{1}{3} + \frac{2}{p}}\lVert \nabla \omega_{\frac{3}{4}} \rVert_{L^{2}}^{1-\frac{2}{p}} + \lVert \partial_{3} u^{3} \rVert_{\mathcal{H}_{\theta}}^{\frac{2}{p}} \lVert \nabla \partial_{3} u^{3} \rVert_{\mathcal{H}_{\theta}}^{1-\frac{2}{p}})\lVert \nabla \partial_{3} u^{3} \rVert_{\mathcal{H}_{\theta}}\lVert u^{3} \rVert_{\dot{H}^{\frac{1}{2} + \frac{2}{p}}}\nonumber\\
&+ (\lVert d_{\frac{3}{4}} \rVert_{L^{2}}^{\frac{1}{3} + \frac{2}{p}}\lVert \nabla d_{\frac{3}{4}} \rVert_{L^{2}}^{1-\frac{2}{p}} + \lVert \partial_{3} b^{3} \rVert_{\mathcal{H}_{\theta}}^{\frac{2}{p}} \lVert \nabla \partial_{3} b^{3} \rVert_{\mathcal{H}_{\theta}}^{1-\frac{2}{p}})\lVert \nabla \partial_{3} b^{3} \rVert_{\mathcal{H}_{\theta}}\lVert u^{3} \rVert_{\dot{H}^{\frac{1}{2} + \frac{2}{p}}}\nonumber\\
\leq& \frac{1}{36} (\lVert \nabla\partial_{3} u^{3} \rVert_{\mathcal{H}_{\theta}}^{2} + \lVert \nabla\partial_{3} b^{3} \rVert_{\mathcal{H}_{\theta}}^{2})\nonumber\\
&+ c(\lVert u^{3} \rVert_{\dot{H}^{\frac{1}{2} + \frac{2}{p}}}^{2} \lVert \omega_{\frac{3}{4}} \rVert_{L^{2}}^{2(\frac{p+6}{3p})}\lVert \nabla \omega_{\frac{3}{4}} \rVert_{L^{2}}^{2(1-\frac{2}{p})} + \lVert u^{3} \rVert_{\dot{H}^{\frac{1}{2} + \frac{2}{p}}}^{p} \lVert \partial_{3} u^{3} \rVert_{\mathcal{H}_{\theta}}^{2}\nonumber\\
&+ \lVert u^{3} \rVert_{\dot{H}^{\frac{1}{2} + \frac{2}{p}}}^{2} \lVert d_{\frac{3}{4}} \rVert_{L^{2}}^{2(\frac{p+6}{3p})}\lVert \nabla d_{\frac{3}{4}} \rVert_{L^{2}}^{2(1-\frac{2}{p})} + \lVert u^{3} \rVert_{\dot{H}^{\frac{1}{2} + \frac{2}{p}}}^{p} \lVert \partial_{3} b^{3} \rVert_{\mathcal{H}_{\theta}}^{2} )\nonumber
\end{align}
by (16) with $A = Id + 2 \partial_{3}^{2} (-\Delta)^{-1}, g = u, h = u$, (17) with $A = Id + 2 \partial_{3}^{2} (-\Delta)^{-1}, g = b, h = u$ and Young's inequalities. 

We treat $R_{1}$ similarly: 
\begin{align}
&(R_{1}(u, b) \lvert \partial_{3} b^{3})_{\mathcal{H}_{\theta}}\\
=& (\sum_{l=1}^{2} u^{l} \partial_{l} b^{3} - b^{l} \partial_{l} u^{3} \lvert \partial_{3}^{2} b^{3})_{\mathcal{H}_{\theta}} + (\sum_{l=1}^{2} u^{l} \partial_{l} \partial_{3} b^{3} - b^{l} \partial_{l}\partial_{3} u^{3} \lvert \partial_{3} b^{3})_{\mathcal{H}_{\theta}}\nonumber\\
\triangleq& (R_{1,1} (u,b) \lvert \partial_{3}^{2} b^{3})_{\mathcal{H}_{\theta}} + (R_{1,2} (u,b) \lvert \partial_{3} b^{3})_{\mathcal{H}_{\theta}}\nonumber
\end{align}
where
\begin{align}
&\lvert (R_{1,1} (u,b) \lvert \partial_{3}^{2} b^{3})_{\mathcal{H}_{\theta}}\rvert\\
\lesssim& \sum_{l=1}^{2}(\lVert u^{l} \partial_{l} b^{3} \rVert_{\mathcal{H}_{\theta}} + \lVert b^{l} \partial_{l} u^{3} \rVert_{\mathcal{H}_{\theta}} )\lVert \partial_{3}^{2} u^{3} \rVert_{\mathcal{H}_{\theta}}\nonumber\\
\lesssim& \sum_{l=1}^{2} (\lVert u^{l} \rVert_{(\dot{B}_{2,1}^{1})_{h}(\dot{B}_{2,1}^{\frac{1}{2} - \frac{2}{p}})_{v}} \lVert b^{3} \rVert_{\dot{H}^{\frac{1}{2} + \frac{2}{p}}} + \lVert b^{l} \rVert_{(\dot{B}_{2,1}^{1})_{h}(\dot{B}_{2,1}^{\frac{1}{2} - \frac{2}{p}})_{v}} \lVert u^{3} \rVert_{\dot{H}^{\frac{1}{2} + \frac{2}{p}}})\lVert \partial_{3}^{2} u^{3} \rVert_{\mathcal{H}_{\theta}}\nonumber\\
\lesssim& (\lVert \omega_{\frac{3}{4}} \rVert_{L^{2}}^{\frac{1}{3} + \frac{2}{p}}\lVert \nabla \omega_{\frac{3}{4}} \rVert_{L^{2}}^{1-\frac{2}{p}} + \lVert \partial_{3} u^{3} \rVert_{\mathcal{H}_{\theta}}^{\frac{2}{p}} \lVert \nabla \partial_{3} u^{3} \rVert_{\mathcal{H}_{\theta}}^{1-\frac{2}{p}})\lVert b^{3} \rVert_{\dot{H}^{\frac{1}{2} + \frac{2}{p}}}\lVert \partial_{3}^{2} u^{3} \rVert_{\mathcal{H}_{\theta}}\nonumber\\
&+ (\lVert d_{\frac{3}{4}} \rVert_{L^{2}}^{\frac{1}{3} + \frac{2}{p}}\lVert \nabla d_{\frac{3}{4}} \rVert_{L^{2}}^{1-\frac{2}{p}} + \lVert \partial_{3} b^{3} \rVert_{\mathcal{H}_{\theta}}^{\frac{2}{p}} \lVert \nabla\partial_{3} b^{3} \rVert_{\mathcal{H}_{\theta}}^{1-\frac{2}{p}})\lVert u^{3} \rVert_{\dot{H}^{\frac{1}{2} + \frac{2}{p}}}\lVert\partial_{3}^{2} u^{3} \rVert_{\mathcal{H}_{\theta}}\nonumber\\
\leq& \frac{1}{36} (\lVert \nabla\partial_{3} u^{3} \rVert_{\mathcal{H}_{\theta}}^{2} + \lVert \nabla\partial_{3} b^{3} \rVert_{\mathcal{H}_{\theta}}^{2})\nonumber\\
&+ c(\lVert b^{3} \rVert_{\dot{H}^{\frac{1}{2} + \frac{2}{p}}}^{2} \lVert \omega_{\frac{3}{4}} \rVert_{L^{2}}^{2(\frac{p+6}{3p})}\lVert \nabla \omega_{\frac{3}{4}} \rVert_{L^{2}}^{2(1-\frac{2}{p})} + \lVert b^{3} \rVert_{\dot{H}^{\frac{1}{2} + \frac{2}{p}}}^{p} \lVert \partial_{3} u^{3} \rVert_{\mathcal{H}_{\theta}}^{2}\nonumber\\
&+ \lVert u^{3} \rVert_{\dot{H}^{\frac{1}{2} + \frac{2}{p}}}^{2} \lVert d_{\frac{3}{4}} \rVert_{L^{2}}^{2(\frac{p+6}{3p})}\lVert \nabla d_{\frac{3}{4}} \rVert_{L^{2}}^{2(1-\frac{2}{p})} + \lVert u^{3} \rVert_{\dot{H}^{\frac{1}{2} + \frac{2}{p}}}^{p} \lVert \partial_{3} b^{3} \rVert_{\mathcal{H}_{\theta}}^{2} )\nonumber
\end{align}
by (94), Lemma 2.2 with $\alpha = \frac{2}{p}$ and Young's inequalities. Next, we work on 
\begin{align}
&\lvert (R_{1,2} (u,b) \lvert \partial_{3}b^{3})_{\mathcal{H}_{\theta}}\rvert\\
\lesssim& \left( \lVert \omega_{\frac{3}{4}} \rVert_{L^{2}}^{\frac{1}{3} + \frac{2}{p}}\lVert \nabla \omega_{\frac{3}{4}} \rVert_{L^{2}}^{1-\frac{2}{p}} + \lVert \partial_{3} u^{3} \rVert_{\mathcal{H}_{\theta}}^{\frac{2}{p}}\lVert \nabla \partial_{3} u^{3} \rVert_{\mathcal{H}_{\theta}}^{1-\frac{2}{p}}\right) \lVert \nabla \partial_{3} b^{3} \rVert_{\mathcal{H}_{\theta}}\lVert b^{3} \rVert_{\dot{H}^{\frac{1}{2} + \frac{2}{p}}}\nonumber\\
&+  \left( \lVert d_{\frac{3}{4}} \rVert_{L^{2}}^{\frac{1}{3} + \frac{2}{p}}\lVert \nabla d_{\frac{3}{4}} \rVert_{L^{2}}^{1-\frac{2}{p}} + \lVert \partial_{3} b^{3} \rVert_{\mathcal{H}_{\theta}}^{\frac{2}{p}}\lVert \nabla \partial_{3} b^{3} \rVert_{\mathcal{H}_{\theta}}^{1-\frac{2}{p}}\right) \lVert \nabla \partial_{3} u^{3} \rVert_{\mathcal{H}_{\theta}}\lVert b^{3} \rVert_{\dot{H}^{\frac{1}{2} + \frac{2}{p}}}\nonumber\\
\leq& \frac{1}{36} (\lVert \nabla\partial_{3} u^{3} \rVert_{\mathcal{H}_{\theta}}^{2} + \lVert \nabla\partial_{3} b^{3} \rVert_{\mathcal{H}_{\theta}}^{2})\nonumber\\
&+ c(\lVert b^{3} \rVert_{\dot{H}^{\frac{1}{2} + \frac{2}{p}}}^{2} \lVert \omega_{\frac{3}{4}} \rVert_{L^{2}}^{2(\frac{p+6}{3p})}\lVert \nabla \omega_{\frac{3}{4}} \rVert_{L^{2}}^{2(1-\frac{2}{p})} + \lVert b^{3} \rVert_{\dot{H}^{\frac{1}{2} + \frac{2}{p}}}^{p} \lVert \partial_{3} u^{3} \rVert_{\mathcal{H}_{\theta}}^{2}\nonumber\\
&+ \lVert b^{3} \rVert_{\dot{H}^{\frac{1}{2} + \frac{2}{p}}}^{2} \lVert d_{\frac{3}{4}} \rVert_{L^{2}}^{2(\frac{p+6}{3p})}\lVert \nabla d_{\frac{3}{4}} \rVert_{L^{2}}^{2(1-\frac{2}{p})} + \lVert b^{3} \rVert_{\dot{H}^{\frac{1}{2} + \frac{2}{p}}}^{p} \lVert \partial_{3} b^{3} \rVert_{\mathcal{H}_{\theta}}^{2} )\nonumber
\end{align}
by (16) with $A = Id, g = b, h = b$ and (17) with $A = Id, g =u, h = b$ and Young's inequalities. 

We finally work on 
\begin{align}
&(Q_{3}(u,b) \lvert \partial_{3}u^{3})_{\mathcal{H}_{\theta}} + (R_{2}(u,b) \lvert \partial_{3} b^{3})_{\mathcal{H}_{\theta}}\\
=& (-u^{h} \cdot\nabla_{h} \partial_{3} u^{3}+ b^{h} \cdot \nabla_{h} \partial_{3} b^{3} \lvert \partial_{3} u^{3})_{\mathcal{H}_{\theta}} + (-u^{3} \partial_{3}^{2} u^{3} + b^{3} \partial_{3}^{2} b^{3}  \lvert \partial_{3}u^{3})_{\mathcal{H}_{\theta}}\nonumber\\
&+ (-u^{h} \cdot\nabla_{h} \partial_{3} b^{3} + b^{h} \cdot\nabla_{h} \partial_{3} u^{3} \lvert \partial_{3}b^{3})_{\mathcal{H}_{\theta}} + (-u^{3}\partial_{3}^{2} b^{3} + b^{3}\partial_{3}^{2} u^{3}\lvert \partial_{3}b^{3})_{\mathcal{H}_{\theta}}\nonumber\\
\triangleq& (Q_{3,1} (u,b) \lvert \partial_{3} u^{3})_{\mathcal{H}_{\theta}} + (Q_{3,2} (u,b) \lvert \partial_{3} u^{3})_{\mathcal{H}_{\theta}} + (R_{2,1}(u,b) \lvert \partial_{3} b^{3})_{\mathcal{H}_{\theta}} + (R_{2,2}(u,b) \lvert \partial_{3} b^{3})_{\mathcal{H}_{\theta}}\nonumber
\end{align}
where 
\begin{align}
&(Q_{3,1} (u,b) \lvert \partial_{3} u^{3})_{\mathcal{H}_{\theta}} + (R_{2,1}(u,b) \lvert \partial_{3} b^{3})_{\mathcal{H}_{\theta}}\\
\lesssim& (\lVert \omega_{\frac{3}{4}} \rVert_{L^{2}}^{\frac{1}{3} + \frac{2}{p}}\lVert \nabla \omega_{\frac{3}{4}} \rVert_{L^{2}}^{1-\frac{2}{p}} + \lVert \partial_{3} u^{3} \rVert_{\mathcal{H}_{\theta}}^{\frac{2}{p}} \lVert \nabla \partial_{3} u^{3} \rVert_{\mathcal{H}_{\theta}}^{1-\frac{2}{p}}\nonumber \\
&\hspace{15mm} + \lVert d_{\frac{3}{4}} \rVert_{L^{2}}^{\frac{1}{3} + \frac{2}{p}}\lVert \nabla d_{\frac{3}{4}} \rVert_{L^{2}}^{1-\frac{2}{p}} + \lVert \partial_{3} b^{3} \rVert_{\mathcal{H}_{\theta}}^{\frac{2}{p}} \lVert \nabla \partial_{3} b^{3} \rVert_{\mathcal{H}_{\theta}}^{1-\frac{2}{p}}
)\nonumber\\
&\times (\lVert \nabla \partial_{3} u^{3} \rVert_{\mathcal{H}_{\theta}}\lVert u^{3} \rVert_{\dot{H}^{\frac{1}{2} + \frac{2}{p}}} 
+ \lVert \nabla \partial_{3} b^{3} \rVert_{\mathcal{H}_{\theta}}\lVert u^{3} \rVert_{\dot{H}^{\frac{1}{2} + \frac{2}{p}}}\nonumber\\
& \hspace{15mm} + \lVert \nabla \partial_{3} b^{3} \rVert_{\mathcal{H}_{\theta}} \lVert b^{3} \rVert_{\dot{H}^{\frac{1}{2} + \frac{2}{p}}} + \lVert \nabla \partial_{3} u^{3} \rVert_{\mathcal{H}_{\theta}} \lVert b^{3} \rVert_{\dot{H}^{\frac{1}{2} + \frac{2}{p}}} )\nonumber
\end{align}
by (16) with $A = Id, (g,h) = (u,u), (b,b)$ and (17) with $A = Id, (g,h) = (b,u), (u,b)$. On the other hand, 
\begin{align}
&(Q_{3,2} (u,b) \lvert \partial_{3} u^{3})_{\mathcal{H}_{\theta}} + (R_{2,2}(u,b) \lvert \partial_{3} b^{3})_{\mathcal{H}_{\theta}}\\
\lesssim& (\lVert u^{3}\rVert_{(\dot{B}_{2,2}^{\frac{2}{p}})_{h}(\dot{B}_{2,1}^{\frac{1}{2}})_{v}}\lVert \partial_{3}^{2} u^{3} \rVert_{
\dot{H}^{-\frac{1}{2} + \theta, -\theta}
}\nonumber\\
& \hspace{30mm} + \lVert b^{3}\rVert_{(\dot{B}_{2,2}^{\frac{2}{p}})_{h}(\dot{B}_{2,1}^{\frac{1}{2}})_{v}}\lVert \partial_{3}^{2} b^{3} \rVert_{\dot{H}^{-\frac{1}{2} + \theta, -\theta}})\lVert \partial_{3} u^{3} \rVert_{\dot{H}^{\frac{1}{2} + \theta - \frac{2}{p}, - \theta}}\nonumber\\
&+ (\lVert u^{3}\rVert_{(\dot{B}_{2,2}^{\frac{2}{p}})_{h}(\dot{B}_{2,1}^{\frac{1}{2}})_{v}}\lVert \partial_{3}^{2} b^{3} \rVert_{\dot{H}^{-\frac{1}{2} + \theta, -\theta}}\nonumber\\
& \hspace{30mm} + \lVert b^{3}\rVert_{(\dot{B}_{2,2}^{\frac{2}{p}})_{h}(\dot{B}_{2,1}^{\frac{1}{2}})_{v}}\lVert \partial_{3}^{2} u^{3} \rVert_{\dot{H}^{-\frac{1}{2} + \theta, -\theta}}) 
\lVert \partial_{3} b^{3} \rVert_{\dot{H}^{\frac{1}{2} + \theta - \frac{2}{p}, - \theta}}\nonumber\\
\lesssim& \lVert u^{3} \rVert_{\dot{H}^{\frac{1}{2} + \frac{2}{p}}}(\lVert \partial_{3}^{2} u^{3} \rVert_{\mathcal{H}_{\theta}} \lVert \partial_{3} u^{3} \rVert_{\mathcal{H}_{\theta}}^{\frac{2}{p}} \lVert \nabla_{h} \partial_{3} u^{3} \rVert_{\mathcal{H}_{\theta}}^{1-\frac{2}{p}} + \lVert \partial_{3}^{2} b^{3} \rVert_{\mathcal{H}_{\theta}} \lVert \partial_{3} b^{3} \rVert_{\mathcal{H}_{\theta}}^{\frac{2}{p}} \lVert \nabla_{h} \partial_{3} b^{3} \rVert_{\mathcal{H}_{\theta}}^{1-\frac{2}{p}} )\nonumber\\
&+ \lVert b^{3} \rVert_{\dot{H}^{\frac{1}{2} + \frac{2}{p}}}(\lVert \partial_{3}^{2} b^{3} \rVert_{\mathcal{H}_{\theta}} \lVert \partial_{3} u^{3} \rVert_{\mathcal{H}_{\theta}}^{\frac{2}{p}} \lVert \nabla_{h} \partial_{3} u^{3} \rVert_{\mathcal{H}_{\theta}}^{1-\frac{2}{p}} + \lVert \partial_{3}^{2} u^{3} \rVert_{\mathcal{H}_{\theta}} \lVert \partial_{3} b^{3} \rVert_{\mathcal{H}_{\theta}}^{\frac{2}{p}} \lVert \nabla_{h} \partial_{3} b^{3} \rVert_{\mathcal{H}_{\theta}}^{1-\frac{2}{p}} )\nonumber
\end{align}
where we used (95), Lemma 2.7 with $p = q = 2, \alpha = \frac{1}{2}, s = \frac{1}{2} + \frac{2}{p}$ and horizontal Gagliardo-Nirenberg inequality. Hence, in sum of (54), (55) in (53) we obtain 
\begin{align}
&(Q_{3}(u,b) \lvert \partial_{3}u^{3})_{\mathcal{H}_{\theta}} + (R_{2}(u,b) \lvert \partial_{3} b^{3})_{\mathcal{H}_{\theta}}\\
\lesssim& [\left( \lVert \omega_{\frac{3}{4}} \rVert_{L^{2}}^{2} + \lVert d_{\frac{3}{4}} \rVert_{L^{2}}^{2}\right)^{\frac{p+6}{6p}} \left( \lVert \nabla \omega_{\frac{3}{4}} \rVert_{L^{2}}^{2} + \lVert \nabla d_{\frac{3}{4}} \rVert_{L^{2}}^{2} \right)^{\frac{1}{2} (1-\frac{2}{p})}\nonumber \\
& + \left( \lVert \partial_{3} u^{3} \rVert_{\mathcal{H}_{\theta}}^{\frac{2}{p}} + \lVert \partial_{3} b^{3} \rVert_{\mathcal{H}_{\theta}}^{\frac{2}{p}} \right) \left( \lVert \nabla \partial_{3}u^{3} \rVert_{\mathcal{H}_{\theta}}^{1-\frac{2}{p}} + \lVert \nabla \partial_{3} b^{3} \rVert_{\mathcal{H}_{\theta}}^{1-\frac{2}{p}} \right) ]\nonumber\\
&\times \left( \lVert \nabla\partial_{3} u^{3} \rVert_{\mathcal{H}_{\theta}} + \lVert \nabla \partial_{3} b^{3} \rVert_{\mathcal{H}_{\theta}} \right) \left( \lVert u^{3} \rVert_{\dot{H}^{\frac{1}{2} + \frac{2}{p}}} +  \lVert b^{3} \rVert_{\dot{H}^{\frac{1}{2} + \frac{2}{p}}}\right)\nonumber\\
&+ \left( \lVert u^{3} \rVert_{\dot{H}^{\frac{1}{2} + \frac{2}{p}}} + \lVert b^{3} \rVert_{\dot{H}^{\frac{1}{2} + \frac{2}{p}}} \right) \left( \lVert \nabla \partial_{3} u^{3} \rVert_{\mathcal{H}_{\theta}} + \lVert \nabla \partial_{3} b^{3} \rVert_{\mathcal{H}_{\theta}}\right)^{2(1-\frac{1}{p})}\nonumber\\
&\times\left(\lVert \partial_{3} u^{3} \rVert_{\mathcal{H}_{\theta}} + \lVert \partial_{3} b^{3} \rVert_{\mathcal{H}_{\theta}} \right)^{\frac{2}{p}}\nonumber\\
\leq& \frac{1}{3} (\lVert \nabla \partial_{3} u^{3} \rVert_{\mathcal{H}_{\theta}}^{2} + \lVert \nabla \partial_{3} b^{3} \rVert_{\mathcal{H}_{\theta}}^{2}) \nonumber\\
&+ c [\left( \lVert \omega_{\frac{3}{4}} \rVert_{L^{2}}^{2} + \lVert d_{\frac{3}{4}} \rVert_{L^{2}}^{2} \right)^{\frac{p+6}{3p}} \left( \lVert \nabla \omega_{\frac{3}{4}} \rVert_{L^{2}}^{2} + \lVert \nabla d_{\frac{3}{4}} \rVert_{L^{2}}^{2} \right)^{1-\frac{2}{p}} \left( \lVert u^{3} \rVert_{\dot{H}^{\frac{1}{2} + \frac{2}{p}}}^{2} + \lVert b^{3} \rVert_{\dot{H}^{\frac{1}{2} + \frac{2}{p}}}^{2} \right) \nonumber\\
&+ \left( \lVert \partial_{3} u^{3} \rVert_{\mathcal{H}_{\theta}}^{2} + \lVert \partial_{3}b^{3} \rVert_{\mathcal{H}_{\theta}}^{2} \right) \left( \lVert u^{3} \rVert_{\dot{H}^{\frac{1}{2} + \frac{2}{p}}}^{p} + \lVert b^{3} \rVert_{\dot{H}^{\frac{1}{2} + \frac{2}{p}}}^{p}\right) ]\nonumber
\end{align}
by Young's inequalities. Therefore, in sum of (44)-(52), (56) into (42), we obtain 
\begin{align}
&\partial_{t} (\lVert \partial_{3} u^{3} \rVert_{\mathcal{H}_{\theta}}^{2} + \lVert \partial_{3} b^{3} \rVert_{\mathcal{H}_{\theta}}^{2}) + \lVert \nabla \partial_{3} u^{3} \rVert_{\mathcal{H}_{\theta}}^{2} + \lVert \nabla \partial_{3} b^{3} \rVert_{\mathcal{H}_{\theta}}^{2}\\
&-  c (\lVert u^{3} \rVert_{\dot{H}^{\frac{1}{2} + \frac{2}{p}}}^{p} + \lVert b^{3} \rVert_{\dot{H}^{\frac{1}{2} + \frac{2}{p}}}^{p})(\lVert \partial_{3} u^{3} \rVert_{\mathcal{H}_{\theta}}^{2} + \lVert \partial_{3} b^{3} \rVert_{\mathcal{H}_{\theta}}^{2})\nonumber\\
\lesssim& 
 \lVert u^{3} \rVert_{\dot{H}^{\frac{1}{2} + \frac{2}{p}}}(\lVert \omega_{\frac{3}{4}} \rVert_{L^{2}}^{2(\frac{p+3}{3p})}\lVert \nabla \omega_{\frac{3}{4}} \rVert_{L^{2}}^{2(1-\frac{1}{p})} + \lVert d_{\frac{3}{4}} \rVert_{L^{2}}^{2(\frac{p+3}{3p})}\lVert \nabla d_{\frac{3}{4}} \rVert_{L^{2}}^{2(1-\frac{1}{p})})\nonumber\\
&+ (\lVert u^{3} \rVert_{\dot{H}^{\frac{1}{2} + \frac{2}{p}}}^{2} + \lVert b^{3} \rVert_{\dot{H}^{\frac{1}{2} + \frac{2}{p}}}^{2})(\lVert \omega_{\frac{3}{4}} \rVert_{L^{2}}^{2} + \lVert d_{\frac{3}{4}} \rVert_{L^{2}}^{2})^{\frac{p+6}{3p}}(\lVert \nabla \omega_{\frac{3}{4}} \rVert_{L^{2}}^{2} + \lVert \nabla d_{\frac{3}{4}} \rVert_{L^{2}}^{2})^{1-\frac{2}{p}}.\nonumber
\end{align}
Gronwall's type argument using (8) and that $\lVert f\rVert_{\dot{H}^{\frac{1}{2}}} \lesssim \lVert \nabla\times f\rVert_{L^{\frac{3}{2}}}$ for any $f$ such that $\nabla\cdot f= 0$ by Sobolev embedding of $\dot{W}^{\frac{1}{2}, \frac{3}{2}}(\mathbb{R}^{3}) \hookrightarrow L^{2}(\mathbb{R}^{3})$ and continuity of Riesz transform in $L^{\frac{3}{2}}(\mathbb{R}^{3})$ completes the proof of Proposition 3.2.
\end{proof} 

We fix for $T < T^{\ast}$
\begin{equation}
e(T) \triangleq c\exp\left(c\int_{0}^{T}\lVert u^{3} \rVert_{\dot{H}^{\frac{1}{2} + \frac{2}{p}}}^{p} + \lVert b \rVert_{\dot{H}^{\frac{1}{2} + \frac{2}{p}}}^{p} + \lVert   b\rVert_{L^{p_{1}}}^{\frac{2p_{1}}{p_{1} - 3}} + \lVert \nabla b\rVert_{L^{p_{2}}}^{\frac{2p_{2}}{2p_{2} - 3}}d\tau\right).
\end{equation}

\begin{proposition}
Under the hypothesis of Theorem 1.1, for $\theta \in (\frac{1}{2} - \frac{2}{p}, \frac{1}{6})$, the solution to the MHD system (1a)-(1c) satisfies for any $t \leq T$
\begin{align}
&(\lVert \omega_{\frac{3}{4}}  \rVert_{L^{2}}^{2(\frac{p+3}{3})} + \lVert d_{\frac{3}{4}}  \rVert_{L^{2}}^{2(\frac{p+3}{3})})(t) + \lVert \nabla \omega_{\frac{3}{4}} \rVert_{L_{t}^{2} L^{2}}^{2(\frac{p+3}{3})} + \lVert \nabla d_{\frac{3}{4}} \rVert_{L_{t}^{2} L^{2}}^{2(\frac{p+3}{3})}\\
\lesssim& \exp({e(T)})[\lVert \Omega_{0}\rVert_{L^{\frac{3}{2}}}^{\frac{p+3}{2}} +  \lVert j_{0}\rVert_{L^{\frac{3}{2}}}^{\frac{p+3}{2}}]\nonumber
\end{align}
and 
\begin{align}
&(\lVert \partial_{3} u^{3} \rVert_{\mathcal{H}_{\theta}}^{2} + \lVert \partial_{3} b^{3} \rVert_{\mathcal{H}_{\theta}}^{2})(t) + \int_{0}^{t} \lVert \nabla \partial_{3} u^{3} \rVert_{\mathcal{H}_{\theta}}^{2} + \lVert \nabla \partial_{3} b^{3} \rVert_{\mathcal{H}_{\theta}}^{2}d\tau\\
\lesssim& \exp(e(T))[\lVert \Omega_{0} \rVert_{L^{\frac{3}{2}}}^{2} + \lVert j_{0} \rVert_{L^{\frac{3}{2}}}^{2}].\nonumber
\end{align}

\end{proposition}

\begin{proof}
For $t \leq T$ we let 
\begin{align*}
&III_{1}(t) \triangleq \left(\int_{0}^{t} \lVert u^{3} \rVert_{\dot{H}^{\frac{1}{2} + \frac{2}{p}}}(\lVert \omega_{\frac{3}{4}} \rVert_{L^{2}}^{2(\frac{p+3}{3p})}\lVert \nabla \omega_{\frac{3}{4}} \rVert_{L^{2}}^{2(1-\frac{1}{p})} + \lVert d_{\frac{3}{4}} \rVert_{L^{2}}^{2(\frac{p+3}{3p})}\lVert \nabla d_{\frac{3}{4}} \rVert_{L^{2}}^{2(1-\frac{1}{p})}) d\tau  \right)^{\frac{3}{4}},\\
&III_{2}(t)\\
\triangleq& \left(\int_{0}^{t} (\lVert u^{3} \rVert_{\dot{H}^{\frac{1}{2} + \frac{2}{p}}}^{2} + \lVert b^{3} \rVert_{\dot{H}^{\frac{1}{2} + \frac{2}{p}}}^{2})(\lVert \omega_{\frac{3}{4}} \rVert_{L^{2}}^{2} + \lVert d_{\frac{3}{4}} \rVert_{L^{2}}^{2})^{\frac{p+6}{3p}}(\lVert \nabla \omega_{\frac{3}{4}} \rVert_{L^{2}}^{2} + \lVert \nabla d_{\frac{3}{4}} \rVert_{L^{2}}^{2})^{1-\frac{2}{p}}d\tau  \right)^{\frac{3}{4}}.
\end{align*}
By Proposition 3.2 we have 
\begin{align}
&e(T) \left(\int_{0}^{t} \lVert \nabla \partial_{3} u^{3} \rVert_{\mathcal{H}_{\theta}}^{2} + \lVert \nabla \partial_{3} b^{3} \rVert_{\mathcal{H}_{\theta}}^{2} d\tau \right)^{\frac{3}{4}}\\
\leq& e(T) \left(\lVert \Omega_{0} \rVert_{L^{\frac{3}{2}}}^{\frac{3}{2}} + \lVert j_{0} \rVert_{L^{\frac{3}{2}}}^{\frac{3}{2}} + III_{1}(t) + III_{2}(t)\right)\nonumber
\end{align}
as $(a+b)^{\frac{3}{4}} \approx a^{\frac{3}{4}} + b^{\frac{3}{4}}$. We estimate
\begin{align}
e(T) III_{1}(t) \leq& e(T) \left(\int_{0}^{t} \lVert u^{3} \rVert_{\dot{H}^{\frac{1}{2} + \frac{2}{p}}}^{p}\lVert \omega_{\frac{3}{4}} \rVert_{L^{2}}^{2(\frac{p+3}{3})}d\tau  \right)^{\frac{3}{4p}} \left(\int_{0}^{t} \lVert \nabla \omega_{\frac{3}{4}} \rVert_{L^{2}}^{2}d\tau \right)^{\frac{3}{4} (1-\frac{1}{p})}\nonumber\\
&+ e(T) \left(\int_{0}^{t} \lVert u^{3} \rVert_{\dot{H}^{\frac{1}{2} + \frac{2}{p}}}^{p}\lVert d_{\frac{3}{4}} \rVert_{L^{2}}^{2(\frac{p+3}{3})}d\tau \right)^{\frac{3}{4p}}\left(\int_{0}^{t} \lVert \nabla d_{\frac{3}{4}} \rVert_{L^{2}} ^{2} d\tau \right)^{\frac{3}{4}(1-\frac{1}{p})}\nonumber\\
\leq& \frac{1}{9} \int_{0}^{t} \lVert \nabla \omega_{\frac{3}{4}} \rVert_{L^{2}}^{2} + \lVert \nabla d_{\frac{3}{4}} \rVert_{L^{2}}^{2} d\tau\\
&+  e(T) (\int_{0}^{t} \lVert u^{3} \rVert_{\dot{H}^{\frac{1}{2} + \frac{2}{p}}}^{p} \left(\lVert \omega_{\frac{3}{4}} \rVert_{L^{2}}^{2(\frac{p+3}{3})} + 
\lVert d_{\frac{3}{4}} \rVert_{L^{2}}^{2(\frac{p+3}{3})}\right)
d\tau)^{\frac{3}{p+3}} \nonumber
\end{align}
by H$\ddot{o}$lder's and Young's inequalities. Similarly, 
\begin{align}
&e(T) III_{2}(t)\\
\leq& e(T)\left(\int_{0}^{t} (\lVert u^{3} \rVert_{\dot{H}^{\frac{1}{2} + \frac{2}{p}}}^{p} + \lVert b^{3} \rVert_{\dot{H}^{\frac{1}{2} + \frac{2}{p}}}^{p}) (\lVert \omega_{\frac{3}{4}} \rVert_{L^{2}}^{2} + \lVert d_{\frac{3}{4}} \rVert_{L^{2}}^{2} )^{\frac{p+6}{6}}d\tau\right)^{\frac{3}{2p}} \nonumber\\
&\times \left(\int_{0}^{t} \lVert \nabla \omega_{\frac{3}{4}} \rVert_{L^{2}}^{2} + \lVert \nabla d_{\frac{3}{4}} \rVert_{L^{2}}^{} d\tau \right)^{\frac{3}{4}(1-\frac{2}{p})}\nonumber\\
\leq& \frac{1}{9} \int_{0}^{t} \lVert \nabla \omega_{\frac{3}{4}} \rVert_{L^{2}}^{2} + \lVert \nabla d_{\frac{3}{4}}\rVert_{L^{2}}^{2} d\tau\nonumber\\
&+ e(T) \left(\int_{0}^{t} (\lVert u^{3} \rVert_{\dot{H}^{\frac{1}{2} + \frac{2}{p}}}^{p} + \lVert b^{3} \rVert_{\dot{H}^{\frac{1}{2} + \frac{2}{p}}}^{p})(\lVert \omega_{\frac{3}{4}} \rVert_{L^{2}}^{2} + \lVert d_{\frac{3}{4}} \rVert_{L^{2}}^{2})^{\frac{p+6}{6}}d\tau \right)^{\frac{6}{p+6}}\nonumber\\
\leq& \frac{1}{9} \int_{0}^{t} \lVert \nabla \omega_{\frac{3}{4}} \rVert_{L^{2}}^{2} + \lVert \nabla d_{\frac{3}{4}}\rVert_{L^{2}}^{2} d\tau\nonumber\\
&+ e(T) \left(\int_{0}^{t} (\lVert u^{3} \rVert_{\dot{H}^{\frac{1}{2} + \frac{2}{p}}}^{p} + \lVert b^{3} \rVert_{\dot{H}^{\frac{1}{2} + \frac{2}{p}}}^{p}  )(\lVert \omega_{\frac{3}{4}} \rVert_{L^{2}}^{2} + \lVert d_{\frac{3}{4}} \rVert_{L^{2}}^{2})^{\frac{p+3}{3}}d\tau\right)^{\frac{3}{p+3}}\nonumber
\end{align}
where we used H$\ddot{o}$lder's and Young's inequalities and that 
\begin{align*}
& c\exp\left(c_{1}\int_{0}^{T}\lVert u^{3} \rVert_{\dot{H}^{\frac{1}{2} + \frac{2}{p}}}^{p} + \lVert b \rVert_{\dot{H}^{\frac{1}{2} + \frac{2}{p}}}^{p} + \lVert   b\rVert_{L^{p_{1}}}^{\frac{2p_{1}}{p_{1} - 3}} + \lVert \nabla b\rVert_{L^{p_{2}}}^{\frac{2p_{2}}{2p_{2} - 3}}d\tau\right)\\
&\times \left(\int_{0}^{t} \lVert u^{3} \rVert_{\dot{H}^{\frac{1}{2} + \frac{2}{p}}}^{p} + \lVert b^{3} \rVert_{\dot{H}^{\frac{1}{2} + \frac{2}{p}}}^{p}d\tau\right)^{\frac{3p}{(p+3)(p+6)}}\\
\leq& c\exp\left(c_{2}\int_{0}^{T}\lVert u^{3} \rVert_{\dot{H}^{\frac{1}{2} + \frac{2}{p}}}^{p} + \lVert b \rVert_{\dot{H}^{\frac{1}{2} + \frac{2}{p}}}^{p} + \lVert   b\rVert_{L^{p_{1}}}^{\frac{2p_{1}}{p_{1} - 3}} + \lVert \nabla b\rVert_{L^{p_{2}}}^{\frac{2p_{2}}{2p_{2} - 3}}d\tau\right)
\end{align*}
for $c_{2} > c_{1}$ sufficiently large. Therefore, by (62) and (63) applied to (61), 
\begin{align}
&e(T) \left(\int_{0}^{t} \lVert \nabla \partial_{3} u^{3} \rVert_{\mathcal{H}_{\theta}}^{2} + \lVert \nabla \partial_{3} b^{3} \rVert_{\mathcal{H}_{\theta}}^{2} d\tau \right)^{\frac{3}{4}}\\
\leq& \frac{2}{9} \int_{0}^{t} \lVert \nabla \omega_{\frac{3}{4}} \rVert_{L^{2}}^{2} + \lVert \nabla d_{\frac{3}{4}}\rVert_{L^{2}}^{2} d\tau \nonumber\\
&+ e(T) [\lVert \Omega_{0} \rVert_{L^{\frac{3}{2}}}^{\frac{3}{2}} + \lVert j_{0} \rVert_{L^{\frac{3}{2}}}^{\frac{3}{2}}\nonumber\\
& \hspace{15mm} + \left(\int_{0}^{t} (\lVert u^{3} \rVert_{\dot{H}^{\frac{1}{2} + \frac{2}{p}}}^{p} + \lVert b^{3} \rVert_{\dot{H}^{\frac{1}{2} + \frac{2}{p}}}^{p})(\lVert \omega_{\frac{3}{4}} \rVert_{L^{2}}^{2(\frac{p+3}{3})} + \lVert d_{\frac{3}{4}} \rVert_{L^{2}}^{2(\frac{p+3}{3})})d\tau \right)^{\frac{3}{p+3}}
].\nonumber
\end{align}
This leads to 
\begin{align}
&\frac{2}{3} (\lVert \omega_{\frac{3}{4}} \rVert_{L^{2}}^{2} + \lVert d_{\frac{3}{4}} \rVert_{L^{2}}^{2})(t) + \frac{5}{9} \int_{0}^{t} \lVert \nabla \omega_{\frac{3}{4}} \rVert_{L^{2}}^{2} + \lVert \nabla d_{\frac{3}{4}} \rVert_{L^{2}}^{2} d\tau \\
\leq& \frac{2}{9} \int_{0}^{t} \lVert \nabla \omega_{\frac{3}{4}}\rVert_{L^{2}}^{2} + \lVert \nabla d_{\frac{3}{4}}\rVert_{L^{2}}^{2} d\tau\nonumber \\
&+ e(T) [\lVert \Omega_{0} \rVert_{L^{\frac{3}{2}}}^{\frac{3}{2}} + \lVert j_{0} \rVert_{L^{\frac{3}{2}}}^{\frac{3}{2}}\nonumber\\
& \hspace{15mm} + \left(\int_{0}^{t} (\lVert u^{3} \rVert_{\dot{H}^{\frac{1}{2} + \frac{2}{p}}}^{p} + \lVert b^{3} \rVert_{\dot{H}^{\frac{1}{2} + \frac{2}{p}}}^{p})(\lVert \omega_{\frac{3}{4}} \rVert_{L^{2}}^{2(\frac{p+3}{3})} + \lVert d_{\frac{3}{4}} \rVert_{L^{2}}^{2(\frac{p+3}{3})})d\tau \right)^{\frac{3}{p+3}}
]\nonumber
\end{align}
by Proposition 3.1, (64) and that $\lVert \omega_{\frac{3}{4}}(0) \rVert_{L^{2}}^{2} = \int \lvert \omega_{0}\rvert^{\frac{3}{2}} \leq \int \lvert \Omega_{0} \rvert^{\frac{3}{2}} = \lVert \Omega_{0} \rVert_{L^{\frac{3}{2}}}^{\frac{3}{2}}$. After absorbing, we take powers $\frac{p+3}{3}$ and use $(a+b)^{\frac{p+3}{3}} \approx a^{\frac{p+3}{3}} + b^{\frac{p+3}{3}}$ to obtain 
\begin{align*}
&\left(\lVert \omega_{\frac{3}{4}}  \rVert_{L^{2}}^{2(\frac{p+3}{3})} + \lVert d_{\frac{3}{4}}  \rVert_{L^{2}}^{2(\frac{p+3}{3})}\right)(t)  + \lVert \nabla \omega_{\frac{3}{4}} \rVert_{L_{t}^{2} L^{2}}^{2(\frac{p+3}{3})} + \lVert \nabla d_{\frac{3}{4}} \rVert_{L_{t}^{2} L^{2}}^{2(\frac{p+3}{3})}\\
\leq& e(T) [\lVert \Omega_{0}\rVert_{L^{\frac{3}{2}}}^{\frac{p+3}{2}} +  \lVert j_{0}\rVert_{L^{\frac{3}{2}}}^{\frac{p+3}{2}}\\
& \hspace{15mm} + \int_{0}^{t} (\lVert u^{3} \rVert_{\dot{H}^{\frac{1}{2} + \frac{2}{p}}}^{p} + \lVert b^{3} \rVert_{\dot{H}^{\frac{1}{2} + \frac{2}{p}}}^{p})(\lVert \omega_{\frac{3}{4}} \rVert_{L^{2}}^{2(\frac{p+3}{3})} + \lVert d_{\frac{3}{4}} \rVert_{L^{2}}^{2(\frac{p+3}{3})})d\tau ].
\end{align*}
Thus, Gronwall's type argument using again $e(T) \int_{0}^{t} \lVert u^{3} \rVert_{\dot{H}^{\frac{1}{2} + \frac{2}{p}}}^{p} + \lVert b^{3} \rVert_{\dot{H}^{\frac{1}{2} + \frac{2}{p}}}^{p} d\tau \lesssim e(T)$ by changing the constant in the exponent leads to (59). 

Next, by our Proposition 3.2 
\begin{align}
&(\lVert \partial_{3} u^{3} \rVert_{\mathcal{H}_{\theta}}^{2} + \lVert \partial_{3} b^{3} \rVert_{\mathcal{H}_{\theta}}^{2})(t) + \int_{0}^{t} \lVert \nabla \partial_{3} u^{3} \rVert_{\mathcal{H}_{\theta}}^{2} + \lVert \nabla \partial_{3} b^{3} \rVert_{\mathcal{H}_{\theta}}^{2}d\tau\\
\leq& e(T) [ \lVert \Omega_{0} \rVert_{L^{\frac{3}{2}}}^{2} + \lVert j_{0} \rVert_{L^{\frac{3}{2}}}^{2}\nonumber\\
&+ \left(\int_{0}^{t} \lVert u^{3} \rVert_{\dot{H}^{\frac{1}{2} + \frac{2}{p}}}^{p} d\tau \right)^{\frac{1}{p}} \sup_{\tau \in [0,t]} \left( \lVert \omega_{\frac{3}{4}} \rVert_{L^{2}}^{2(\frac{p+3}{3p})} + \lVert d_{\frac{3}{4}} \rVert_{L^{2}}^{2(\frac{p+3}{3p})}\right)(\tau)\nonumber\\
&\times \left(\int_{0}^{t} \lVert \nabla \omega_{\frac{3}{4}} \rVert_{L^{2}}^{2} + \lVert \nabla d_{\frac{3}{4}} \rVert_{L^{2}}^{2} d\tau \right)^{1-\frac{1}{p}}\nonumber\\
&+ \left(\int_{0}^{t} \lVert u^{3} \rVert_{\dot{H}^{\frac{1}{2} + \frac{2}{p}}}^{p} + \lVert b^{3} \rVert_{\dot{H}^{\frac{1}{2} + \frac{2}{p}}}^{p}d\tau \right)^{\frac{2}{p}}\sup_{\tau \in [0,t]} \left( \lVert \omega_{\frac{3}{4}} \rVert_{L^{2}}^{2} + \lVert d_{\frac{3}{4}} \rVert_{L^{2}}^{2} \right)^{\frac{p+6}{3p}}(\tau)\nonumber\\
&\times\left(\int_{0}^{t} \lVert \nabla \omega_{\frac{3}{4}} \rVert_{L^{2}}^{2} + \lVert \nabla d_{\frac{3}{4}} \rVert_{L^{2}}^{2} d\tau \right)^{1-\frac{2}{p}}]\nonumber\\
\leq& 
e(T) [\lVert \Omega_{0} \rVert_{L^{\frac{3}{2}}}^{2} + \lVert j_{0} \rVert_{L^{\frac{3}{2}}}^{2}\nonumber\\
&+ \lVert u^{3} \rVert_{L_{t}^{p} \dot{H}^{\frac{1}{2} + \frac{2}{p}}} \left(e^{e(T)}[\lVert \Omega_{0}\rVert_{L^{\frac{3}{2}}}^{(\frac{p+3}{2})} +  \lVert j_{0}\rVert_{L^{\frac{3}{2}}}^{(\frac{p+3}{2})}]\right)^{\frac{1}{p} + \frac{3}{p}(\frac{p-1}{p+3})}\nonumber\\
&+ \left(\lVert u^{3} \rVert_{L_{t}^{p} \dot{H}^{\frac{1}{2} + \frac{2}{p}}} + \lVert b^{3} \rVert_{L_{t}^{p} \dot{H}^{\frac{1}{2} + \frac{2}{p}}}\right)^{2} 
\left(e^{e(T)}[\lVert \Omega_{0}\rVert_{L^{\frac{3}{2}}}^{(\frac{p+3}{2})} +  \lVert j_{0}\rVert_{L^{\frac{3}{2}}}^{(\frac{p+3}{2})}]\right)^{\frac{p+6}{p(p+3)} + \frac{3}{p}(\frac{p-2}{p+3})}]\nonumber\\
\lesssim& \exp(e(T))[\lVert \Omega_{0} \rVert_{L^{\frac{3}{2}}}^{2} + \lVert j_{0} \rVert_{L^{\frac{3}{2}}}^{2}]\nonumber
\end{align}
by H$\ddot{o}$lder's inequalities and (59). This completes the proof of Proposition 3.3.
\end{proof}

\section{Blow-up criterion}

\begin{proposition}

Let $u, b\in C(0, T^{\ast}; \dot{H}^{\frac{1}{2}}(\mathbb{R}^{3})) \cap L^{2}(0, T^{\ast}; \dot{H}^{\frac{3}{2}}(\mathbb{R}^{3}))$ solve the MHD system (1a)-(1c). If $T^{\ast} < \infty$, then for any $p_{k,l} \in (1, \infty), k, l \in \{1, 2, 3\}$, 
\begin{equation*}
\sum_{k, l = 1}^{3} \int_{0}^{T^{\ast}} \lVert \partial_{l} u^{k}  \rVert_{\mathcal{B}_{p_{k,l}}}^{p_{k,l}} + \lVert \partial_{l} b^{k}  \rVert_{\mathcal{B}_{p_{k,l}}}^{p_{k,l}}
 d\tau = \infty.
\end{equation*}

\end{proposition}

\begin{proof}
We apply $\dot{\Delta}_{j}$ on (1a)-(1b), take $L^{2}$-inner products with $\dot{\Delta}_{j} u, \dot{\Delta}_{j} b$ respectively, sum the two equations, multiplying by $2^{j}$ and sum over $j \in \mathbb{Z}$ to obtain 
\begin{align}
&\frac{1}{2}\partial_{t} (\lVert u\rVert_{\dot{H}^{\frac{1}{2}}}^{2} + \lVert b\rVert_{\dot{H}^{\frac{1}{2}}}^{2}) + \lVert \nabla u\rVert_{\dot{H}^{\frac{1}{2}}}^{2} + \lVert \nabla b\rVert_{\dot{H}^{\frac{1}{2}}}^{2}\\
=& -(u\cdot\nabla u \lvert u)_{\dot{H}^{\frac{1}{2}}}  +(b\cdot\nabla b \lvert u)_{\dot{H}^{\frac{1}{2}}} -(u\cdot\nabla b \lvert b)_{\dot{H}^{\frac{1}{2}}} +(b\cdot\nabla u \lvert b)_{\dot{H}^{\frac{1}{2}}}.\nonumber
\end{align}
Now we show that  
\begin{align}
&\lvert (u\cdot\nabla u\lvert u)_{\dot{H}^{\frac{1}{2}}} + (u\cdot\nabla b\lvert b)_{\dot{H}^{\frac{1}{2}}} - (b\cdot\nabla b\lvert u)_{\dot{H}^{\frac{1}{2}}} - (b\cdot\nabla u \lvert b)_{\dot{H}^{\frac{1}{2}}}\rvert \\
\lesssim& \sum_{k,l=1}^{3} (\lVert \partial_{l} u^{k} \rVert_{\mathcal{B}_{p_{k,l}}} + \lVert \partial_{l} b^{k} \rVert_{\mathcal{B}_{p_{k,l}}})(\lVert u\rVert_{\dot{H}^{\frac{1}{2}}}^{\frac{2}{p_{k,l}}} + \lVert b\rVert_{\dot{H}^{\frac{1}{2}}}^{\frac{2}{p_{k,l}}})(\lVert \nabla u\rVert_{\dot{H}^{\frac{1}{2}}}^{2(1-\frac{1}{p_{k,l}})} + \lVert \nabla b\rVert_{\dot{H}^{\frac{1}{2}}}^{2(1-\frac{1}{p_{k,l}})}).\nonumber
\end{align}

Here, making use of the structure of the MHD system becomes important as otherwise, the proof leads us to a non-favorable condition of $p_{l,k} > 2$. 

By Lemma 8.1 of [12], we already have 
\begin{equation}
\lvert (u\cdot\nabla u\lvert u)_{\dot{H}^{\frac{1}{2}}}\rvert \lesssim \sum_{k,l=1}^{3} \lVert \partial_{l} u^{k} \rVert_{\mathcal{B}_{p_{k,l}}}\lVert u\rVert_{\dot{H}^{\frac{1}{2}}}^{\frac{2}{p_{k,l}}}\lVert \nabla u\rVert_{\dot{H}^{\frac{1}{2}}}^{2(1-\frac{1}{p_{k,l}})}.
\end{equation}
Firstly we work on 
\begin{align}
& \lvert (u^{l} \partial_{l} b^{k} \lvert b^{k} )_{\dot{H}^{\frac{1}{2}}}\rvert\\
\leq& \sum_{j} 2^{j} \lvert (\dot{\Delta}_{j} T(u^{l}, \partial_{l} b^{k}) \lvert \dot{\Delta}_{j} b^{k}) \rvert + \sum_{j} 2^{j} \lvert (\dot{\Delta}_{j} T(\partial_{l} b^{k}, u^{l}) \lvert \dot{\Delta}_{j} b^{k}) \rvert\nonumber\\
&+ \sum_{j} 2^{j} \lvert (\dot{\Delta}_{j} R(u^{l}, \partial_{l} b^{k}) \lvert \dot{\Delta}_{j} b^{k}) \rvert \triangleq IV_{1} + IV_{2} + IV_{3}\nonumber
\end{align}
due to Bony's paraproduct decomposition (6). We start with 
\begin{align}
IV_{1} 
\leq & \sum_{j} 2^{j}\lvert  \int \dot{S}_{j-1} u^{l} \dot{\Delta}_{j} \partial_{l} b^{k} \dot{\Delta}_{j} b^{k}\rvert  + \lvert \int \sum_{\lvert j-j'\rvert \leq 4} [\dot{\Delta}_{j}, \dot{S}_{j'-1} u^{l} ] \dot{\Delta}_{j'} \partial_{l} b^{k} \dot{\Delta}_{j} b^{k}\rvert\\
&+ \lvert \int  \sum_{\lvert j-j'\rvert \leq 4} (\dot{S}_{j'-1} u^{l} - \dot{S}_{j-1} u^{l} ) \dot{\Delta}_{j} \dot{\Delta}_{j'} \partial_{l} b^{k} \dot{\Delta}_{j} b^{k}\rvert  \triangleq IV_{1,1} + IV_{1,2} + IV_{1,3}.\nonumber
\end{align}
We have due to divergence-free property, $IV_{1,1} =0$. Secondly, 
\begin{align}
IV_{1,2} 
\leq& \sum_{j} 2^{j} \sum_{\lvert j-j'\rvert \leq 4}\lVert [ \dot{\Delta}_{j}, \dot{S}_{j-1} u^{l} ] \dot{\Delta}_{j'} \partial_{l} b^{k} \rVert_{L^{2}} \lVert \dot{\Delta}_{j} b^{k} \rVert_{L^{2}}\\
\lesssim& \sum_{j} 2^{j} 2^{-j} \sum_{\lvert j-j'\rvert \leq 4}\lVert \nabla \dot{S}_{j'-1} u^{l} \rVert_{L^{\infty}} \lVert \dot{\Delta}_{j'} \partial_{l} b^{k} \rVert_{L^{2}} \lVert \dot{\Delta}_{j} b^{k} \rVert_{L^{2}} \nonumber\\
\lesssim& \sum_{j} \sum_{l'=1}^{3} \lVert \partial_{l'} \dot{S}_{j-1} u^{l} \rVert_{L^{\infty}} \lVert \dot{\Delta}_{j} \partial_{l} b^{k} \rVert_{L^{2}} \lVert \dot{\Delta}_{j} b^{k} \rVert_{L^{2}}\nonumber\\
\lesssim& \sum_{l'=1}^{3} \sum_{j} 2^{j(2-\frac{2}{p_{l,l'} })}\nonumber\\
&\times\sum_{j' \leq j-2}2^{(j'-j)(2-\frac{2}{p_{l,l'} })}2^{j'(-2 + \frac{2}{p_{l,l'} })}\lVert \dot{\Delta}_{j'} \partial_{l'} u^{l} \rVert_{L^{\infty}}\lVert \dot{\Delta}_{j} \partial_{l} b^{k} \rVert_{L^{2}} \lVert \dot{\Delta}_{j} b^{k} \rVert_{L^{2}}\nonumber\\
\lesssim& \sum_{l'=1}^{3} \lVert \partial_{l'} u^{l} \rVert_{\mathcal{B}_{p_{l,l'} }}\sum_{j} (2^{-\frac{j}{2}} \lVert \dot{\Delta}_{j} \partial_{l} b^{k} \rVert_{L^{2}})^{\frac{1}{p_{l,l'} }}(2^{\frac{j}{2}} \lVert \dot{\Delta}_{j} b^{k} \rVert_{L^{2}})^{\frac{1}{p_{l,l'} }}\nonumber\\
&\hspace{25mm} \times(2^{\frac{j}{2}} \lVert \dot{\Delta}_{j} \partial_{l} b^{k} \rVert_{L^{2}})^{1-\frac{1}{p_{l,l'} }}(2^{\frac{3j}{2}} \lVert \dot{\Delta}_{j} b^{k} \rVert_{L^{2}})^{1-\frac{1}{p_{l,l'} }}\nonumber\\
\lesssim& \sum_{l'=1}^{3} \lVert \partial_{l'} u^{l} \rVert_{\mathcal{B}_{p_{l,l'} }}\lVert b\rVert_{\dot{H}^{\frac{1}{2}}}^{\frac{2}{p_{l,l'} }}\lVert \nabla b\rVert_{\dot{H}^{\frac{1}{2}}}^{2(1-\frac{1}{p_{l,l'} })}\nonumber
\end{align}
where we used H$\ddot{o}$lder's inequality, a commutator estimate (cf. Lemma 2.97 [1] and also [28]) and Young's inequality for convolution. Thirdly, 
\begin{align}
IV_{1,3} \leq& \sum_{j} 2^{j} \sum_{\lvert j-j'\rvert \leq 4} \lVert (\dot{S}_{j'-1} u^{l} - \dot{S}_{j-1} u^{l} ) \dot{\Delta}_{j} \dot{\Delta}_{j'} \partial_{l} b^{k} \rVert_{L^{2}} \lVert \dot{\Delta}_{j} b^{k} \rVert_{L^{2}}\\
\lesssim& \sum_{j} 2^{j} \sum_{\lvert j-j'\rvert \leq 4, j'' \in [j-1, j'-1]}\lVert \sum_{l'=1}^{3} 2^{-j''}\tilde{\Delta}_{j''}^{l'} \dot{\Delta}_{j''}\partial_{l'} u^{l} \rVert_{L^{\infty}} \lVert \dot{\Delta}_{j} \dot{\Delta}_{j'} \partial_{l} b^{k} \rVert_{L^{2}} \lVert \dot{\Delta}_{j} b^{k} \rVert_{L^{2}}\nonumber\\
\lesssim& \sum_{j} 2^{j} \sum_{l'=1}^{3} 2^{-j} \lVert \dot{\Delta}_{j} \partial_{l'} u^{l} \rVert_{L^{\infty}} \lVert \dot{\Delta}_{j} \partial_{l} b^{k} \rVert_{L^{2}} \lVert \dot{\Delta}_{j} b^{k} \rVert_{L^{2}}\nonumber\\
\lesssim& \sum_{l'=1}^{3} \lVert \partial_{l'} u^{l} \rVert_{\mathcal{B}_{p_{l,l'} }}\sum_{j} (2^{-\frac{j}{2}} \lVert \dot{\Delta}_{j} \partial_{l} b^{k} \rVert_{L^{2}})^{\frac{1}{p_{l,l'} }}(2^{\frac{j}{2}} \lVert \dot{\Delta}_{j} b^{k} \rVert_{L^{2}})^{\frac{1}{p_{l,l'} }}\nonumber\\
&\hspace{20mm} \times (2^{\frac{j}{2}} \lVert \dot{\Delta}_{j} \partial_{l} b^{k} \rVert_{L^{2}})^{1-\frac{1}{p_{l,l'} }}(2^{\frac{3j}{2}} \lVert \dot{\Delta}_{j} b^{k} \rVert_{L^{2}})^{1-\frac{1}{p_{l,l'} }}\nonumber\\
\lesssim& \sum_{l'=1}^{3} \lVert \partial_{l'} u^{l} \rVert_{\mathcal{B}_{p_{l, l'}}}\lVert b\rVert_{\dot{H}^{\frac{1}{2}}}^{\frac{2}{p_{l,l'} }}\lVert \nabla b\rVert_{\dot{H}^{\frac{1}{2}}}^{2(1-\frac{1}{p_{l,l'} })}\nonumber
\end{align}
by H$\ddot{o}$lder's inequality; we also used the fact that $\lvert j-j'\rvert \leq 4$ and $j'' \in [j-1, j'-1]$ imply that we can assume these indices are all $j$ modifying constants.  Therefore, (72), (73) in (71) imply
\begin{equation}
IV_{1} \lesssim \sum_{l'=1}^{3} \lVert \partial_{l'} u^{l} \rVert_{\mathcal{B}_{p_{l, l'}}}\lVert b\rVert_{\dot{H}^{\frac{1}{2}}}^{\frac{2}{p_{l,l'} }}\lVert \nabla b\rVert_{\dot{H}^{\frac{1}{2}}}^{2(1-\frac{1}{p_{l,l'} })}.
\end{equation}
Next, 
\begin{align}
IV_{2} \lesssim& \sum_{j} \sum_{\lvert j- j'\rvert \leq 4} 2^{j} \lVert \dot{S}_{j'-1} \partial_{l} b^{k} \rVert_{L^{\infty}} \lVert \dot{\Delta}_{j'} u^{l} \rVert_{L^{2}} \lVert \dot{\Delta}_{j} b^{k} \rVert_{L^{2}}\\
\lesssim& \sum_{j} 2^{j} \lVert \dot{S}_{j-1} \partial_{l} b^{k} \rVert_{L^{\infty}} \lVert \dot{\Delta}_{j} u^{l} \rVert_{L^{2}} \lVert \dot{\Delta}_{j} b^{k} \rVert_{L^{2}}\nonumber\\
\lesssim&\sum_{j} 2^{j(3-\frac{2}{p_{k,l}})}\sum_{j'\leq j-2} 2^{(j'-j)(2-\frac{2}{p_{k,l}})}2^{j'(-2 + \frac{2}{p_{k,l}})}\lVert \dot{\Delta}_{j'} \partial_{l} b^{k} \rVert_{L^{\infty}} \lVert \dot{\Delta}_{j} u^{l} \rVert_{L^{2}} \lVert \dot{\Delta}_{j} b^{k} \rVert_{L^{2}}\nonumber\\
\lesssim& \lVert \partial_{l} b^{k} \rVert_{\mathcal{B}_{p_{k,l}}}\sum_{j} (2^{\frac{j}{2}} \lVert \dot{\Delta}_{j} u^{l} \rVert_{L^{2}})^{\frac{1}{p_{k,l}}} (2^{\frac{j}{2}} \lVert \dot{\Delta}_{j} b^{k} \rVert_{L^{2}})^{\frac{1}{p_{k,l}}}\nonumber\\
&\hspace{20mm} \times(2^{\frac{3j}{2}} \lVert \dot{\Delta}_{j} u^{l} \rVert_{L^{2}})^{1-\frac{1}{p_{k,l}}} (2^{\frac{3j}{2}} \lVert \dot{\Delta}_{j} b^{k} \rVert_{L^{2}})^{1-\frac{1}{p_{k,l}}}\nonumber\\
\lesssim& \lVert \partial_{l} b^{k} \rVert_{\mathcal{B}_{p_{k,l}}}(\lVert u\rVert_{\dot{H}^{\frac{1}{2}}}^{\frac{2}{p_{k,l}}} + \lVert b\rVert_{\dot{H}^{\frac{1}{2}}}^{\frac{2}{p_{k,l}}})(\lVert \nabla u\rVert_{\dot{H}^{\frac{1}{2}}}^{2(1-\frac{1}{p_{k,l}})} + \lVert \nabla b\rVert_{\dot{H}^{\frac{1}{2}}}^{2(1-\frac{1}{p_{k,l}})})\nonumber
\end{align}
where we used H$\ddot{o}$lder's and Young's inequality for convolution. Finally, we first write by divergence-free condition, 
\begin{align*}
IV_{3} = \sum_{j} 2^{j} \lvert (\dot{\Delta}_{j} R(u^{l}, \partial_{l} b^{k}) \lvert \dot{\Delta}_{j} b^{k}) \rvert= \sum_{j} 2^{j} \lvert (\dot{\Delta}_{j} \sum_{j' \geq j-\delta}\partial_{l}(\dot{\Delta}_{j'} u^{l} \tilde{\dot{\Delta}}_{j'} b^{k}) \lvert \dot{\Delta}_{j} b^{k} ) \rvert 
\end{align*}
for some $\delta \in \mathbb{Z}^{+}0$ and $
\dot{\Delta}_{j'} u^{l} \approx \sum_{l'=1}^{3} \tilde{\dot{\Delta}}_{j'}^{l'} \dot{\Delta}_{j'} u^{l} \approx \sum_{l'=1}^{3} 2^{-j'} \tilde{\dot{\Delta}}_{j'}^{l'} \dot{\Delta}_{j'} \partial_{l'} u^{l}$ so that 
\begin{align}
IV_{3} 
\lesssim& \sum_{j} 2^{j} \lVert \dot{\Delta}_{j} \sum_{j' \geq j-\delta}\partial_{l} (\sum_{l'=1}^{3} 2^{-j'} \tilde{\dot{\Delta}}_{j'}^{l'} \dot{\Delta}_{j'} \partial_{l'} u^{l} \tilde{\dot{\Delta}}_{j'} b^{k}) \rVert_{L^{2}} \lVert \dot{\Delta}_{j} b^{k} \rVert_{L^{2}}\\
\lesssim& \sum_{j} 2^{j} \sum_{j' \geq j - \delta} \sum_{l'=1}^{3} 2^{j-j'} \lVert \tilde{\dot{\Delta}}_{j'}^{l'} \dot{\Delta}_{j'} \partial_{l'} u^{l} \rVert_{L^{\infty}} \lVert \tilde{\dot{\Delta}}_{j'} b^{k} \rVert_{L^{2}} \lVert \dot{\Delta}_{j} b^{k} \rVert_{L^{2}}\nonumber\\
\lesssim& \sum_{l'=1}^{3} \lVert \partial_{l'} u^{l} \rVert_{\mathcal{B}_{p_{l,l'} }}\sum_{j} \sum_{j' \geq j-\delta} 2^{(j-j')(\frac{1}{2} + \frac{1}{p_{l,l'} })}\nonumber\\
&\times (2^{\frac{j'}{2}}\lVert \tilde{\dot{\Delta}}_{j'} b\rVert_{L^{2}})^{\frac{1}{p_{l,l'} }}(2^{\frac{3j'}{2}} \lVert \tilde{\dot{\Delta}}_{j'} b\rVert_{L^{2}})^{1-\frac{1}{p_{l,l'} }} (2^{\frac{j}{2}}\lVert \dot{\Delta}_{j} b\rVert_{L^{2}})^{\frac{1}{p_{l,l'} }}(2^{\frac{3j}{2}} \lVert \dot{\Delta}_{j} b\rVert_{L^{2}})^{1-\frac{1}{p_{l,l'} }}\nonumber\\
\lesssim& \sum_{l'=1}^{3} \lVert \partial_{l'} u^{l} \rVert_{\mathcal{B}_{p_{l,l'} }} \lVert (2^{\frac{j}{2}} \lVert \tilde{\dot{\Delta}}_{j} b\rVert_{L^{2}})^{\frac{1}{p_{l,l'} }}(2^{\frac{3j}{2}} \lVert \tilde{\dot{\Delta}}_{j} b\rVert_{L^{2}})^{1-\frac{1}{p_{l,l'} }}\rVert_{l^{2}}\nonumber\\
&\times \lVert (2^{\frac{j}{2}} \lVert \dot{\Delta}_{j} b\rVert_{L^{2}})^{\frac{1}{p_{l,l'} }}(2^{\frac{3j}{2}} \lVert \dot{\Delta}_{j} b\rVert_{L^{2}})^{1-\frac{1}{p_{l,l'} }}\rVert_{l^{2}}\nonumber\\
\lesssim& \sum_{l'=1}^{3} \lVert \partial_{l'} u^{l} \rVert_{\mathcal{B}_{p_{l,l'} }}\lVert b\rVert_{\dot{H}^{\frac{1}{2}}}^{\frac{2}{p_{l,l'} }}\lVert \nabla b\rVert_{\dot{H}^{\frac{1}{2}}}^{2(1-\frac{1}{p_{l,l'} })}\nonumber
\end{align}
by H$\ddot{o}$lder's, Bernstein's and Young's inequality for convolution. Thus, from (74)-(76) applied to (70) 
\begin{align}
&\lvert (u^{l} \partial_{l} b^{k} \lvert b^{k})_{\dot{H}^{\frac{1}{2}}}\rvert\\
\lesssim& \sum_{k,l=1}^{3} (\lVert \partial_{l} u^{k} \rVert_{\mathcal{B}_{p_{k,l}}} + \lVert \partial_{l} b^{k} \rVert_{\mathcal{B}_{p_{k,l}}}) ( \lVert u\rVert_{\dot{H}^{\frac{1}{2}}}^{\frac{2}{p_{k,l}}} + \lVert b\rVert_{\dot{H}^{\frac{1}{2}}}^{\frac{2}{p_{k,l}}})(\lVert \nabla u\rVert_{\dot{H}^{\frac{1}{2}}}^{2(1-\frac{1}{p_{k,l}})} + \lVert \nabla b\rVert_{\dot{H}^{\frac{1}{2}}}^{2(1-\frac{1}{p_{k,l}})}).\nonumber
\end{align}
Next, we work on 
\begin{align}
(b^{l}\partial_{l} b^{k} \lvert u^{k})_{\dot{H}^{\frac{1}{2}}}
=& \sum_{j} 2^{j} ( \dot{\Delta}_{j} T(b^{l}, \partial_{l} b^{k}) \lvert \dot{\Delta}_{j} u^{k}) + \sum_{j} 2^{j}(\dot{\Delta}_{j} T(\partial_{l} b^{k}, b^{l}) \lvert \dot{\Delta}_{j} u^{k})\\
&\hspace{30mm} + \sum_{j} 2^{j} (\dot{\Delta}_{j} R(b^{l}, \partial_{l}b^{k}) \lvert \dot{\Delta}_{j} u^{k}) \triangleq V_{1} + V_{2} + V_{3}\nonumber
\end{align}
due to Bony's paraproduct decomposition (6). Let us work on $V_{1}$ subsequently together with $VI_{1}$ to be defined below. We now estimate similarly to $IV_{2}$ in (75) 
\begin{align}
V_{2} 
\lesssim& \sum_{j} 2^{j} \lVert \dot{S}_{j-1} \partial_{l} b^{k} \rVert_{L^{\infty}} \lVert \dot{\Delta}_{j} b^{l} \rVert_{L^{2}} \lVert \dot{\Delta}_{j} u^{k} \rVert_{L^{2}}\\
\lesssim& \lVert \partial_{l} b^{k} \rVert_{\mathcal{B}_{p_{k,l}}}\sum_{j}(2^{\frac{j}{2}}\lVert \dot{\Delta}_{j} b^{l} \rVert_{L^{2}})^{\frac{1}{p_{k,l}}}(2^{\frac{j}{2}}\lVert \dot{\Delta}_{j} u^{k} \rVert_{L^{2}})^{\frac{1}{p_{k,l}}}
\nonumber\\
& \hspace{20mm} \times(2^{\frac{3j}{2}}\lVert \dot{\Delta}_{j} b^{l} \rVert_{L^{2}})^{1-\frac{1}{p_{k,l}}}(2^{\frac{3j}{2}}\lVert \dot{\Delta}_{j} u^{k} \rVert_{L^{2}})^{1-\frac{1}{p_{k,l}}}\nonumber\\
\lesssim& \lVert \partial_{l} b^{k} \rVert_{\mathcal{B}_{p_{k,l}}}(\lVert b\rVert_{\dot{H}^{\frac{1}{2}}}^{\frac{2}{p_{k,l}}} + \lVert u\rVert_{\dot{H}^{\frac{1}{2}}}^{\frac{2}{p_{k,l}}})(\lVert \nabla b\rVert_{\dot{H}^{\frac{1}{2}}}^{2(1-\frac{1}{p_{k,l}})} + \lVert \nabla u\rVert_{\dot{H}^{\frac{1}{2}}}^{2(1-\frac{1}{p_{k,l}})})\nonumber
\end{align}
by H$\ddot{o}$lder's and Young's inequality for convolution. Next, as done in (76), by divergence-free condition, and writing $
\dot{\Delta}_{j'} b^{l} \approx \sum_{l'=1}^{3} \tilde{\dot{\Delta}}_{j'}^{l'} \dot{\Delta}_{j'} b^{l} \approx \sum_{l'=1}^{3} 2^{-j'} \tilde{\dot{\Delta}}_{j'}^{l'} \dot{\Delta}_{j'} \partial_{l'} b^{l}$, we estimate 
\begin{align}
V_{3} \approx& \sum_{j}2^{j} (\dot{\Delta}_{j} \sum_{j' \geq j-\delta} \partial_{l} ((\sum_{l'=1}^{3} 2^{-j'} \tilde{\dot{\Delta}}_{j'}^{l'} \dot{\Delta}_{j'} \partial_{l'} b^{l}) \tilde{\dot{\Delta}}_{j'} b^{k}) \lvert \dot{\Delta}_{j} b^{k})\\
\lesssim& \sum_{j} 2^{j} \sum_{j' \geq j-\delta} \sum_{l'=1}^{3} 2^{j-j'}\lVert \tilde{\dot{\Delta}}_{j'}^{l'} \dot{\Delta}_{j'} \partial_{l'} b^{l} \rVert_{L^{\infty}} \lVert \tilde{\dot{\Delta}}_{j'} b^{k} \rVert_{L^{2}} \lVert \dot{\Delta}_{j} b^{k} \rVert_{L^{2}}\nonumber\\
\lesssim& \sum_{l'=1}^{3} \lVert \partial_{l'} b^{l} \rVert_{\mathcal{B}_{p_{l,l'} }}\sum_{j} \sum_{j' \geq j-\delta} 2^{(j-j')(\frac{1}{2} + \frac{1}{p_{l,l'} })}\nonumber\\
&\times (2^{\frac{j'}{2}} \lVert \tilde{\dot{\Delta}}_{j'} b^{k}\rVert_{L^{2}})^{\frac{1}{p_{l,l'} }} (2^{\frac{3j'}{2}} \lVert \tilde{\dot{\Delta}}_{j'} b^{k}\rVert_{L^{2}})^{1-\frac{1}{p_{l,l'} }} (2^{\frac{j}{2}} \lVert \dot{\Delta}_{j} b^{k}\rVert_{L^{2}})^{\frac{1}{p_{l,l'} }} (2^{\frac{3j}{2}} \lVert \dot{\Delta}_{j} b^{k}\rVert_{L^{2}})^{1-\frac{1}{p_{l,l'} }}\nonumber\\
\lesssim& \sum_{l'=1}^{3} \lVert \partial_{l'} b^{l} \rVert_{\mathcal{B}_{p_{l,l'} }}\lVert (2^{\frac{j}{2}} \lVert \dot{\Delta}_{j} b^{k}\rVert_{L^{2}})^{\frac{1}{p_{l,l'} }} (2^{\frac{3j}{2}} \lVert \dot{\Delta}_{j} b^{k}\rVert_{L^{2}})^{1-\frac{1}{p_{l,l'} }}\rVert_{l^{2}}^{2}\nonumber\\
\lesssim& \sum_{l'=1}^{3} \lVert \partial_{l'} b^{l} \rVert_{\mathcal{B}_{p_{l,l'} }}\lVert b\rVert_{\dot{H}^{\frac{1}{2}}}^{\frac{2}{p_{l,l'} }} \lVert \nabla b\rVert_{\dot{H}^{\frac{1}{2}}}^{2(1-\frac{1}{p_{l,l'} })}\nonumber
\end{align}
by H$\ddot{o}$lder's, Bernstein's and Young's inequality for convolution. 

Finally, we consider 
\begin{align}
(b^{l} \partial_{l} u^{k} \lvert b^{k})_{\dot{H}^{\frac{1}{2}}}=& \sum_{j} 2^{j} (\dot{\Delta}_{j} T(b^{l}, \partial_{l} u^{k}) \lvert \dot{\Delta}_{j} b^{k}) + \sum_{j}2^{j}(\dot{\Delta}_{j} T(\partial_{l} u^{k}, b^{l}) \lvert \dot{\Delta}_{j} b^{k})\\
&+  \sum_{j}2^{j}(\dot{\Delta}_{j} R(b^{l}, \partial_{l} u^{k}) \lvert \dot{\Delta}_{j} b^{k})
\triangleq VI_{1} + VI_{2} + VI_{3}.\nonumber
\end{align}
We now consider $V_{1}$ from (78) along with $VI_{1}$ of (81): 
\begin{align}
V_{1} + VI_{1} 
=& \sum_{j} 2^{j} \int \dot{S}_{j-1} b^{l} \dot{\Delta}_{j} \partial_{l} b^{k} \dot{\Delta}_{j} u^{k} + \sum_{\lvert j-j'\rvert \leq 4} [\dot{\Delta}_{j}, \dot{S}_{j'-1} b^{l}] \dot{\Delta}_{j'}\partial_{l}b^{k} \dot{\Delta}_{j} u^{k}\\
&+ \sum_{\lvert j-j'\rvert \leq 4} (\dot{S}_{j'-1} b^{l} - \dot{S}_{j-1} b^{l}) \dot{\Delta}_{j} \dot{\Delta}_{j'} \partial_{l} b^{k} \dot{\Delta}_{j} u^{k} + \dot{S}_{j-1} b^{l} \dot{\Delta}_{j} \partial_{l} u^{k} \dot{\Delta}_{j} b^{k}\nonumber\\
&+ \sum_{\lvert j-j'\rvert \leq 4} [\dot{\Delta}_{j}, \dot{S}_{j'-1} b^{l}] \dot{\Delta}_{j'} \partial_{l} u^{k} \dot{\Delta}_{j} b^{k}\nonumber\\
&+ \sum_{\lvert j-j'\rvert \leq 4} (\dot{S}_{j'-1} b^{l} - \dot{S}_{j-1} b^{l}) \dot{\Delta}_{j} \dot{\Delta}_{j'} \partial_{l} u^{k} \dot{\Delta}_{j} b^{k} \triangleq \sum_{i=1}^{6} (V_{1} + VI_{1})_{i}.\nonumber
\end{align}
We make use of that, together due to the divergence-free property of $b$, $
(V_{1} + VI_{1})_{1} + (V_{1} + VI_{1})_{4} = 0$. Now we work similarly to (72) on 
\begin{align}
&(V_{1} + VI_{1})_{2} + (V_{1} + VI_{1})_{5} \\
\lesssim& \sum_{j} 2^{j} 2^{-j} \sum_{\lvert j-j'\rvert \leq 4} \lVert \nabla \dot{S}_{j'-1} b^{l} \rVert_{L^{\infty}} (\lVert \dot{\Delta}_{j'} \partial_{l} b^{k} \rVert_{L^{2}} \lVert \dot{\Delta}_{j} u^{k} \rVert_{L^{2}} + \lVert \dot{\Delta}_{j'} \partial_{l} u^{k} \rVert_{L^{2}} \lVert \dot{\Delta}_{j} b^{k} \rVert_{L^{2}})\nonumber\\
\lesssim& \sum_{l'=1}^{3} \sum_{j} \lVert \partial_{l'} \dot{S}_{j-1} b^{l} \rVert_{L^{\infty}} (\lVert \dot{\Delta}_{j} \partial_{l} b^{k} \rVert_{L^{2}} \lVert \dot{\Delta}_{j} u^{k} \rVert_{L^{2}} + \lVert \dot{\Delta}_{j} \partial_{l} u^{k} \rVert_{L^{2}} \lVert \dot{\Delta}_{j} b^{k} \rVert_{L^{2}})\nonumber\\
\lesssim& \sum_{l'=1}^{3} \sum_{j}2^{j(2-\frac{2}{p_{l,l'} })}\sum_{j' \leq j-2} 2^{(j'-j)(2-\frac{2}{p_{l,l'} })}2^{j'(-2 + \frac{2}{p_{l,l'} })} \lVert \dot{\Delta}_{j'} \partial_{l'} b^{l} \rVert_{L^{\infty}}\nonumber\\
& \hspace{30mm} \times (\lVert \dot{\Delta}_{j} \partial_{l} b^{k} \rVert_{L^{2}} \lVert \dot{\Delta}_{j} u^{k} \rVert_{L^{2}} + \lVert \dot{\Delta}_{j} \partial_{l} u^{k} \rVert_{L^{2}} \lVert \dot{\Delta}_{j} b^{k} \rVert_{L^{2}})\nonumber\\
\lesssim& \sum_{l'=1}^{3} \lVert \partial_{l'} b^{l} \rVert_{\mathcal{B}_{p_{l,l'} }}\sum_{j}(2^{-\frac{j}{2}}\lVert \dot{\Delta}_{j} \partial_{l} b^{k} \rVert_{L^{2}})^{\frac{1}{p_{l,l'} }}(2^{\frac{j}{2}}\lVert \dot{\Delta}_{j} u^{k} \rVert_{L^{2}})^{\frac{1}{p_{l,l'} }}\nonumber\\
&\hspace{35mm} \times
(2^{\frac{j}{2}}\lVert \dot{\Delta}_{j} \partial_{l} b^{k} \rVert_{L^{2}})^{1-\frac{1}{p_{l,l'} }}(2^{\frac{3j}{2}}\lVert \dot{\Delta}_{j} u^{k} \rVert_{L^{2}})^{1-\frac{1}{p_{l,l'} }}\nonumber\\
&\hspace{30mm} + (2^{-\frac{j}{2}}\lVert \dot{\Delta}_{j} \partial_{l} u^{k} \rVert_{L^{2}})^{\frac{1}{p_{l,l'} }}(2^{\frac{j}{2}}\lVert \dot{\Delta}_{j} b^{k} \rVert_{L^{2}})^{\frac{1}{p_{l,l'} }}\nonumber\\
&\hspace{35mm} \times (2^{\frac{j}{2}}\lVert \dot{\Delta}_{j} \partial_{l} u^{k} \rVert_{L^{2}})^{1-\frac{1}{p_{l,l'} }}(2^{\frac{3j}{2}}\lVert \dot{\Delta}_{j} b^{k} \rVert_{L^{2}})^{1-\frac{1}{p_{l,l'} }}\nonumber\\
\lesssim& \sum_{l'=1}^{3} \lVert \partial_{l'} b^{l} \rVert_{\mathcal{B}_{p_{l,l'} }}(\lVert b\rVert_{\dot{H}^{\frac{1}{2}}}^{\frac{2}{p_{l,l'} }} + \lVert u\rVert_{\dot{H}^{\frac{1}{2}}}^{\frac{2}{p_{l,l'} }})(\lVert \nabla b\rVert_{\dot{H}^{\frac{1}{2}}}^{2(1-\frac{1}{p_{l,l'} })} + \lVert \nabla u\rVert_{\dot{H}^{\frac{1}{2}}}^{2(1-\frac{1}{p_{l,l'} })})\nonumber
\end{align}
by H$\ddot{o}$lder's inequality, commutator estimate used in (72) and Young's inequality for convolution. Next, we work similarly to $IV_{1,3}$ in (73): 
\begin{align}
&(V_{1} + VI_{1})_{3} + (V_{1} + VI_{1})_{6}\\
\lesssim& \sum_{j} 2^{j} \sum_{\lvert j-j'\rvert \leq 4, j'' \in [j-1, j'-1]}\lVert \sum_{l'=1}^{3} 2^{-j''} \tilde{\dot{\Delta}}_{j''}^{l'} \dot{\Delta}_{j''} \partial_{l'} b^{l} \rVert_{L^{\infty}}\nonumber\\
&\hspace{35mm} \times(\lVert \dot{\Delta}_{j} \dot{\Delta}_{j'} \partial_{l} b^{k} \rVert_{L^{2}} \lVert \dot{\Delta}_{j} u^{k} \rVert_{L^{2}} + \lVert \dot{\Delta}_{j} \dot{\Delta}_{j'} \partial_{l} u^{k} \rVert_{L^{2}} \lVert \dot{\Delta}_{j} b^{k} \rVert_{L^{2}} )\nonumber\\
\lesssim& \sum_{l'=1}^{3} \lVert \partial_{l'} b^{l} \rVert_{\mathcal{B}_{p_{l,l'} }}\sum_{j} (2^{-\frac{j}{2}}\lVert \dot{\Delta}_{j} \partial_{l} b^{k} \rVert_{L^{2}})^{\frac{1}{p_{l,l'} }}(2^{\frac{j}{2}}\lVert \dot{\Delta}_{j} u^{k} \rVert_{L^{2}})^{\frac{1}{p_{l,l'} }}\nonumber\\
&\hspace{35mm} \times (2^{\frac{j}{2}}\lVert \dot{\Delta}_{j} \partial_{l} b^{k} \rVert_{L^{2}})^{1-\frac{1}{p_{l,l'} }}
(2^{\frac{3j}{2}}\lVert \dot{\Delta}_{j} u^{k} \rVert_{L^{2}})^{1-\frac{1}{p_{l,l'} }}\nonumber\\
&\hspace{30mm} + (2^{-\frac{j}{2}}\lVert \dot{\Delta}_{j} \partial_{l} u^{k} \rVert_{L^{2}})^{\frac{1}{p_{l,l'} }}(2^{\frac{j}{2}}\lVert \dot{\Delta}_{j} b^{k} \rVert_{L^{2}})^{\frac{1}{p_{l,l'} }}\nonumber\\
&\hspace{35mm} \times (2^{\frac{j}{2}}\lVert \dot{\Delta}_{j} \partial_{l} u^{k} \rVert_{L^{2}})^{1-\frac{1}{p_{l,l'} }}
(2^{\frac{3j}{2}}\lVert \dot{\Delta}_{j} b^{k} \rVert_{L^{2}})^{1-\frac{1}{p_{l,l'} }}\nonumber\\
\lesssim& \sum_{l'=1}^{3} \lVert \partial_{l'} b^{l} \rVert_{\mathcal{B}_{p_{l,l'} }}(\lVert b\rVert_{\dot{H}^{\frac{1}{2}}}^{\frac{2}{p_{l,l'} }} + \lVert u\rVert_{\dot{H}^{\frac{1}{2}}}^{\frac{2}{p_{l,l'} }})(\lVert \nabla b\rVert_{\dot{H}^{\frac{1}{2}}}^{2(1-\frac{1}{p_{l,l'} })} + \lVert \nabla u\rVert_{\dot{H}^{\frac{1}{2}}}^{2(1-\frac{1}{p_{l,l'} })})\nonumber
\end{align}
by H$\ddot{o}$lder's inequality, that we can write $
\dot{\Delta}_{j''} b^{l} \approx \sum_{l'=1}^{3} 2^{-j''} \tilde{\dot{\Delta}}_{j''}^{l'} \dot{\Delta}_{j''} \partial_{l'} b^{l}$ and Young's inequality. Next, we work similarly to $IV_{2}$ in (75) to estimate 
\begin{align}
VI_{2} 
\lesssim& \sum_{j} \sum_{\lvert j-j'\rvert \leq 4} 2^{j} \lVert \dot{S}_{j'-1} \partial_{l} u^{k} \rVert_{L^{\infty}} \lVert \Delta_{j'} b^{l} \rVert_{L^{2}} \lVert \dot{\Delta}_{j} b^{k} \rVert_{L^{2}}\\
\lesssim& \sum_{j}2^{j(3-\frac{2}{p_{k,l}})} \sum_{j' \leq j-2}2^{(j'-j)(2-\frac{2}{p_{k,l}})}2^{j'(-2 + \frac{2}{p_{k,l}})}\lVert \dot{\Delta}_{j'} \partial_{l} u^{k} \rVert_{L^{\infty}} \lVert \dot{\Delta}_{j} b^{l} \rVert_{L^{2}} \lVert \dot{\Delta}_{j} b^{k} \rVert_{L^{2}}\nonumber\\
\lesssim& \lVert \partial_{l} u^{k} \rVert_{\mathcal{B}_{p_{k,l}}} \lVert b\rVert_{\dot{H}^{\frac{1}{2}}}^{\frac{2}{p_{k,l}}}\lVert \nabla b\rVert_{\dot{H}^{\frac{1}{2}}}^{2(1-\frac{1}{p_{k,l}})}\nonumber
\end{align}
by H$\ddot{o}$lder's inequality and Young's inequality for convolution. Finally, we use divergence-free condition and write $
\dot{\Delta}_{j'} b^{l} \approx \sum_{l'=1}^{3} \tilde{\dot{\Delta}}_{j'}^{l'} \dot{\Delta}_{j'} b^{l} \approx \sum_{l'=1}^{3} 2^{-j'} \tilde{\dot{\Delta}}_{j'}^{l'} \dot{\Delta}_{j'} \partial_{l'} b_{l}$ so that similarly to $IV_{3}$ in (76), we can estimate 
\begin{align}
VI_{3} \approx& \sum_{j} 2^{j} (\dot{\Delta}_{j} \sum_{j' \geq j-\delta} \partial_{l} (\sum_{l'=1}^{3} 2^{-j'} \tilde{\dot{\Delta}}_{j'}^{l'} \dot{\Delta}_{j'} \partial_{l'} b^{l} \tilde{\dot{\Delta}}_{j'} u^{k}) \lvert \dot{\Delta}_{j} b^{k})\\
\lesssim& \sum_{l'=1}^{3}\lVert \partial_{l'} b^{l} \rVert_{\mathcal{B}_{p_{l,l'} }}\sum_{j} \sum_{j' \geq j-\delta} 2^{(j-j')(\frac{1}{2} + \frac{1}{p_{l,l'} })}   (2^{\frac{j'}{2}}\lVert \tilde{\dot{\Delta}}_{j'} u^{k} \rVert_{L^{2}})^{\frac{1}{p_{l,l'} }}
\nonumber\\
&\hspace{20mm} \times (2^{\frac{3j'}{2}}\lVert \tilde{\dot{\Delta}}_{j'} u^{k} \rVert_{L^{2}})^{1-\frac{1}{p_{l,l'} }}(2^{\frac{j}{2}}\lVert \dot{\Delta}_{j} b^{k} \rVert_{L^{2}})^{\frac{1}{p_{l,l'} }}
(2^{\frac{3j}{2}}\lVert \dot{\Delta}_{j} b^{k} \rVert_{L^{2}})^{1-\frac{1}{p_{l,l'} }}\nonumber\\
\lesssim& \sum_{l'=1}^{3}\lVert \partial_{l'} b^{l} \rVert_{\mathcal{B}_{p_{l,l'} }} \lVert (2^{\frac{j}{2}}\lVert \tilde{\dot{\Delta}}_{j} u^{k} \rVert_{L^{2}})^{\frac{1}{p_{l,l'} }}
(2^{\frac{3j}{2}}\lVert \tilde{\dot{\Delta}}_{j} u^{k} \rVert_{L^{2}})^{1-\frac{1}{p_{l,l'} }} \rVert_{l^{2}}\nonumber\\
& \hspace{25mm} \times\lVert (2^{\frac{j}{2}}\lVert \dot{\Delta}_{j} b^{k} \rVert_{L^{2}})^{\frac{1}{p_{l,l'} }}
(2^{\frac{3j}{2}}\lVert \dot{\Delta}_{j} b^{k} \rVert_{L^{2}})^{1-\frac{1}{p_{l,l'} }} \rVert_{l^{2}}\nonumber\\
\lesssim& \sum_{l'=1}^{3} \lVert \partial_{l'} b^{l} \rVert_{\mathcal{B}_{p_{l,l'} }}(\lVert u\rVert_{\dot{H}^{\frac{1}{2}}}^{\frac{2}{p_{l,l'} }} + 
\lVert b\rVert_{\dot{H}^{\frac{1}{2}}}^{\frac{2}{p_{l,l'} }})(\lVert \nabla u\rVert_{\dot{H}^{\frac{1}{2}}}^{2(1-\frac{1}{p_{l,l'} })} + 
\lVert \nabla b\rVert_{\dot{H}^{\frac{1}{2}}}^{2(1-\frac{1}{p_{l,l'} })})\nonumber
\end{align}
by H$\ddot{o}$lder's, Bernstein's and Young's inequality of convolution. 

Hence, we obtain (68) due to (69), (77), (79), (80), (82)-(86). Applying (68) to (67), Young's and Gronwall's inequalities complete the proof of the Proposition 4.1. 
\end{proof}

\section{Proof of (4) in Theorem 1.1}

Firstly, for any $p \in (4, 6)$, 
\begin{align*}
\max_{1 \leq l \leq 3} (\lVert \partial_{l} u^{3} \rVert_{\mathcal{B}_{p}} + \lVert\partial_{l} b^{3} \rVert_{\mathcal{B}_{p}})\lesssim& \sup_{j\in\mathbb{Z}} 2^{j(\frac{1}{2} + \frac{2}{p} )}(\lVert \dot{\Delta}_{j} u^{3} \rVert_{L^{2}} + \lVert \dot{\Delta}_{j} b^{3} \rVert_{L^{2}})\\
\lesssim& \lVert u^{3} \rVert_{\dot{H}^{\frac{1}{2} + \frac{2}{p}}} + \lVert b^{3} \rVert_{\dot{H}^{\frac{1}{2} + \frac{2}{p}}}
\end{align*}
by Bernstein's inequality, which implies 
\begin{equation}
\max_{1 \leq l \leq 3} \int_{0}^{T^{\ast}} \lVert \partial_{l} u^{3} \rVert_{\mathcal{B}_{p}}^{p} + \lVert \partial_{l} b^{3} \rVert_{\mathcal{B}_{p}}^{p}d\tau \lesssim \int_{0}^{T^{\ast}} \lVert u^{3} \rVert_{\dot{H}^{\frac{1}{2} + \frac{2}{p}}}^{p} + \lVert b^{3} \rVert_{\dot{H}^{\frac{1}{2} + \frac{2}{p}}}^{p} d\tau \lesssim 1.
\end{equation}
Next, for $T < T^{\ast}$ by (5) for any $p \in (4, 6)$
\begin{align}
&\int_{0}^{T} \lVert \nabla_{h} u^{h} \rVert_{\mathcal{B}_{p}}^{p} + \lVert \nabla_{h} b^{h} \rVert_{\mathcal{B}_{p}}^{p} d\tau\\
\approx& \int_{0}^{T} \lVert \nabla_{h}\nabla_{h}^{\bot} \Delta_{h}^{-1} \omega \rVert_{\mathcal{B}_{p}}^{p} + \lVert \nabla_{h} \nabla_{h}^{\bot} \Delta_{h}^{-1} d \rVert_{\mathcal{B}_{p}}^{p}\nonumber\\
&+ \lVert \nabla_{h}\nabla_{h} \Delta_{h}^{-1} \partial_{3}u^{3} \rVert_{\mathcal{B}_{p}}^{p} + \lVert \nabla_{h} \nabla_{h} \Delta_{h}^{-1} \partial_{3}b^{3} \rVert_{\mathcal{B}_{p}}^{p}d\tau \nonumber
\end{align}
where 
\begin{align}
&\int_{0}^{T} \lVert \nabla_{h} \nabla_{h}^{\bot} \Delta_{h}^{-1} \omega \rVert_{\mathcal{B}_{p}}^{p} +  \lVert \nabla_{h} \nabla_{h}^{\bot} \Delta_{h}^{-1} d \rVert_{\mathcal{B}_{p}}^{p} d\tau\\
\lesssim& \int_{0}^{T} \lVert \omega \rVert_{L^{\frac{3}{2}}}^{p(1-\frac{3}{2p})}
\lVert \nabla \omega \rVert_{L^{\frac{9}{5}}}^{\frac{3}{2}}
+ \lVert d \rVert_{L^{\frac{3}{2}}}^{p(1-\frac{3}{2p})}
\lVert  \nabla d \rVert_{L^{\frac{9}{5}}}^{\frac{3}{2}} d\tau\nonumber\\
\lesssim& \int_{0}^{T} \lVert \omega \rVert_{L^{\frac{3}{2}}}^{p(1-\frac{3}{2p})}
\lVert \nabla \omega_{\frac{3}{4}} \rVert_{L^{2}}^{2}
+ \lVert d \rVert_{L^{\frac{3}{2}}}^{p(1-\frac{3}{2p})}
\lVert  \nabla d_{\frac{3}{4}} \rVert_{L^{\frac{9}{5}}}^{2} d\tau\nonumber\\
\lesssim& \sup_{\tau \in [0,T]} \lVert \omega(\tau) \rVert_{L^{\frac{3}{2}}}^{p(1- \frac{3}{2p})}\int_{0}^{T} \lVert \nabla \omega_{\frac{3}{4}} \rVert_{L^{2}}^{2}d\tau + 
\sup_{\tau \in [0,T]} \lVert d(\tau) \rVert_{L^{\frac{3}{2}}}^{p(1- \frac{3}{2p})}\int_{0}^{T} \lVert \nabla d_{\frac{3}{4}} \rVert_{L^{2}}^{2}d\tau\nonumber
\end{align}
by (11), continuity of Riesz transform in $L^{p}, p \in (1, \infty)$ and (9) while 
\begin{align}
&\int_{0}^{T} \lVert \nabla_{h}\nabla_{h} \Delta_{h}^{-1} \partial_{3}u^{3} \rVert_{\mathcal{B}_{p}}^{p} + \lVert \nabla_{h} \nabla_{h} \Delta_{h}^{-1} \partial_{3}b^{3} \rVert_{\mathcal{B}_{p}}^{p}d\tau \\
\lesssim& \int_{0}^{T}\lvert \sup_{j\in \mathbb{Z}} 2^{j(\frac{1}{2} + \frac{2}{p})} \lVert \dot{\Delta}_{j} u^{3} \rVert_{L^{2}} \rvert^{p} + \lvert \sup_{j\in \mathbb{Z}} 2^{j(\frac{1}{2} + \frac{2}{p})}\lVert \dot{\Delta}_{j} b^{3} \rVert_{L^{2}} \rvert^{p}d\tau\nonumber\\
\lesssim& \int_{0}^{T} \lVert u^{3} \rVert_{\dot{H}^{\frac{1}{2} + \frac{2}{p}}}^{p} +  \lVert b^{3} \rVert_{\dot{H}^{\frac{1}{2} + \frac{2}{p}}}^{p} d\tau\nonumber 
\end{align}
by Bernstein's inequality and continuity of Riesz transform in $L^{p}, p \in (1, \infty)$. Therefore, applying (89) and (90) in (88), by (59) we obtain 
\begin{equation*}
\int_{0}^{T^{\ast}} \lVert \nabla_{h} u^{h} \rVert_{\mathcal{B}_{p}}^{p} + \lVert \nabla_{h} b^{h} \rVert_{\mathcal{B}_{p}}^{p} d\tau \lesssim 1.
\end{equation*} 
Finally, by (5) for any $T < T^{\ast}$
\begin{align}
&\int_{0}^{T} \lVert \partial_{3}u^{h} \rVert_{\mathcal{B}_{2}}^{2} + \lVert \partial_{3} b^{h} \rVert_{\mathcal{B}_{2}}^{2} d\tau\\
\lesssim& \int_{0}^{T}  \lVert \partial_{3} \nabla_{h}^{\bot} \Delta_{h}^{-1} \omega \rVert_{\mathcal{B}_{2}}^{2} +  \lVert \partial_{3} \nabla_{h}^{\bot} \Delta_{h}^{-1} d \rVert_{\mathcal{B}_{2}}^{2} + \lVert \partial_{3} \nabla_{h}\Delta_{h}^{-1} \partial_{3}u^{3} \rVert_{\mathcal{B}_{2}}^{2} + \lVert \partial_{3} \nabla_{h} \Delta_{h}^{-1} \partial_{3}b^{3} \rVert_{\mathcal{B}_{2}}^{2} d\tau\nonumber
\end{align}
where 
\begin{align}
&\int_{0}^{T}  \lVert \partial_{3} \nabla_{h}^{\bot} \Delta_{h}^{-1} \omega \rVert_{\mathcal{B}_{2}}^{2} +  \lVert \partial_{3} \nabla_{h}^{\bot} \Delta_{h}^{-1} d \rVert_{\mathcal{B}_{2}}^{2} d\tau\\
\lesssim& \int_{0}^{T} \lvert \sup_{j\in \mathbb{Z}} 2^{-j} \sum_{k \leq j + \delta, l \leq j+ \delta}\lVert \dot{\Delta}_{j} \dot{\Delta}_{k}^{h} \dot{\Delta}_{l}^{v} \partial_{3} \nabla_{h}^{\bot} \Delta_{h}^{-1} \omega \rVert_{L^{\infty}} \rvert^{2}\nonumber\\
& \hspace{5mm} + \lvert \sup_{j\in \mathbb{Z}} 2^{-j} \sum_{k \leq j + \delta, l \leq j+ \delta}\lVert \dot{\Delta}_{j} \dot{\Delta}_{k}^{h} \dot{\Delta}_{l}^{v} \partial_{3} \nabla_{h}^{\bot} \Delta_{h}^{-1} d \rVert_{L^{\infty}} \rvert^{2}d\tau\nonumber\\
\lesssim&  \int_{0}^{T} \lvert \sup_{j\in \mathbb{Z}} 2^{-j} \sum_{k \leq j+ \delta, l \leq j + \delta} 2^{l(\frac{2}{3})} 2^{\frac{k}{3}} \lVert \dot{\Delta}_{j}\dot{\Delta}_{k}^{h} \dot{\Delta}_{l}^{v} \partial_{3} \omega \rVert_{L^{\frac{3}{2}}} \rvert^{2} \nonumber\\
&\hspace{5mm} +\lvert \sup_{j\in \mathbb{Z}} 2^{-j} \sum_{k \leq j+ \delta, l \leq j + \delta} 2^{l(\frac{2}{3})} 2^{\frac{k}{3}}  \lVert \dot{\Delta}_{j}\dot{\Delta}_{k}^{h} \dot{\Delta}_{l}^{v} \partial_{3} d \rVert_{L^{\frac{3}{2}}} \rvert^{2}d\tau\nonumber \\
\lesssim& \int_{0}^{T} \lVert \partial_{3} \omega \rVert_{L^{\frac{3}{2}}}^{2} + \lVert \partial_{3} d\rVert_{L^{\frac{3}{2}}}^{2} d\tau\nonumber\\
\lesssim& \sup_{\tau \in [0,T]} \lVert \omega_{\frac{3}{4}}(\tau) \rVert_{L^{2}}^{\frac{2}{3}} \int_{0}^{T} \lVert \nabla \omega_{\frac{3}{4}} \rVert_{L^{2}}^{2} d\tau + \sup_{\tau \in [0,T]} \lVert d_{\frac{3}{4}}(\tau) \rVert_{L^{2}}^{\frac{2}{3}} \int_{0}^{T} \lVert \nabla d_{\frac{3}{4}}\rVert_{L^{2}}^{2} d\tau\nonumber
\end{align}
by Bernstein's inequality, and $L^{\frac{3}{2}}(\mathbb{R}^{3}) \subset \dot{B}_{\frac{3}{2}, 2}^{0}$ (cf. [5]) and (9) whereas for fixed $\theta \in (\frac{1}{2} - \frac{2}{p}, \frac{1}{6})$    
\begin{align}
&\int_{0}^{T} \lVert \partial_{3} \nabla_{h}\Delta_{h}^{-1} \partial_{3}u^{3} \rVert_{\mathcal{B}_{2}}^{2} + \lVert \partial_{3} \nabla_{h} \Delta_{h}^{-1} \partial_{3}b^{3} \rVert_{\mathcal{B}_{2}}^{2} d\tau\\
\lesssim& \int_{0}^{T} \lvert \sup_{j \in \mathbb{Z}} 2^{-j} \sum_{k \leq j + \delta, l \leq j + \delta} \lVert \dot{\Delta}_{j}\dot{\Delta}_{k}^{h} \dot{\Delta}_{l}^{v} \partial_{3}^{2} \nabla_{h} \Delta_{h}^{-1} u^{3} \rVert_{L^{\infty}} \rvert^{2}\nonumber\\
&+ \lvert \sup_{j \in \mathbb{Z}} 2^{-j} \sum_{k \leq j + \delta, l \leq j + \delta} \lVert \dot{\Delta}_{j}\dot{\Delta}_{k}^{h} \dot{\Delta}_{l}^{v} \partial_{3}^{2} \nabla_{h} \Delta_{h}^{-1} b^{3} \rVert_{L^{\infty}} \rvert^{2}d\tau\nonumber \\
\lesssim& \int_{0}^{T} \lvert \sup_{j \in \mathbb{Z}} 2^{-j} \sum_{k \leq j + \delta, l \leq j + \delta} 2^{k(\frac{1}{2} - \theta)}2^{l(\frac{1}{2} + \theta)}2^{k(\theta - \frac{1}{2})}2^{-l\theta} \lVert \dot{\Delta}_{j} \dot{\Delta}_{k}^{h} \dot{\Delta}_{l}^{v} \partial_{3}^{2} u^{3} \rVert_{L^{2}} \rvert^{2}\nonumber\\
&+ \lvert \sup_{j \in \mathbb{Z}} 2^{-j} \sum_{k \leq j + \delta, l \leq j + \delta} 2^{k(\frac{1}{2} - \theta)}2^{l(\frac{1}{2} + \theta)}2^{k(\theta - \frac{1}{2})}2^{-l\theta} \lVert \dot{\Delta}_{j} \dot{\Delta}_{k}^{h} \dot{\Delta}_{l}^{v} \partial_{3}^{2} b^{3} \rVert_{L^{2}} \rvert^{2}d\tau\nonumber\\
\lesssim& \int_{0}^{T} \lVert \partial_{3}^{2} u^{3} \rVert_{\mathcal{H}_{\theta}}^{2} + \lVert \partial_{3}^{2} b^{3} \rVert_{\mathcal{H}_{\theta}}^{2}d\tau\nonumber
\end{align}
by Bernstein's inequality. Thus, applying (92) and (93) in (91), by (59) and (60), 
\begin{equation*}
\int_{0}^{T^{\ast}} \lVert \partial_{3}u^{h} \rVert_{\mathcal{B}_{2}}^{2} + \lVert \partial_{3} b^{h} \rVert_{\mathcal{B}_{2}}^{2} d\tau \lesssim 1.
\end{equation*}
Due to Proposition 4.1, this completes the proof of (4). 

\section{Appendix}

\subsection{Local theory of Theorem 1.1}

We let $X^{\pm} \triangleq u \pm b, Y^{\pm} \triangleq \Omega \pm j$ so that from (18a), (18b) and (19)
\begin{equation*}
\partial_{t} Y^{\pm} - \Delta Y^{\pm} + (X^{\mp} \cdot\nabla)Y^{\pm} = (Y^{\pm} \cdot\nabla) X^{\mp} \pm 2M(u,b).
\end{equation*}
By Sobolev embedding of $\dot{W}^{\frac{1}{2}, \frac{3}{2}}(\mathbb{R}^{3}) \hookrightarrow L^{2}(\mathbb{R}^{3})$, and continuity of Riesz transform in $L^{\frac{3}{2}}(\mathbb{R}^{3})$, we have $u_{0}, b_{0} \in \dot{B}_{2,2}^{\frac{1}{2}}(\mathbb{R}^{3})$. By [23] (also [35]), we find a unique solution pair $u, b \in C(0, T^{\ast}; \dot{H}^{\frac{1}{2}}(\mathbb{R}^{3})) \cap L_{loc}^{2} (0, T^{\ast}; \dot{H}^{\frac{3}{2}}(\mathbb{R}^{3}))$. Now by Lemma 3.1 [12] ,  
\begin{align*}
& \frac{2}{3} \lVert Y^{+} \rVert_{L^{\frac{3}{2}}}^{\frac{3}{2}} + \frac{1}{2} \int_{0}^{t} \lvert \nabla Y^{+} \rvert^{2} \lvert Y^{+} \rvert^{-\frac{1}{2}}d\tau\\
=& \frac{2}{3} \lVert Y_{0}^{+}  \rVert_{L^{\frac{3}{2}}}^{\frac{3}{2}} + \int_{0}^{t} \int [(Y^{+} \cdot\nabla) X^{-} + 2M(u,b) ] \lvert Y ^{+} \rvert^{-\frac{1}{2}} Y^{+}  d\tau \\
\leq& \frac{2}{3} \lVert Y_{0}^{+} \rVert_{L^{\frac{3}{2}}}^{\frac{3}{2}} + \int_{0}^{t} \lVert Y^{+} \rVert_{L^{3}} \lVert \nabla X^{-} \rVert_{L^{3}} \lVert Y^{+} \rVert_{L^{\frac{3}{2}}} ^{\frac{1}{2}} + \lVert \nabla u \rVert_{L^{3}} \lVert \nabla b \rVert_{L^{3}} \lVert Y^{+} \rVert_{L^{\frac{3}{2}}} ^{\frac{1}{2}}d\tau\\
\lesssim&  \lVert Y_{0}^{+} \rVert_{L^{\frac{3}{2}}}^{\frac{3}{2}} + \int_{0}^{t} \lVert u\rVert_{\dot{H}^{\frac{3}{2}}}^{2} + \lVert b\rVert_{\dot{H}^{\frac{3}{2}}}^{2} d\tau + \int_{0}^{t} (\lVert u\rVert_{\dot{H}^{\frac{3}{2}}}^{2} + \lVert b\rVert_{\dot{H}^{\frac{3}{2}}}^{2}) \lVert Y^{+} \rVert_{L^{\frac{3}{2}}}^{\frac{3}{2}}d\tau 
\end{align*}
by H$\ddot{o}$lder's inequalities, Sobolev embedding of $\dot{H}^{\frac{1}{2}}(\mathbb{R}^{3}) \hookrightarrow L^{3}(\mathbb{R}^{3})$ and Young's inequality so that for $T < T^{\ast}$
\begin{equation*}
\sup_{t\in [0,T]} \lVert Y^{+}(t) \rVert_{L^{\frac{3}{2}}}^{\frac{3}{2}} + \int_{0}^{T} 
\lvert \nabla Y^{+} \rvert^{2} \lvert Y^{+} \rvert^{-\frac{1}{2}}
 d\tau \lesssim 
(e+ \lVert Y_{0}^{+} \rVert_{L^{\frac{3}{2}}}^{\frac{3}{2}}) e^{\int_{0}^{T} \lVert u\rVert_{\dot{H}^{\frac{3}{2}}}^{2} + \lVert b\rVert_{\dot{H}^{\frac{3}{2}}}^{2}d\tau }.
\end{equation*}
Similar procedure on the equation of $Y^{-}$ gives in sum
\begin{align*}
\sup_{t\in [0, T]}(\lVert \Omega\rVert_{L^{\frac{3}{2}}}^{\frac{3}{2}} + \lVert j&\rVert_{L^{\frac{3}{2}}}^{\frac{3}{2}})(t) + \int_{0}^{T} \lvert \nabla (\Omega + j)\rvert^{2} \lvert \Omega + j\rvert^{-\frac{1}{2}} + \lvert \nabla (\Omega - j)\rvert^{2} \lvert \Omega - j\rvert^{-\frac{1}{2}} d\tau\\
\lesssim& (1+ \lVert \Omega_{0} \rVert_{L^{\frac{3}{2}}}^{\frac{3}{2}} + \lVert j_{0} \rVert_{L^{\frac{3}{2}}}^{\frac{3}{2}})\exp\left(\int_{0}^{T} \lVert u\rVert_{\dot{H}^{\frac{3}{2}}}^{2} + \lVert b\rVert_{\dot{H}^{\frac{3}{2}}}^{2} d\tau\right) \lesssim 1.
\end{align*}
This completes the proof of the local theory of Theorem 1.1. 

\subsection{Additional estimates}

Here we prove two additional estimates: 
\begin{align}
\lVert fg\rVert_{\dot{H}^{s-1, \sigma_{1} + \sigma_{2} - \frac{1}{2}}}\lesssim \lVert f\rVert_{(\dot{B}_{2,1}^{1})_{h} (\dot{B}_{2,1}^{\sigma_{1}})_{v}}\lVert g\rVert_{\dot{H}^{s-1, \sigma_{2}}}, 
\end{align}
where $\sigma_{1} < \frac{1}{2}, \sigma_{2} < \frac{1}{2}, \sigma_{1} + \sigma_{2} > 0, 2 > s > 0$ and 
\begin{align}
\lVert fg\rVert_{\dot{H}^{s_{1} + s_{2} - 1, \sigma - \frac{1}{2}}} \lesssim \lVert f\rVert_{(\dot{B}_{2,2}^{s_{1}})_{h} (\dot{B}_{2,1}^{\frac{1}{2}})_{v}} \lVert g\rVert_{\dot{H}^{s_{2}, \sigma - \frac{1}{2}}}
\end{align}
where $s_{1} < 1, s_{2} < 1, s_{1} + s_{2} > 0, 1 > \sigma > 0$. 
Since these are standard applications, we only sketch (94); the proof of (95) is similar. 

Due to the following horizontal and vertical Bony paraproduct decompositions 
\begin{align*}
&T^{h} (f,g) \triangleq \sum_{m} \dot{S}_{m-1}^{h} f \dot{\Delta}_{m}^{h} g , \hspace{2mm} R^{h}(f,g) \triangleq \sum_{m} \dot{\Delta}_{m}^{h} f \tilde{\dot{\Delta}}_{m}^{h}g, \hspace{2mm} \tilde{T^{h}} (f,g) \triangleq T^{h}(g,f),\\
&T^{v} (f,g) \triangleq \sum_{n} \dot{S}_{n-1}^{v} f \dot{\Delta}_{n}^{v} g , \hspace{3mm} R^{v}(f,g) \triangleq \sum_{n} \dot{\Delta}_{n}^{v} f \tilde{\dot{\Delta}}_{n}^{v}g, \hspace{4mm} \tilde{T^{v}} (f,g) \triangleq T^{v}(g,f),
\end{align*}
we can write $fg = (T^{h} + R^{h} + \tilde{T^{h}}) (T^{v} + R^{v} + \tilde{T^{v}}) (f,g)$ in nine parts: e.g. 
\begin{align*}
R^{h} T^{v} (f,g) = \sum_{m,n} \dot{\Delta}_{m}^{h} \dot{S}_{n-1}^{v} f \tilde{\dot{\Delta}}_{m}^{h} \dot{\Delta}_{n}^{v} g.
\end{align*}
Let us estimate this term:   
\begin{align*}
& \lVert \dot{\Delta}_{k}^{h} \dot{\Delta}_{l}^{v}\sum_{k \leq m + \delta, \lvert l-n\rvert \leq 4} \dot{\Delta}_{m}^{h} \dot{S}_{n-1}^{v} f \tilde{\dot{\Delta}}_{m}^{h} \dot{\Delta}_{n}^{v} g \rVert_{L^{2}}\\
\lesssim& 2^{k} \sum_{k \leq m + \delta} \lVert \dot{\Delta}_{m}^{h} \dot{S}_{l-1}^{v} f \tilde{\dot{\Delta}}_{m}^{h} \dot{\Delta}_{l}^{v} g\rVert_{L_{h}^{1}L_{v}^{2}}\\
\lesssim& 2^{k} \sum_{k \leq m + \delta} \sum_{l' \leq l-2}\lVert \dot{\Delta}_{m}^{h} \dot{\Delta}_{l'}^{v} f\rVert_{L_{h}^{2} L_{v}^{\infty}}\lVert \tilde{\dot{\Delta}}_{m}^{h} \dot{\Delta}_{l}^{v} g\rVert_{L^{2}}\\
\lesssim& 2^{k} \sum_{k \leq m + \delta} \sum_{l' \leq l-2} 2^{\frac{l'}{2}}\lVert \dot{\Delta}_{m}^{h} \dot{\Delta}_{l'}^{v} f\rVert_{L^{2}} \lVert \tilde{\dot{\Delta}}_{m}^{h} \dot{\Delta}_{l}^{v} g\rVert_{L^{2}}
\end{align*}
by Bernstein's and H$\ddot{o}$lder's inequalities. Thus, 
\begin{align*}
& 2^{k(s-1 )}2^{l(\sigma_{1} + \sigma_{2} - \frac{1}{2})}\lVert \dot{\Delta}_{k}^{h} \dot{\Delta}_{l}^{v} R^{h} T^{v}(f,g) \rVert_{L^{2}}\\
\lesssim& \sum_{\alpha = -1}^{1}\sum_{k \leq m + \delta} \sum_{l' \leq l-2} 2^{(k-m)s}2^{(l' -l) (\frac{1}{2} - \sigma_{1})}2^{m} 2^{l'\sigma_{1}}\lVert \dot{\Delta}_{m}^{h} \dot{\Delta}_{l'}^{v} f\rVert_{L^{2}}\\
&\hspace{30mm} \times 2^{\alpha (s-1)}2^{(m-\alpha)(s-1)}2^{l\sigma_{2}}\lVert \dot{\Delta}_{m-\alpha}^{h} \dot{\Delta}_{l}^{v} g\rVert_{L^{2}}.
\end{align*}
We take $l^{2}$-norm in $k$ now to obtain 
\begin{align*}
&\left\lVert \left(2^{k(s-1 )}2^{l(\sigma_{1} + \sigma_{2} - \frac{1}{2})}\lVert \dot{\Delta}_{k}^{h} \dot{\Delta}_{l}^{v} R^{h} T^{v}(f,g) \rVert_{L^{2}}\right)_{k} \right\rVert_{l^{2}}\\
\lesssim& \sum_{\alpha = -1}^{1}\left\lVert 
\left(\sum_{l' \leq l-2} 2^{(l' -l) (\frac{1}{2} - \sigma_{1})}2^{k} 2^{l'\sigma_{1}}\lVert \dot{\Delta}_{k}^{h} \dot{\Delta}_{l'}^{v} f\rVert_{L^{2}} 2^{(k-\alpha)(s-1)}2^{l\sigma_{2}}\lVert \dot{\Delta}_{k-\alpha}^{h} \dot{\Delta}_{l}^{v} g\rVert_{L^{2}}\right)_{k} 
\right\rVert_{l^{1}}
\end{align*}
by Young's inequality for convolution. We now take $l^{2}$-norm in $l$ and use Minkowski's inequality to obtain 
\begin{align*}
&\left\lVert \left(2^{k(s-1 )}2^{l(\sigma_{1} + \sigma_{2} - \frac{1}{2})}\lVert \dot{\Delta}_{k}^{h} \dot{\Delta}_{l}^{v} R^{h} T^{v}(f,g) \rVert_{L^{2}}\right)_{k,l} \right\rVert_{l^{2}}\\
\lesssim& \sum_{\alpha = -1}^{1}\left \lVert \left(
\left\lVert 
\left(\sum_{l' \leq l-2} 2^{(l' -l) (\frac{1}{2} - \sigma_{1})}2^{k} 2^{l'\sigma_{1}}\lVert \dot{\Delta}_{k}^{h} \dot{\Delta}_{l'}^{v} f\rVert_{L^{2}} 2^{(k-\alpha)(s-1)}2^{l\sigma_{2}}\lVert \dot{\Delta}_{k-\alpha}^{h} \dot{\Delta}_{l}^{v} g\rVert_{L^{2}}\right)_{l} 
\right\rVert_{l^{2}}\right)_{k} \right\rVert_{l^{1}}\\
\lesssim& \sum_{\alpha = -1}^{1}\left\lVert \left(
\left\lVert 
\left( 2^{k} 2^{l\sigma_{1}}\lVert \dot{\Delta}_{k}^{h} \dot{\Delta}_{l}^{v} f\rVert_{L^{2}}\right)_{l} \right\rVert_{l^{2}}  \left\lVert \left( 2^{(k-\alpha)(s-1)}2^{l\sigma_{2}}\lVert \dot{\Delta}_{k-\alpha}^{h} \dot{\Delta}_{l}^{v} g\rVert_{L^{2}}\right)_{l} 
\right\rVert_{l^{2}}\right)_{k} \right\rVert_{l^{1}}\\
\lesssim& \lVert f\rVert_{(\dot{B}_{2,1}^{1})_{h} (\dot{B}_{2,1}^{\sigma_{1}})_{v}}\lVert g\rVert_{\dot{H}^{s-1, \sigma_{2}}}
\end{align*}
by H$\ddot{o}$lder's inequality and Young's inequality for convolution. 

The other terms are similar and we refer to e.g. [19, 29] for details.

\end{document}